\documentclass[leqno,12pt,twoside]{article}

\usepackage{amssymb}
\usepackage{euscript}
\usepackage[dvips]{graphicx}
\usepackage{flafter}
\usepackage{pstricks}
\usepackage[latin1]{inputenc}
\usepackage{amsmath}
\usepackage{amsfonts}
\usepackage{fancyhdr}

\usepackage{verbatim}

\usepackage{mathrsfs} 

\evensidemargin\oddsidemargin

\begin{document}

\baselineskip=18pt \setcounter{page}{1}

\renewcommand{\theequation}{\thesection.\arabic{equation}}
\newtheorem{theorem}{Theorem}[section]
\newtheorem{lemma}[theorem]{Lemma}
\newtheorem{proposition}[theorem]{Proposition}
\newtheorem{corollary}[theorem]{Corollary}
\newtheorem{remark}[theorem]{Remark}
\newtheorem{fact}[theorem]{Fact}
\newtheorem{problem}[theorem]{Problem}
\newtheorem{example}[theorem]{Example}
\newtheorem{question}[theorem]{Question}
\newtheorem{conjecture}[theorem]{Conjecture}
\newtheorem{definition}[theorem]{Definition}
\newcommand{\eqnsection}{
\renewcommand{\theequation}{\thesection.\arabic{equation}}
    \makeatletter
    \csname  @addtoreset\endcsname{equation}{section}
    \makeatother}
\eqnsection

\def\r{{\mathbb R}}
\def\e{{\mathbb E}}
\def\p{{\mathbb P}}
\def\P{{\bf P}}
\def\E{{\bf E}}
\def\Q{{\hat {\bf P}}}
\def\z{{\mathbb Z}}
\def\N{{\mathbb N}}
\def\T{{\mathbb T}}
\def\G{{\mathbb G}}
\def\L{{\mathbb L}}

\def\deg{\chi}

\def\ee{\mathrm{e}}
\def\d{\, \mathrm{d}}
\def\S{\mathscr{S}}

\def\I{\Omega}



\pagestyle{fancy}

\fancyhead{}
\fancyhead[CE]{\small MINIMUM OF A BRANCHING RANDOM WALK}
\fancyhead[CO]{\small E. AIDEKON}
\renewcommand{\headrulewidth}{0pt}

\thispagestyle{plain}

\centerline{\Large \bf Convergence in law of the minimum of a }

\bigskip

\centerline{\Large \bf branching random walk}

\bigskip
\medskip

\centerline{Elie A\"\i d\'ekon \footnote{ Supported in part by the Netherlands Organisation for scientific Research (NWO).
}}

\bigskip

\centerline{\it Eindhoven University of Technology}

\bigskip
\bigskip
\bigskip
\bigskip
\bigskip
\bigskip

{\leftskip=2truecm \rightskip=2truecm \baselineskip=15pt \small

\noindent{\slshape\bfseries Summary.} We consider the minimum of a super-critical branching random walk. In \cite{addario-berry-reed}, Addario-Berry and Reed proved the tightness of the minimum centered around its mean value. We show that a convergence in law holds, giving the analog of a well-known result of Bramson \cite{bramson83} in the case of the branching Brownian motion.

} 

\vglue50pt
\section{Introduction}

We consider a branching random walk defined as follows. The process starts with one particle located at $0$. At time $1$, the particle dies and gives birth to a point process $\mathcal L$. Then, at each time $n\in\mathbb{N}$, the particles of generation $n$ die and give birth to independent copies of the point process $\mathcal L$, translated to their position. If $\T$ is the genealogical tree of the process, we see that $\T$ is a Galton-Watson tree, and we denote by $|x|$ the generation of the vertex $x\in\T$ (the ancestor is the only particle at generation $0$). For each $x\in\T$, we denote by $V(x)\in\r$ its position on the real line. With this notation, $(V(x),|x|=1)$ is distributed as $\mathcal L$. The collection of positions $(V(x),\,x\in\T)$ defines our branching random walk. \\

We assume that we are in the boundary case (in the sense of \cite{biggins-kyprianou05})
\begin{equation} \label{cond-bound}
\E\left[\sum_{|x|=1} 1\right]>1,\qquad \E\left[\sum_{|x|=1}\ee^{-V(x)}\right]=1,\qquad \E\left[\sum_{|x|=1}V(x)\ee^{-V(x)}\right]=0.
\end{equation}

\noindent Every branching random walk satisfying mild assumptions can be reduced to this case by some renormalization. Notice that we may have $\sum_{|x|=1} 1 =\infty$ with positive probability. We are interested in the minimum at time $n$
$$
M_n := \min\{\, V(x),\, |\,x|=n \}
$$

\noindent where $\min \emptyset:=\infty$. Writing for $y\in\r\cup\{\pm \infty\}$, $y_+:=\max(y,0)$, we introduce the random variables
\begin{equation}\label{def:X}
X:=\sum_{|x|=1} \ee^{-V(x)},  \qquad \tilde X:=\sum_{|x|=1} V(x)_+\ee^{-V(x)}.
\end{equation}

\noindent We assume throughout the remainder of the paper, including in the statements of Theorems, Lemmas etc. that
\begin{itemize}
\item the distribution of $\mathcal L$ is non-lattice,
\item we have
\begin{eqnarray}\label{cond-2}
\E\left[\sum_{|x|=1}V(x)^2\ee^{-V(x)}\right] < \infty, \\
\label{cond-ln} \E\left[ X (\ln_+ X)^2\right]<\infty,\qquad \E\left[ \tilde X \ln_+ \tilde X \right]<\infty.
\end{eqnarray}
\end{itemize}

These assumptions are discussed after Theorem \ref{main}. Under (\ref{cond-bound}), the minimum $M_n$ goes to infinity, as it can be easily seen from the fact that $\sum_{|u|=n}\ee^{-V(u)}$ goes to zero (\cite{lyons}). The law of large numbers for the speed of the minimum goes back to the works of Hammersley \cite{hammersley}, Kingman \cite{kingman} and Biggins \cite{biggins76}, and  we know that  ${M_n \over n}$ converges almost surely to $0$ in the boundary case. The second order was recently found separately by Hu and Shi \cite{yzpolymer}, and Addario-Berry and Reed \cite{addario-berry-reed}, and is proved to be equal to ${3\over 2} \ln n$ in probability, though there exist almost sure fluctuations (Theorem 1.2 in \cite{yzpolymer}). In \cite{addario-berry-reed}, the authors computed the expectation of $M_n$ to within $O(1)$, and  showed, under suitable assumptions, that the sequence of the minimum is tight around its mean. Through recursive equations, Bramson and Zeitouni \cite{bramson-zeitouni} obtained the tightness of $M_n$ around its median, when assuming some properties on the decay of the tail distribution. In the particular case where the step distribution is log-concave, the convergence in law of $M_n$ around its median was proved earlier by Bachmann \cite{bachmann}. The aim of this paper is to get the convergence of the minimum $M_n$ centered around $\frac32 \ln n$ for a general class of branching random walks. This is the analog of the seminal work from Bramson \cite{bramson83}, to which our approach bears some resemblance. To state our result, we introduce the {\it derivative martingale}, defined for any $n\ge 0$ by
\begin{equation}\label{def:dWn}
D_n := \sum_{|x|=n}V(x)\ee^{-V(x)}.
\end{equation}

\noindent From \cite{biggins-kyprianou04} (and Proposition \ref{l:Dbeta} in the Appendix), we know that the martingale converges almost surely to some limit $D_{\infty}$, which is strictly positive on the set of non-extinction of $\T$. Notice that under (\ref{cond-bound}), the tree $\T$ has a positive probability to survive.

\bigskip

\begin{theorem}\label{main}
There exists a constant $C^*\in (0,\infty)$ such that for any real $x$,
\begin{equation}
\lim_{n\to\infty}\P\left(M_n \ge {3\over 2} \ln n +x \right) = \E\left[\ee^{-C^*\ee^{x}D_\infty}\right].
\end{equation}
\end{theorem}

\medskip
\noindent {\bf Remark 1}.  We can see our theorem as the analog of the result of Lalley and Sellke \cite{lalley-sellke} in the case of the branching Brownian motion : the minimum converges to a random shift of the Gumbel distribution.\\

\noindent {\bf Remark 2}.  The condition of non-lattice distribution is necessary since it is hopeless to have a convergence in law around $\frac32 \ln n$ in general. We do not know if an analogous result holds in the lattice case. If (\ref{cond-2}) does not hold, we can expect, under suitable conditions, to have still a convergence in law but centered around $\kappa \ln n$ for some constant $\kappa  \neq 3/2$. This comes from the different behaviour of the probability to remain positive for one-dimensional random walks with infinite variance. Finally, the condition (\ref{cond-ln}) appears naturally for $D_\infty$ not being identically zero (see \cite{biggins-kyprianou04}, Theorem 5.2). \\

\bigskip

The proof of the theorem is divided into three steps. First, we look at the tail distribution of the minimum $M_n^{\rm kill}$ of the branching random walk killed below zero, i.e $M_n^{\rm kill}:= \min\{V(x),\, V(x_k)\ge 0,\, \forall 0\le k\le |x|\}$, where $x_k$ denotes the ancestor of $x$ at generation $k$.
\begin{proposition}\label{p:tailminkillmain}
There exists a constant $C_1>0$ such that
$$
\limsup_{z\to\infty}\limsup_{n\to\infty}\Big| \ee^z\P\left(M_n^{\rm kill} < {3\over 2} \ln n-z\right) -C_1\Big|=0.
$$
\end{proposition}

\bigskip

\noindent This allows us to get the tail distribution of $M_n$ in a second stage.
\begin{proposition}\label{p:tailminmain}
We have
$$
\limsup_{z\to\infty}\limsup_{n\to\infty}\Big| {\ee^z\over z}\P\left(M_n < {3\over 2} \ln n-z\right) -C_1c_0\Big|=0
$$

\noindent where $C_1$ is the constant in Proposition \ref{p:tailminkillmain}, and $c_0>0$ is defined in (\ref{c0}).
\end{proposition}

\noindent Looking at the set of particles that cross a high level $A>0$ for the first time, we then deduce the theorem for the constant $C^*=C_1c_0$.

\bigskip

The paper is organized as follows. Section \ref{s:many-to-one} introduces a useful and well-known tool, the many-to-one lemma. Then, Sections \ref{s:tailminkill}, \ref{s:tailmin} and \ref{s:main} contain respectively the proofs of Proposition \ref{p:tailminkillmain}, Proposition \ref{p:tailminmain} and Theorem \ref{main}. A sum-up of the notation used in the paper can be found in Appendix \ref{s:notation}.

\bigskip

Throughout the paper, $(c_i)_{i\ge 0}$ denote positive constants. We say that $a_n\sim b_n$ as $n\to \infty$ if $\lim_{n\to\infty} {a_n\over b_n}=1$. We write $\E[f,A]$ for $\E[f{\bf 1}_A]$, and we set $\sum_{\emptyset}:=0$, $\prod_{\emptyset} :=1$.

\section{The many-to-one lemma}
\label{s:many-to-one}

For $a\in\r$, we denote by $\P_a$ the probability distribution associated to the branching random walk starting from $a$, and $\E_a$ the corresponding expectation. Under (\ref{cond-bound}), there exists a centered random walk $(S_n,n\ge 0)$ such that for any $n\ge 1$, $a\in\r$ and any measurable function $g: \r^n \to [0, \, \infty)$,
\begin{equation}
    \E_a\Big[ \sum_{|x|=n} g(V(x_1), \cdots, V(x_n)) \Big]
    =
    \E_a \Big[ \ee^{S_n-a} g(S_1, \cdots, S_n)\Big]
    \label{many-to-one}
\end{equation}

\noindent where, under $\P_a$, we have $S_0=a$ almost surely. We will write $\P$ and $\E$ instead of $\P_0$ and $\E_0$ for brevity. In particular, under (\ref{cond-2}), $S_1$ has a finite variance $\sigma^2:=\E[S_1^2]=\E[\sum_{|x|=1} V(x)^2\ee^{-V(x)}]$. Equation (\ref{many-to-one}) is called in the literature the many-to-one lemma and can be seen as a consequence of Proposition \ref{p:spine} below.

\bigskip

\subsection{Lyons' change of measure}
\label{s:spine}
We introduce the {\it additive martingale}
\begin{equation}\label{def:Wn}
W_n:= \sum_{|u|=n} \ee^{-V(u)}.
\end{equation}

\noindent The fact that $W_n$ is a martingale comes from the branching property together with the assumption that $\E[\sum_{|x|=1}\ee^{-V(x)}]=1$. From \cite{lyons}, we know that $W_n$ converges almost surely as $n\to\infty$ to $0$ under our assumption (\ref{cond-bound}). For any $n\ge 0$, let  $\mathscr F_n$ denote the $\sigma$-algebra generated by the positions $(V(x),\,|x|\le n)$ up to time $n$, and $\mathscr F_\infty:=\bigvee_{n\ge 0}\mathscr F_n$. For any $a\in\r$, the Kolmogorov extension theorem guarantees that there exists a probability measure $\Q_a$ on $\mathscr F_\infty$ such that for any $n\ge 0$,
\begin{equation}\label{defQ}
\Q_a \, |_{\mathscr{F}_n} =\ee^a  W_n \bullet \P_a \, |_{\mathscr{F}_n}.
\end{equation}

\noindent We will write $\Q$ instead of $\Q_0$. We associate to the probability $\Q_a$ the expectation $\hat \E_a$. 
\medskip

We introduce the point process $\hat{ \mathcal L}$ with Radon-Nykodim derivative $\sum_{i\in \mathcal L} \ee^{-V(i)}$ with respect to the law of $\mathcal L$ and we consider the following process. At time $0$, the population is composed of one particle $w_0$ located at $V(w_0)=0$. Then, at each step $n$, particles of generation $n$ die and give birth to independent point processes distributed as $\mathcal L$, except for the particle $w_n$ which generates a point process distributed as $\hat {\mathcal L}$. The particle $w_{n+1}$ is chosen among the children of $w_n$ with probability proportional to $\ee^{-V(x)}$ for each child $x$ of $w_n$. This defines a branching random walk $\hat{\mathcal{B}}$ with a marked ray $(w_n)_{n\ge 0}$, which we call the {\it spine}. On the space of marked branching random walks, let $\hat {\mathscr F}_n$ be the $\sigma$-algebra generated by the positions $(V(x),|x|\le n)$ and the marked ray (or spine) $(w_k,k\le n)$ up to time $n$. Then, $\hat{\mathcal{B}}$ is measurable with respect to $\hat{\mathscr F}_\infty:=\bigvee_{n\ge 0} \hat{\mathscr F}_n$. We call $\mathcal{B}$ the natural projection of $\hat{\mathcal{B}}$ on the space of branching random walks without marked rays; in other words $ {\mathcal{B}}$ is obtained from $\hat{\mathcal{B}}$ by forgetting the identity of the spine. In particular, $ {\mathcal{B}}$ is measurable with respect to $\mathscr F_\infty$.  Notice that $\mathcal B$ is a branching random walk with immigration. We use the notation $a+ \hat {\mathcal{B}}$ or $a+{\mathcal{B}}$ to denote the branching random walk which positions are translated by $a$. 
\begin{proposition}\label{p:spine1}{\bf (\cite{lyons})}
 Under
 $\Q_a$, the branching random walk has the distribution of $ a+{\mathcal{B}}$.
\end{proposition}

Hence we will identify from now on our branching random walk under $\Q_a$ to the marked branching random walk $a+\hat{\mathcal B}$. Notice that by doing so, we introduce in our branching random walk a marked particle, the spine, and we extend the probability $\Q_a$ to $\hat{\mathscr F}_\infty$. We stress that in the filtration $(\mathscr F_n,n\ge 0)$, we do not know the identity of the spine. For $\ell\ge 1$, we call  $\I(w_\ell)$ the siblings of the spine at generation $\ell$: they are the vertices which share the same parent as $w_\ell$.
We will often use the $\sigma$-algebra 
\begin{eqnarray}\label{def:Gell}
\hat{\mathcal G}_\ell &:=& \sigma\{w_j,V(w_j),\Omega(w_j),(V(u))_{u\in \Omega(w_j)},j\in [1,\ell]\},\\
\hat{\mathcal G}_\infty &:=& \sigma\{w_j,V(w_j),\Omega(w_j),(V(u))_{u\in \Omega(w_j)},j\ge 1\}\label{def:Ginfty}
\end{eqnarray}
associated to the positions of the spine and its siblings, respectively up to time $\ell$ and up to time $\infty$. 

\begin{proposition}\label{p:spine}{\bf (\cite{lyons})}
 (i) For any $|x|=n$, we have
\begin{equation}
    \Q_a
    \{ w_n =x \, | \, \mathscr{F}_n \}
    =
     {\ee^{-V(x)}
     \over W_n}.
    \label{wn}
\end{equation}
(ii) The  process of the positions of the spine $(V(w_n), \, n\ge 0)$ under $\Q_a$ has the
 distribution of the centered random
 walk $(S_n, \, n\ge 0)$ under $\P_a$.
\end{proposition}

\bigskip

This change of probability was used in \cite{lyons}. We  refer to \cite{lyons-pemantle-peres} for the case of the Galton--Watson tree, to \cite{chauvin-rouault} for the analog for the branching Brownian motion, and to \cite{biggins-kyprianou04} for spine decompositions in various types of branching. Before closing this section, we collect some elementary facts about centered random walks with finite variance. We recall that we deal with non-lattice random walks. \\

There exists a constant $\alpha_1>0$ such that for any $x\ge 0$ and $n\ge 1$
\begin{equation}\label{eq:kozlov}
\P_x(\min_{j\le n} S_j \ge 0) \le \alpha_1 (1+x) n^{-1/2}.
\end{equation}

\noindent There exists a constant $\alpha_2>0$ such that for any $b\ge a\ge 0$, $x\ge 0$ and $n\ge 1$
\begin{equation}\label{lemmaA1}
\P_x(S_n \in [a,b],\, \min_{j\le n} S_j \ge 0) \le \alpha_2(1+x)(1+b-a)(1+b)n^{-3/2}.
\end{equation}

\noindent Let $0<\lambda<1$. There exists a constant $\alpha_3=\alpha_3(\lambda)>0$ such that for any $b\ge a\ge 0$, $x,y\ge 0$ and $n\ge 1$
\begin{eqnarray}\label{lemmaA3}
&& \P_x(S_n \in [y+a,y+b],\, \min_{j\le n} S_j \ge 0,\, \min_{\lambda n \le j\le n} S_j \ge y) \\
&\le& \alpha_3 (1+x)(1+b-a)(1+b)n^{-3/2}.\nonumber
\end{eqnarray}

\noindent Let $(a_n,\,n\ge 0)$ be a non-negative sequence such that $\lim_{n\to \infty} {a_n \over n^{1/2}}=0$. There exists a constant $\alpha_4>0$ such that for any $a\in[0, a_n]$ and $n\ge 1$
\begin{equation}\label{lemmaA3unif}
\P(S_n \in [a,a+1],\, \min_{j\le n} S_j \ge 0,\, \min_{n/2 < j\le n} S_j \ge a) \ge \alpha_4 n^{-3/2}.
\end{equation}

\noindent Equation (\ref{eq:kozlov}) is Theorem 1a, p.415 of \cite{feller}. Equations (\ref{lemmaA1}) and (\ref{lemmaA3}) are for example Lemmas 2.2 and 2.4 in \cite{ezratio}. Equation (\ref{lemmaA3unif}) is Lemma 4.3 of \cite{ezsimple}: even if the uniformity in $a\in[0,a_n]$ is not stated there, it follows directly from the proof.

\subsection{A convergence in law for the one-dimensional random walk}
\label{s:convlaw}
We recall that $(S_n)_{n\ge 0}$ is a non-lattice centered random walk under $\P$, with finite variance $\E[S_1^2]=\sigma^2\in (0,\infty)$. We introduce its renewal function $R(x)$ which is zero if $x<0$, $1$ if $x=0$, and for $x>0$
\begin{equation}\label{def:R(x)}
R(x):= \sum_{k\ge 0} \P(S_k \ge -x,\, S_k<\min_{0\le j\le k-1} S_j).
\end{equation}

\noindent If $H_n$ denotes the $n$-th strict descending ladder height (where by strict descending ladder height, we mean any $S_k$ such that $S_k<\min_{0\le j\le k-1} S_j$), then we observe that for $x\ge 0$,
\begin{equation}\label{eq:R(x)}
R(x)=\sum_{n\ge 0} \P(H_n\ge -x)
\end{equation}
\noindent which is $\E\left[\mbox{number of strict descending ladder heights which are } \ge -x\right]$. Similarly, we define $R_-(x)$ as the renewal function associated to $-S$. Since $\E[S_1]=0$ and $\E[S_1^2]<\infty$, we have that $\E[|H_1|]<\infty$ (see Theorem 1, Section XVIII.5 p.612 in \cite{feller}).  Then the renewal theorem \cite{feller}, p.360 implies that there exists  $c_0>0$ such that,
\begin{equation}\label{c0}
\lim_{x\to\infty} {R(x)\over x}=c_0.
\end{equation}

\noindent Moreover, there exist $C_-, C_+>0$ such that
\begin{eqnarray}
\label{eq:kozlov+} \P\left(\min_{1\le i\le n} S_i \ge 0\right) &\sim&  {C_+ \over \sqrt{n}}, \\
\label{eq:kozlov-} \P\left(\max_{1\le i\le n} S_i \le 0\right) &\sim&  {C_- \over \sqrt{n}}
\end{eqnarray}

\noindent as $n\to\infty$ (Theorem 1a, Section XII.7 p.415 of \cite{feller}).

\begin{lemma}\label{l:RWconvlaw}
Let $(r_n)_{n\ge 0}$ and $(\lambda_n)_{n\ge 0}$ be two sequences of numbers resp. in $\r_+$ and in $(0,1)$ and such that resp. $\lim_{n\to\infty} {r_n \over n^{1/2}}=0$, and $0<\liminf_{n\to\infty} \lambda_n\le \limsup_{n\to\infty} \lambda_n<1$. Let $F:\r_+\to\r$ be a Riemann integrable function. We suppose that there exists a non-increasing function ${\overline F}:\r_+\to \r$ such that $|F(x)| \le {\overline F}(x)$ for any $x\ge 0$ and $\int_{x\ge 0} x {\overline F}(x)<\infty$. Then, as $n\to\infty$,
\begin{equation}\label{eq:convRW}
\E\left[F(S_n-y),\min_{k\in[0,n]}S_k\ge 0,\min_{k\in [\lambda_n n,n]} S_k \ge y\right]
\sim
{C_-C_+ \sqrt{\pi} \over \sigma \sqrt{2}} n^{-3/2} \int_{x\ge 0}F(x)R_-(x)dx
\end{equation}

\noindent uniformly in $y\in[0,r_n]$.
\end{lemma}

\noindent {\it Proof}. Let $\varepsilon>0$. Since $|F(x)|\le {\overline F}(x)$ and ${\overline F}$ is non-increasing, we have for any integer $M\ge 1$,
\begin{eqnarray*}
&&\E\left[|F(S_n-y)|,\min_{k\in[0,n]}S_k\ge 0,\min_{k\in [\lambda_n n,n]} S_k \ge y, S_n\ge y+M\right]\\
&\le&
\sum_{j\ge M} {\overline F}(j) \P\left(\min_{k\in[0,n]} S_k\ge 0,\min_{k\in [\lambda_n n,n]} S_k \ge y, S_n\in [y+j,y+j+1)\right).
\end{eqnarray*}

\noindent For $ j\ge 1$, we have by (\ref{lemmaA3}) and the fact that $\limsup_{n\to\infty} \lambda_n<1 $,
\begin{eqnarray*}
\P\left(\min_{k\in[0,n]} S_k\ge 0,\min_{k\in [\lambda_n n,n]} S_k \ge y, S_n\in [y+j,y+j+1)\right)
\le c_1 {j\over n^{3/2}}.
\end{eqnarray*}

\noindent This yields that
\begin{eqnarray*}
\E\left[|F(S_n-y)|,\min_{k\in[0,n]}S_k\ge 0,\min_{k\in [\lambda_n n,n]} S_k \ge y, S_n\ge y+M\right]
\le {c_1\over n^{3/2}} \sum_{j\ge M} {\overline F}(j) j
\end{eqnarray*}

\noindent which is less than $\varepsilon n^{-3/2}$ for $M\ge 1$ large enough by the assumption that $\int_{x\ge 0}x{\overline F}(x)dx<\infty$. Therefore, we can restrict to $F$ with compact support. By approximating $F$ by scale functions ($F$ is Riemann integrable by assumption), we only prove (\ref{eq:convRW}) for $F(x)={\bf 1}_{\{x\in[0,a] \}}$, where $a\ge 0$. Let $a\ge 0$ be a fixed constant in the remainder of the proof. We have for such $F$
\begin{eqnarray*}
&&\E\left[F(S_n-y),\min_{k\in[0,n]}S_k\ge 0,\,\min_{k\in [\lambda_n n,n]} S_k \ge y\right]\\
&=&
\P(\min_{k\in[0,n]}S_k\ge 0,\,\min_{k\in [\lambda_n n,n]} S_k \ge y,\,S_n\le y+a).
\end{eqnarray*}

\noindent Let
$$
\phi_{y,a,n}(x):=\P_{x}\left(\min_{k\in [0,(1-\lambda_n)n]} S_k \ge y,\, S_{(1-\lambda_n)n}\le y+a \right).
$$

\noindent For $F(x)={\bf 1}_{\{x\in[0,a] \}}$, applying the Markov property at time $\lambda_n n$ (we assume that $\lambda_n n$ is integer for simplicity), we obtain that
\begin{equation}
\E\left[F(S_n-y),\min_{k\in[0,n]}S_k\ge 0,\,\min_{k\in [\lambda_n n,n]} S_k \ge y\right]
=
\E\left[  \phi_{y,a,n}(S_{\lambda_n n}) ,\min_{k\in[0,\lambda_n n]} S_k\ge 0 \right] \label{eq:conv0}.
\end{equation}

\noindent We estimate $\phi_{y,a,n}(x)$. Reversing time, we notice that
\begin{equation}\label{eq:convrever}
\phi_{y,a,n}(x)
=
\P\left(\min_{k\in [0,(1-\lambda_n)n]} (-S_k) \ge -S_{(1-\lambda_n)n} -(x-y)\ge -a\right).
\end{equation}

\noindent We introduce the strict descending ladder heights and times $(H_\ell^-,T_\ell^-)$ of $-S$ defined by $H_0^-:=0$, $T_0^-:=0$ and for any $\ell\ge 0$,
\begin{eqnarray*}
T_{\ell+1}^- &:=& \min\{k\ge T_\ell^- +1\,:\, (-S_k) < H_\ell^- \},\\
H_{\ell+1}^- &:=& -S_{T_{\ell+1}^-}.
\end{eqnarray*}

\noindent Since $\E[S_1]=0$ (and $\sigma>0$), we have $T_\ell^-<\infty$ for any $\ell \ge 0$ almost surely. Similarly to equation (\ref{eq:R(x)}), we have now $R_-(x) = \sum_{\ell\ge 0} \P(H_{\ell}^- \ge -x)$. Splitting the right-hand side of (\ref{eq:convrever}) depending on the value of the time $\ell$ for which $H_\ell^- = \min_{k\in[0,(1-\lambda_n) n]} (-S_k)$, we then have
\begin{eqnarray}\label{eq:conva}
&& \phi_{y,a,n}(x) \\
&=& \sum_{\ell\ge 0}\P\left(T_\ell^-\le (1-\lambda_n)n, H_\ell^- \ge -S_{(1-\lambda_n)n} -(x-y) \ge -a, \min_{k\in[T_\ell^-,(1-\lambda_n)n]} (-S_k) \ge H_\ell^- \right). \nonumber
\end{eqnarray}

\noindent By the strong Markov property at time $T_\ell^-$, we see that for any $h\in [-a,0]$ and $t\in [0,(1-\lambda_n )n)]$,
\begin{eqnarray*}
&&\P\left( H_\ell^- \ge -S_{(1-\lambda_n)n} -(x-y) \ge -a, \min_{k\in[T_\ell^-,(1-\lambda_n)n]} (-S_k) \ge H_\ell^-\,\Big|\, (H_\ell^-,\,T_\ell^-)=(h,t) \right)\\
&=&
{\bf 1}_{\{h\ge -a\}}\P\left(\min_{j\in[0,(1-\lambda_n)n-t]} (-S_j)\ge 0 , -S_{(1-\lambda_n)n-t}\in [(x-y)-a-h,(x-y)]\right).
\end{eqnarray*}

\noindent  Let $\psi(x):= x\ee^{-x^2/2}{\bf 1}_{\{x\ge 0\}}$. By Theorem 1 of \cite{caravenna} and equation (\ref{eq:kozlov-}), we check that
\begin{eqnarray*}
&&{\bf 1}_{\{h\ge -a\}}\P\left(\min_{j \in[0,(1-\lambda_n)n-t]} (-S_j)\ge 0 , -S_{(1-\lambda_n)n-t}\in [(x-y)-a-h,(x-y)]\right)\\
&=&
 {\bf 1}_{\{h\ge -a\}}{ C_- \over \sigma (1-\lambda_n) n}(h+a)\psi\left({x \over \sigma\sqrt{(1-\lambda_n)n}}\right) + {\bf 1}_{\{h\ge -a\}}o(n^{-1})
\end{eqnarray*}

\noindent uniformly in  $x\in \r$, $t\le n^{1/2}$, $h\in[-a,0]$ and $y\in[0,r_n] $. Here we used the fact that $\limsup_{n\to\infty} \lambda_n<1$. We mention that the cut-off $t\le n^{1/2}$ is arbitrary since the statement is valid for any $t=o(n)$. To deal with $t\in[n^{1/2},(1-\lambda_n)n]$, we see that
\begin{eqnarray*}
&&{\bf 1}_{\{h\ge -a\}}\P\left(\min_{j \in[0,(1-\lambda_n)n-t]} (-S_j)\ge 0 , -S_{(1-\lambda_n)n-t}\in [(x-y)-a-h,(x-y)]\right)\\
&=& {\bf 1}_{\{h\ge -a\}}O(h+a+1)\left((1-\lambda_n)n -t+1\right)^{-1}
\end{eqnarray*}

\noindent again by Theorem 1 of \cite{caravenna}. The last equation is valid uniformly in $x,y\in \r$, $t\in[0,(1-\lambda_n)n]$ and $h\in[-a,0]$. Going back to (\ref{eq:conva}), this implies that, for any $x\in\mathbb{R}$ and $y\in[0,r_n]$,
\begin{eqnarray*}
\phi_{y,a,n}(x)
&=&
o(n^{-1})+ {C_- \over \sigma (1-\lambda_n)n} \psi\left({x \over \sigma\sqrt{(1-\lambda_n)n}}\right) \sum_{\ell\ge 0} \E\left[ (H_\ell^- +a){\bf 1}_{\{H_\ell^- \ge -a,\,T_\ell^- \le n^{1/2}\}} \right]\\
&& + \,O(1)\sum_{\ell\ge 0} \E\left[ {H_\ell^- + a +1 \over (1-\lambda_n)n-T_\ell^-+1}{\bf 1}_{\{H_\ell^- \ge -a,\,T_\ell^-  \in (n^{1/2}, (1-\lambda_n)n]\}} \right] \\
\end{eqnarray*}

\noindent where we used the fact that $\sum_{\ell\ge 0}\P(H_\ell^-\ge -a)=R_-(a)=O(1)$ since $a$ is a constant. Observe that
$$
\sum_{\ell\ge 0} \E\left[ (H_\ell^-+a){\bf 1}_{\{H_\ell^- \ge -a,\,T_\ell^-> n^{1/2}\}} \right]
\le 
a \sum_{\ell\ge 0} \P(H_\ell^- \ge -a,\,T_\ell^-> n^{1/2})=o(1)
$$

\noindent as $n\to \infty$ by dominated convergence. Therefore,
\begin{eqnarray*}
\phi_{y,a,n}(x)
&=&
o(n^{-1})+{ C_- \over \sigma(1-\lambda_n) n} \psi\left({x \over \sigma\sqrt{(1-\lambda_n)n}}\right) \sum_{\ell\ge 0} \E\left[ (H_\ell^- + a){\bf 1}_{\{H_\ell^- \ge -a\}} \right]\\
&& + \,O(1)\sum_{\ell\ge 0} \E\left[ {H_\ell^- +a +1\over (1-\lambda_n)n-T_\ell^- + 1}{\bf 1}_{\{H_\ell^- \ge -a,\,T_\ell^- \in (n^{1/2}, (1-\lambda_n)n]\}} \right].
\end{eqnarray*}

\noindent  We want to show that the last term is $o(n^{-1})$ as well. We observe that
$$
\E\left[ {H_\ell^- +a +1\over (1-\lambda_n)n-T_\ell^-+1}{\bf 1}_{\{H_\ell^- \ge -a,\,T_\ell^- \in (n^{1/2}, (1-\lambda_n)n]\}} \right]
\le 
(a+1) \E\left[ { {\bf 1}_{\{H_\ell^- \ge -a,\,T_\ell^- \in (n^{1/2}, (1-\lambda_n)n]\}}\over (1-\lambda_n)n-T_\ell^-+1} \right]
.
$$

\noindent Since $\sum_{\ell\ge 0}\P\left(H_\ell^- \ge -a,\,T_\ell^- = k \right)\le \P(S_k\in[0,a],\,\min_{j\le k} S_j\ge 0)$, we obtain by (\ref{lemmaA1}) that 
$$
\sum_{\ell\ge 0}\P\left(H_\ell^- \ge -a,\,T_\ell^-=k \right)
\le 
\alpha_2(1+a)^2 k^{-3/2} 
$$

\noindent which yields that
\begin{eqnarray*}
&&\sum_{\ell\ge 0} \E\left[ {H_\ell^- + a +1\over (1-\lambda_n)n-T_\ell^- + 1}{\bf 1}_{\{H_\ell^-  \ge -a,\,T_\ell^- \in (n^{1/2}, (1-\lambda_n)n]\}} \right]\\
&\le&
 \alpha_2 (1+a)^3 \sum_{k=\lfloor n^{1/2}\rfloor + 1}^{\lfloor (1-\lambda_n)n \rfloor} k^{-3/2}{1\over (1-\lambda_n)n -k+1}
=
o(n^{-1}) 
\end{eqnarray*}

\noindent as we require. Therefore
$$
\phi_{y,a,n}(x)
=
o(n^{-1})+{ C_- \over \sigma(1-\lambda_n) n} \psi\left({x \over \sigma\sqrt{(1-\lambda_n)n}}\right) \sum_{\ell\ge 0} \E\left[ (H_\ell^- + a){\bf 1}_{\{H_\ell^-  \ge -a\}} \right]
$$

\noindent uniformly in $x\ge 0$ and $y\in[0,r_n]$. By (\ref{eq:kozlov}), we know that $\P(\min_{k\in [0,n]} S_k\ge 0)\le \alpha_1 n^{-1/2}$. It follows from equation (\ref{eq:conv0})  that 
\begin{eqnarray*}
&&\E\left[F(S_n-y),\min_{k\in[0,n]}S_k\ge 0,\,\min_{k\in [(1-\lambda_n)n,n]} S_k \ge y\right]\\
&=&
o(n^{-3/2}) + {C_- \over \sigma (1-\lambda_n)n}
\E\left[\psi\left({S_{\lambda_n n} \over \sigma\sqrt{(1-\lambda_n)n}}\right),\,\min_{k\in[0,\lambda_n n]} S_k\ge 0 \right]\sum_{\ell\ge 0} \E\left[ (H_\ell^- + a){\bf 1}_{\{H_\ell^- \ge -a\}} \right] .
\end{eqnarray*}

\noindent We know (see \cite{bolthausen}) that $S_n/(\sigma n^{1/2})$ conditioned on $\min_{k\in[0,n]} S_k$ being non-negative converges to the Rayleigh distribution. Therefore, 
\begin{eqnarray*}
\lim_{n\to\infty} \E\left[\psi\left({S_{\lambda_n n} \over \sigma\sqrt{(1-\lambda_n)n}}\right)\,\Big|\, \min_{k\in[0,\lambda_n n]} S_k\ge 0 \right] 
&=& \int_{x\ge 0} \psi\left(x\sqrt{\lambda_n \over 1-\lambda_n}\right)\psi(x)dx \\
&=& \sqrt{\lambda_n} (1-\lambda_n)\sqrt{\pi\over 2}.
\end{eqnarray*}

\noindent In view of (\ref{eq:kozlov+}), we get that, as $n\to\infty$,
$$
\E\left[\psi\left({S_{\lambda_n n} \over \sigma\sqrt{(1-\lambda_n)n}}\right),\,\min_{k\in[0,\lambda_n n]} S_k\ge 0 \right]\sim 
{C_+ (1-\lambda_n)\over \sqrt{n}} \sqrt{{\pi\over 2}}.
$$

\noindent We end up with
\begin{eqnarray*}
&&\E\left[F(S_n-y),\min_{k\in[0,n]}S_k\ge 0,\,\min_{k\in [\lambda_n n,n]} S_k \ge y\right]\\
&=&
o(n^{-3/2}) + { C_-C_+ \over \sigma n^{3/2}}\sqrt{{\pi\over 2}}\sum_{\ell\ge 0} \E\left[ (H_\ell^-+a){\bf 1}_{\{H_\ell^- \ge -a\}} \right]
\end{eqnarray*}

\noindent uniformly in $y\in[0,r_n]$. We recall that $\sum_{\ell\ge 0} \P(H_\ell^- \ge -a)=R_-(a)$ by definition and we took $F(x)={\bf 1}_{[0,a]}(x)$. By Fubini's theorem, it follows that $$\sum_{\ell\ge 0}  \E\left[ (H_\ell^-+a){\bf 1}_{\{H_\ell^- \ge -a\}} \right]=
\int_{x\ge 0}F(x)R_-(x)dx,$$

\noindent which completes the proof. \hfill $\Box$

\section{The minimum of a killed branching random walk}
\label{s:tailminkill}
It turns out to be useful to study first the killed branching random walk. Let 
$$
\T^{\rm{kill}} :=\{ u\in \T\,:\, V(u_k)\ge 0,\, \forall \, 0 \le k \le |u| \}
$$

\noindent be the set of individuals that stay above $0$.  We investigate the behaviour of the minimal position 
\begin{equation}\label{defminkill}
M_n^{\rm kill} := \min \{V(u),\,|u|^{\rm{kill}}=n \}
\end{equation}

\noindent where we write $|u|^{ \rm{kill} }$ to say that $u\in \T^{ \rm{kill}}$ and $|u|=n$.  If $M_n^{\rm kill}<\infty$, i.e. if the killed branching random walk survives until time $n$, we denote by $m^{{\rm kill},(n)}$ a vertex chosen uniformly in the set $\{u: |u|^{{\rm kill}}=n,\, V(u)=M_n^{\rm kill}\}$ of the particles that achieve the minimum.  It will be convenient to use the following notation, for $z\ge 0$:
\begin{eqnarray}
\label{def:an} a_n(z) &:=& {3\over 2} \ln n-z,\\
\label{def:In} I_n(z) &:=& [a_n(z)-1,a_n(z)),
\end{eqnarray}

\noindent and for $z\ge 0$, $0\le k\le n$ and $\lambda\in (0,1)$,
\begin{equation}\label{def:dk}
 d_k(n,z,\lambda) :=
\begin{cases}
0 , &\hbox{if $0\le k\le \lambda n$,} \cr
\max(a_n(z+1),0), &\hbox{if $\lambda n < k \le n$.}\cr
\end{cases}
\end{equation}

\noindent We will see later that, as $n\to \infty$, conditionally on being in $I_n(z)$, a particle that achieves the minimum at time $n$ did not cross the curve $k\to d_k(n,z+L,\lambda)$ with probability tending to $1$ when the constant $L$ goes to $\infty$ (and $\lambda$ is any constant in $(0,1)$). The section is devoted to the proof of the following proposition.
\begin{proposition}\label{p:tailminkill}
For any $\varepsilon>0$, there exist a real $A\ge 0$ and an integer $N\ge 1$ such that for any $n\ge N$ and $z\in [A,(3/2)\ln(n)-A]$,
$$
\Big|\,\ee^{z}\,\P(M_n^{\rm kill} \in I_n(z)) -C_2\Big| \le \varepsilon
$$
where $C_2$ is some positive constant.
\end{proposition}

\begin{corollary}\label{c:tailminkill}
Let $C_1:= {C_2 \over 1-\ee^{-1}}$. For any $\varepsilon>0$, there exist a real $A\ge 0$ and an integer $N\ge 1$ such that for any $n\ge N$ and $z\in[A,(3/2)\ln(n)-A]$,
$$
\Big|\,\ee^{z}\,\P\left(M_n^{\rm kill} < {3\over 2} \ln n -z\right) -C_1\Big| \le \varepsilon.
$$
\end{corollary}

Proposition \ref{p:tailminkillmain} immediately follows from Corollary \ref{c:tailminkill}. Assuming that Proposition \ref{p:tailminkill} holds, let us see how it implies the corollary.\\

\noindent {\it Proof of Corollary \ref{c:tailminkill}}. Let $\varepsilon>0$. We have by equation (\ref{many-to-one}), for any integer $n\ge 1$ and any real $r\ge 0$,
\begin{eqnarray*}
\E\left[ \sum_{|u|^{{\rm kill}}=n}{\bf 1}_{\{ V(u)\le  r \}} \right]
&=&
\E\left[\ee^{S_n},\, S_n \le r,\, \min_{0\le j\le n} S_j \ge 0  \right]\\
&\le&
\ee^{r} \P\left(S_n \le  r,\, \min_{0\le j\le n} S_j \ge 0\right).
\end{eqnarray*}

\noindent By (\ref{lemmaA1}), we have $\P\left(S_n \le  r,\, \min_{0\le j\le n} S_j \ge 0\right) \le c_2 {(1+r)^2 \over n^{3/2}}$. We deduce that 
\begin{equation}\label{eq:c(r)}
\P(M_n^{\rm kill} \le r) \le c(r) n^{-3/2}
\end{equation} 

\noindent with $c(r):=c_2\ee^{r}(1+ r)^2$. Let $A_1$ and $N_1$ be as in Proposition \ref{p:tailminkill}. We have for $n\ge N_1$ and $z\in[A_1,(3/2)\ln(n)-A_1]$,
$$
\Big|\,\P(M_n^{\rm kill} \in I_n(z)) -C_2\ee^{-z}\Big| \le \varepsilon\ee^{-z}.
$$

\noindent Summing this equation over  $z+k$ such that $z+k\in[A_1,(3/2)\ln(n)-A_1]$, we get that for any $n\ge N_1$ and $z\in[A_1,(3/2)\ln(n)-A_1]$,
$$
\Big|\,\P(M_n^{\rm kill} \in [r_{n,z},a_n(z))) -C_2\ee^{-z}\sum_{k=0}^{\lfloor a_n(z+A_1)\rfloor} \ee^{-k}\Big| \le \varepsilon\ee^{-z}\sum_{k=0}^{\lfloor a_n(z+A_1)\rfloor} \ee^{-k}
\le  {\varepsilon\over 1-\ee^{-1}} \ee^{-z}
$$

\noindent where  $r_{n,z}:=a_n(z)-\lfloor a_n(z+A_1)\rfloor\le A_1$. By (\ref{eq:c(r)}), we get that
$$
\Big|\,\P(M_n^{\rm kill} < a_n(z)) -C_2\ee^{-z}\sum_{k=0}^{\lfloor a_n(z+A_1)\rfloor} \ee^{-k}\Big| \le {\varepsilon\over 1-\ee^{-1}} \ee^{-z}+ c(A_1+1)n^{-3/2} .
$$

\noindent Let $A_2\ge A_1$ large enough such that $\sum_{k>A_2-A_1} \ee^{-k} \le \varepsilon $, and $c(A_1+1)\le \varepsilon \ee^{A_2}$. Then, for any $n\ge N_1$ and $z\in[A_1,(3/2)\ln(n)-A_2]$
\begin{eqnarray*}
\Big|\,\P(M_n^{\rm kill} < a_n(z)) - {C_2\over 1-\ee^{-1}}\ee^{-z}\Big| &\le& \varepsilon((1-\ee^{-1})^{-1} +1)\ee^{-z} +\varepsilon \ee^{A_2}n^{-3/2}\\
&\le&
\varepsilon((1-\ee^{-1})^{-1} +1)\ee^{-z} +\varepsilon \ee^{-z}
\end{eqnarray*}
which completes the proof. \hfill $\Box$

\subsection{Tightness of the minimum}
\label{ss:tightness}
Our aim is now to prove Proposition \ref{p:tailminkill}. In other words, we want to estimate the probability of the event $\{M^{{\rm kill}}_n \in I_n(z)\}$. The first lemma gives information on the path of particles located in $I_n(z)$. 

\begin{lemma}\label{l:tightL}
Let $0<\lambda<1$. There exist constants $c_3,c_4>0$ such that for any $n\ge 1$, $L\ge 0$, $x\ge 0$ and $z\ge 0$,
\begin{eqnarray}\label{eq:tightL}
 \\
\P_x\Big(\exists u\in \T^{{\rm kill}}\,:\,|u|=n ,\, V(u) \in I_n(z),\, \min_{k\in[\lambda n,n]} V(u_k) \in I_n(z+L)\Big)
\le c_{3}(1+x)\ee^{-{c_4L}}\ee^{-x-z}. \nonumber
\end{eqnarray}
\end{lemma}

\noindent {\it Proof}. Let $E$ be the event in (\ref{eq:tightL}), and write $d_k=d_k(n,z+L,\lambda)$ as defined in (\ref{def:dk}). Considering the time when the minimum $\min_{k\in[\lambda n,n]} V(u_k)$ is reached, we observe that $E \subset \bigcup_{k\in[\lambda n,n]} E_k$ where we define  $E_k:= \bigcup_{|u|=n} E_k(u) $ and for any $u\in \T$ with $|u|=n$,
$$
E_k(u) := \Big\{V(u_\ell) \ge d_\ell\,,\, \forall\,0\le \ell \le n,\, V(u)\in I_n(z),\,V(u_k)\in I_n(z+L)\Big\}.
$$

\noindent  Similarly, let
$$
E_k(S):=\Big\{S_\ell \ge d_\ell,  \forall\, 0\le \ell \le n,\,S_n \in I_n(z),\,S_k\in I_n(z+L)\Big\}.
$$

\noindent We notice that $\P_x(E_k) \le \E_x\left[ \sum_{|u|=n} {\bf 1}_{ E_k(u)}\right]$ which is $\E_x\left[\ee^{S_n-x}{\bf 1}_{ E_k(S) }  \right]$ by (\ref{many-to-one}).

\noindent In particular,
\begin{equation}\label{eq:EKS}
\P_x(E_k) \le n^{3/2}\ee^{-x-z}\P_x(E_k(S)).
\end{equation}

\noindent We need to estimate $\P_x(E_k(S))$. By the Markov property at time $k$,
\begin{eqnarray*}
\P_x(E_k(S))
&\le&
\P_x\left(S_\ell \ge d_\ell\,,\,\forall \, 0\le \ell \le k, S_k \in I_n(z+L)\right)\\
\nonumber &&\qquad \times \, \P\left(S_{n-k} \in[L-1,L+1],\, \min_{\ell \in[0,n-k]} S_\ell \ge -1\right).
\end{eqnarray*}

\noindent For the second term of the right-hand side, we know from (\ref{lemmaA1}) that there exists a constant $c_5>0$  such that
\begin{equation}\label{eq:tightL1a}
\P\left(S_{n-k} \in[L-1,L+1],\, \min_{\ell \in[0,n-k]} S_\ell \ge -1\right) \le c_5 (n-k+1)^{-3/2} (1+L) .
\end{equation}

\noindent To bound the first term, our argument depends on the value of $k$. Suppose that ${\lambda+1\over 2} n \le k\le n$. We have by (\ref{lemmaA3})
\begin{equation}\label{eq:tightL1b}
\P_x\left(S_\ell \ge d_\ell,\,\forall\, 0\le \ell \le k, S_k \in I_n(z+L)\right) \le c_6{(1+x)\over n^{3/2}}.
\end{equation}

\noindent If $\lambda n \le k< {\lambda+1\over 2} n$, we simply write
\begin{eqnarray}\nonumber
\P_x\left(S_\ell \ge d_\ell,\,\forall\, 0\le \ell \le k, S_k \in I_n(z+L)\right)
&\le&
\P_x\left(S_k \in I_n(z+L),\, \min_{\ell\in[0,k]} S_\ell \ge 0\right) \\
\label{eq:tightL1c}&\le&
c_7 (1+x) \ln(n) n^{-3/2}
\end{eqnarray}

\noindent by (\ref{lemmaA1}) . From (\ref{eq:tightL1a}), (\ref{eq:tightL1b}) and (\ref{eq:tightL1c}), there exists a constant $c_8>0$ such that
$$
\sum_{k\in[\lambda n,n-a]} \P_x(E_k(S)) \le c_8 (1+x)(1+L) n^{-3/2}\left(  a^{-1/2} +{\ln(n)\over \sqrt{n}}\right)
$$

\noindent for any $a\ge 1$. By (\ref{eq:EKS}), this yields that
\begin{equation} \label{Eksum1}
\sum_{k\in[\lambda n,n-a]} \P_x(E_k) \le c_8 (1+x)(1+L)\ee^{-x-z} \left(  a^{-1/2} +{\ln(n)\over \sqrt{n}}\right)
 .
\end{equation}

\noindent It remains to bound $\P_x(E_k)$ for $n-a < k\le n$. We observe that
$$
\P_x(E_k)
\le
\P_x\left(\exists\, |u|=k\,:\, V(u_\ell) \ge d_\ell\,,\,\forall\, 0\le \ell \le k, V(u)\in I_n(z+L) \right).
$$

\noindent By an application of (\ref{many-to-one}), we have
$$
\P_x(E_k) \le
n^{3/2}\ee^{-x-z-L}\P_x\left(S_\ell \ge d_\ell\,,\,\forall\, 0\le \ell \le k, S_k\in I_n(z+L)\right)
$$

\noindent which is $\le c_9 \ee^{-x-z-L}(1+x)$ by (\ref{lemmaA3}) (for $k\ge (1+\lambda)n/2$ for example).  It follows that, for $a\in [1, (1-\lambda)n/2]$,
\begin{equation}\label{Eksum2}
\sum_{k\in[n-a,n]} \P_x(E_k) \le c_9 (1+a) (1+x)\ee^{-x-z-L}.
\end{equation}

\noindent Equations (\ref{Eksum1}) and (\ref{Eksum2}) yield that, for any $n\ge 1$, $z,L\ge 0$ and $a\in [1, (1-\lambda)n/2]$,
\begin{equation}\label{eq:tightLend}
\P_x(E)
\le
\sum_{k\in[\lambda n,n]} \P_x(E_k)
\le
(1+x)\ee^{-x-z}\Big\{ c_8 (1+L)  \left(  a^{-1/2} +{\ln(n)\over \sqrt{n}}\right) + c_9 (1+a)\ee^{-L}\Big\}.
\end{equation}

\noindent Notice that (\ref{eq:tightL}) holds if $L> (3/2)\ln n$ since the left-hand side is $0$. If $L\le (3/2)\ln n$, take $a=\max(1,\alpha \ee^{\beta L})$ with $\alpha,\beta>0$ small enough and use (\ref{eq:tightLend}) to complete the proof. \hfill $\Box$

\bigskip

We deduce the following corollary. 
\begin{corollary}\label{cor:tightnesskill-upper}
We have for any $z,x\ge 0$ and any integer $n\ge 1$,
$$
\P_x\left(M_n^{{\rm kill}} \le a_n(z)\right) \le c_{10}(1+x)\ee^{-x-z}.
$$ 
\end{corollary}

\noindent Notice that by equation (\ref{many-to-one}),
\begin{eqnarray*}
\P(\exists \, u \in \T\,:\, V(u)\le -y)
&\le& \sum_{n\ge 0} \E\left[ \sum_{|u|=n}{\bf 1}_{\{ V(u)\le -y,V(u_k)>-y,\forall\, k<n \}} \right]\\
\nonumber &=&
\sum_{n\ge 0} \E\left[\ee^{S_n},S_n\le -y,S_k >-y \, \forall k<n \right]\\
\nonumber &\le& \ee^{-y}.
\end{eqnarray*}

\noindent Therefore, taking $x=y$ and $z=0$ in the corollary yields the following result.
\begin{corollary}\label{c:infinite}
We have for any $y\ge 0$ and any integer $n\ge 1$,
$$
\P\left(M_n \le a_n(y)\right) \le \left(1+c_{10} (1+y)\right)\ee^{-y}.
$$
\end{corollary}

\bigskip 

Recall that $a_n(z):={3\over 2} \ln n -z$, $I_n(z):=[a_n(z)-1,a_n(z))$ and $d_k(n,z+L,\lambda)$ is defined in (\ref{def:dk}).
\begin{definition}\label{def:ZzL}
For $u\in \T$, we say that $u\in \mathcal Z^{z,L}_n$ if
$$
|u|=n,\,V(u)\in I_n(z) \mbox{ and } V(u_k)\ge d_k\left(n,z+L,1/2\right)\,\forall\, k\le n.
$$
\end{definition}

\begin{figure}[h]
\begin{center}
\resizebox{11 cm}{6 cm}{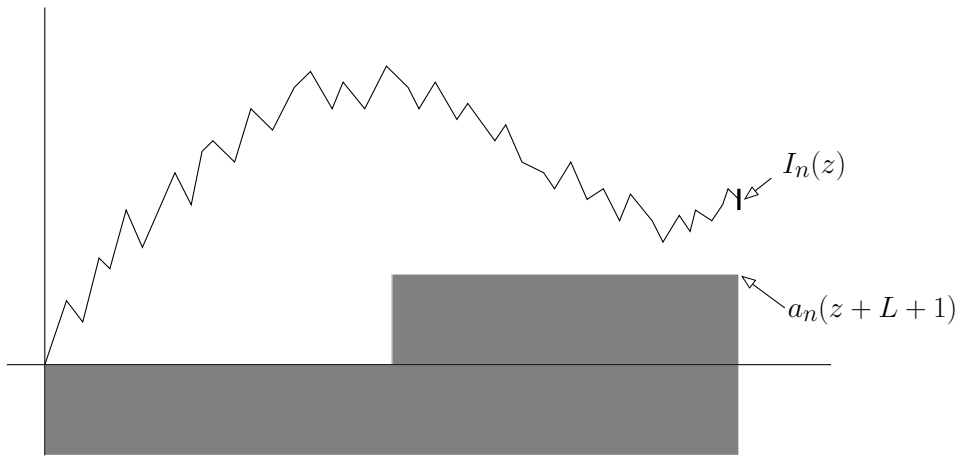}
\caption{ path of a vertex in $\mathcal Z^{z,L}_n$}
\label{f:ZzL}
\end{center}
\end{figure}

\noindent Notice that if $u\in \mathcal Z_n^{z,L}$, then necessarily $u\in \T^{{\rm kill}}$. In words, $u\in\mathcal Z^{z,L}_n$ means that a particle is located around $\frac32\ln n-z$, and did not cross the curve $k\to d_k(n,z+L,1/2)$. We deduce from Lemma \ref{l:tightL} that for any $\varepsilon>0$, there exists $L_0>0$ such that for any $n\ge 1$, $L\ge L_0$ and $z\ge 0$,
\begin{equation}\label{cor:tightL}
\P\left(\exists\, u \in \T^{{\rm kill}}\,:\, |u|=n,\,u\notin \mathcal Z^{z,L}_n,\, V(u) \in I_n(z) \right) \le \varepsilon \ee^{-z}.
\end{equation}

\noindent Equivalently, with high probability, any particle of the killed branching random walk located around $\frac32\ln n-z$ stayed above the curve $k\to d_k(n,z+L,1/2)$. We show now that $\P(M_n^{\rm kill} \le\frac32 \ln n -z)$ has an exponential decay as $z\to\infty$. Corollary \ref{cor:tightnesskill-upper} gives an upper bound. The following lemma gives a lower bound.

\begin{lemma}\label{l:tightnesskill}
There exists $c_{11}>0$ such that for any $n\ge 1$ and $z\in [0,(3/2)\ln n-1]$
$$
\P\left(M_n^{\rm kill} < {3\over 2} \ln n -z\right) \ge c_{11}\ee^{-z}.
$$
\end{lemma}

\noindent {\it Proof}. The proof relies on a second moment argument. Let $z\in [0,(3/2)\ln n-1]$ and $n\ge 1$. For $1\le k \le n$, let
$$
e_k=e_k^{(n)}:=
\begin{cases}
k^{1/12}, &\hbox{if $1\le k\le {n\over 2}$,} \cr
(n-k)^{1/12}, &\hbox{if ${n\over 2}<k\le n$}\cr
\end{cases}
$$

\noindent and write for brevity $d_k=d_k(n,z,1/2)$. In order to have good bounds in our second moment argument, we will restrict to 'good' vertices which do not have 'too many' descendants. This leads us to the following definition. We say that $|u|=n$ is a $z$-good vertex if $u \in \mathcal Z^{z,\,0}_n$ and
\begin{equation}\label{def:good}
\sum_{v\in \Omega(u_k)}\ee^{-(V(v)-d_k)}\Big(1+(V(v)-d_{k})_+\Big) \le B\ee^{-e_k} \qquad \forall\, 1\le k\le n,
\end{equation}

\noindent where $\Omega(y)$ stands for the set of siblings of $y$, i.e the particles $x\neq y$ which share the same parent as $y$ in the tree $\T$. The number $B>0$ is a constant that we will fix later on. The reason of such a definition becomes clear in the computation of the second moment in (\ref{eq:goodmoment2}). Such conditions on the behavior of the children off the path of the spine in a second moment argument are not new, and were already used in \cite{gantert-hu-shi}.

Remember the probability measure $\Q$ that we introduced in Section \ref{s:spine}, which is associated to the expectation $\hat{\E}$. We recall that $w_n$ is the spine at generation $n$, and we know from Proposition \ref{p:spine} (ii) that $(V(w_k),k\ge 0)$ under $\Q$ has the law of the centered random walk $(S_k,k\ge 0)$ under $\P$. By (\ref{lemmaA3unif}) with $a_n=(3/2)\ln n$, there exists $c_{13}>0$ such that $\Q(w_n\in\mathcal Z^{z,\,0}_n)\ge 2c_{13}n^{-3/2}$. Then, by Lemma \ref{l:goodvertex}, we can choose $B>0$ such that for any $n\ge 1$ and $z\in[0,(3/2)\ln n-1]$
$$
\Q(w_n \mbox{ is a $z$-good vertex} ) \ge c_{13}n^{-3/2}.
$$

\noindent Let ${\rm Good}_n$ be the number of $z$-good vertices at generation $n$. We have by definition of the measure $\Q$ then Proposition \ref{p:spine} (i),
\begin{eqnarray*}\nonumber
\E\left[ {\rm Good}_n \right]
&=&
\hat{\E}\left[ {1\over W_n}\sum_{|u|=n} {\bf 1}_{\{ u \;\mbox{\footnotesize is a $z$-good vertex} \}} \right]\\
\nonumber &=& \hat{\E}\left[ \ee^{V(w_n)},\, w_n \;\mbox{\small is a $z$-good vertex} \right].
\end{eqnarray*}

\noindent On the event that $w_n \in \mathcal Z^{z,0}_n$, we have that $V(w_n)\ge (3/2)\ln(n)-z-1$. Therefore,
\begin{equation}
\E\left[ {\rm Good}_n \right]
\ge
n^{3/2}\ee^{-z-1}\Q\left( w_n\; \mbox{\small is a $z$-good vertex} \right)
\label{eq:killmoment1} \ge
c_{13} \ee^{-z-1}.
\end{equation}

\noindent We look at the second moment. We use again Proposition \ref{p:spine} (i) to see that
\begin{eqnarray*}
\E\left[ ({\rm Good}_n)^2 \right]
&=&
\hat{\E}\left[\ee^{V(w_n)}{\rm Good}_n,\, w_n \; \mbox{\small is a $z$-good vertex}  \right]\\
&\le&
n^{3/2}\ee^{-z}\hat{\E}\left[{\rm Good}_n,\, w_n \; \mbox{\small is a $z$-good vertex}  \right]
\end{eqnarray*}

\noindent since $V(w_n)\le (3/2)\ln(n)-z$ when $w_n \in \mathcal Z^{{z,0}}_n$. Let $Y_n$ be the number of vertices $u$ such that $u\in \mathcal Z^{z,0}_n$. We notice that $Y_n\ge {\rm Good}_n$, hence
$$
\E\left[ ({\rm Good}_n)^2 \right]
\le
n^{3/2}\ee^{-z}\hat {\E}\left[Y_n,\, w_n \; \mbox{\small is a $z$-good vertex}  \right].
$$

\noindent We decompose $Y_n$ along the spine. We get
$$
Y_n = {\bf 1}_{\{ w_n \in \mathcal Z^{z,\,0}_n \}} + \sum_{k=1}^{n} \sum_{u\in \Omega(w_k)} Y_n(u)
$$

\noindent where $Y_n(u)$ is the number of vertices $v$ which are descendants of $u$ and such that $v \in \mathcal Z^{z,\,0}_n$. Therefore,
\begin{eqnarray}\label{eq:goodmoment2bis}
\E\left[ ({\rm Good}_n)^2 \right]
&\le&
n^{3/2}\ee^{-z} \Big( \Q(w_n \; \mbox{is a $z$-good vertex})\\
&& +\sum_{k=1}^{n} \hat{\E}\left[ \sum_{u\in \Omega(w_k)} Y_n(u),w_n \; \mbox{\small is a $z$-good vertex}\right]\Big).\nonumber
\end{eqnarray}

\noindent Recall from (\ref{def:Ginfty}) that $\hat{\mathcal G}_\infty$ is the $\sigma$-algebra generated by the spine and its siblings. Recall that the branching random walk rooted at $u\in \Omega(w_k)$ has the same law under $\P$ and $\Q$. For $u\in \Omega(w_k)$, we have $Y_n(u)=0$ if there exists $j\le |u|$ such that $V(u_j) < d_j$. Otherwise, we have by (\ref{many-to-one}),
\begin{eqnarray*}
\hat{\E}\left[Y_n(u)\,|\, \hat{\mathcal G}_\infty\right]
&=&
\E_{V(u)}\left[\sum_{|v|=n-k} {\bf 1}_{\{ V(v_j)\ge d_{k+j},\,\forall \,0\le j\le n-k,\,V(v)\in I_n(z) \}}\right]\\
&=&
\ee^{-V(u)}\E_{V(u)}\left[\ee^{S_{n-k}}, S_j\ge d_{k+j},\forall\, 0\le j\le n-k,S_{n-k}\in I_n(z)\right].
\end{eqnarray*}

\noindent Consequently,
\begin{eqnarray*}
\hat{\E}[Y_n(u)\,|\, \hat{\mathcal G}_\infty]
&\le&
n^{3/2}\ee^{-z-V(u)}\P_{V(u)}\left(S_j\ge d_{k+j},\forall \,0\le j\le n-k,S_{n-k}\in I_n(z)\right)\\
&=:&
n^{3/2}\ee^{-z-V(u)} {\rm p}(V(u),k,n,z)
\end{eqnarray*}

\noindent the latter inequality consisting the definition of ${\rm p}(V(u),k,n,z)$. Hence, equation (\ref{eq:goodmoment2bis}) gives that
\begin{eqnarray}\nonumber
\E\left[ ({\rm Good}_n)^2 \right]
&\le&
n^{3/2}\ee^{-z} \Big( \Q(w_n \; \mbox{is a $z$-good vertex}) +\\
 && n^{3/2}\ee^{-z}\sum_{k=1}^{n} \hat{\E}\left[ \sum_{u\in \Omega(w_k)}\ee^{-V(u)} {\rm p}(V(u),k,n,z) ,w_n \; \mbox{\small is a $z$-good vertex}\right]\Big).\label{eq:goodmoment2}
\end{eqnarray}

\noindent We want to bound ${\rm p}(r,k,n,z)$ for $r\in \r$. We have to split the cases  $k\le n/2$ and $n/2<k\le n$. Suppose first that $k\le n/2$. Then ${\rm p}(r,k,n,z)=0$ if $r<0$. If $r\ge 0$, we apply (\ref{lemmaA3}) to see that for any $n\ge 1$, $k\le n/2$, $r\ge 0$ and $z\ge 0$, 
$$
 {\rm p}(r,k,n,z) \le c_{14} (r+1)n^{-3/2}.
$$

\noindent This implies that, for any  $n\ge 1$, $k\le n/2$, $r\ge 0$ and $z\ge 0$, 
\begin{eqnarray*}
&&\sum_{k=1}^{\lfloor {n\over 2} \rfloor} \hat{\E}\left[ \sum_{u\in \Omega(w_k)} \ee^{-V(u)} {\rm p}(V(u),k,n,z) ,w_n\; \mbox{\small is a $z$-good vertex}\right]\\
&\le&
c_{14}n^{-3/2} \sum_{k=1}^{\lfloor {n\over 2} \rfloor} \hat{\E}\left[ \sum_{u\in \Omega(w_k)} \ee^{-V(u)} (1+V(u)_+) ,w_n \;\mbox{ \small is a $z$-good vertex}\right]\\
&\le&
c_{14}B n^{-3/2}\sum_{k=1}^{\lfloor {n\over 2} \rfloor} \ee^{-e_k} \Q\left(w_n \;\mbox{\small is a $z$-good vertex}\right)
\end{eqnarray*}

\noindent where the last inequality comes from the property (\ref{def:good}) satisfied by a good vertex. When $n/2<k\le n$, we simply write ${\rm p}(r,k,n,z)\le 1$ and we get

\begin{eqnarray*}
&& \sum_{k=\lfloor {n\over 2} \rfloor+1}^{n} \hat{\E}\left[ \sum_{u\in \Omega(w_k)} \ee^{-V(u)} {\rm p}(V(u),k,n,z) ,\;w_n \;\mbox{\small is a $z$-good vertex}\right]\\
&\le&
\sum_{k=\lfloor {n\over 2} \rfloor+1}^{n} \hat{\E}\left[ \sum_{u\in \Omega(w_k)} \ee^{-V(u)} ,\;w_n \;\mbox{\small is a $z$-good vertex}\right]\\
&=&
n^{-3/2}\ee^{z+1} \sum_{k=\lfloor {n\over 2} \rfloor+1}^n \hat{\E}\left[ \sum_{u\in \Omega(w_k)} \ee^{-(V(u)-d_k)} ,\;w_n \;\mbox{\small is a $z$-good vertex}\right]\\
&\le&
B n^{-3/2}\ee^{z+1}\sum_{k=\lfloor {n\over 2} \rfloor+1}^n \ee^{-e_k} \Q\left(w_n \; \mbox{\small is a $z$-good vertex}\right)
\end{eqnarray*}

\noindent by (\ref{def:good}). Going back to (\ref{eq:goodmoment2}), we deduce that for any $z\ge 0$ and $n\ge 1$,
\begin{eqnarray*}
\E\left[ ({\rm Good}_n)^2 \right]
&\le&
 n^{3/2}\ee^{-z} \Big\{ 1+ c_{15}\sum_{k=1}^{n} \ee^{-e_k} \Big\}\Q\left(w_n\; \mbox{\small is a $z$-good vertex}\right)\\
 &\le&
 c_{16}  n^{3/2}\ee^{-z}\Q\left(w_n \;\mbox{\small is a $z$-good vertex}\right).
\end{eqnarray*}

\noindent Now, observe that $\Q(w_n \; \mbox{\small is a $z$-good vertex}) \le \Q(w_n \in \mathcal Z^{z,0}_n) \le c_{17}n^{-3/2} $ by Definition \ref{def:ZzL} and equation (\ref{lemmaA3}). Hence
\begin{equation}\label{eq:killmoment2}
\E\left[ ({\rm Good}_n)^2 \right] \le c_{18}\ee^{-z}.
\end{equation}

\noindent By the Paley-Zygmund inequality, we have $\P({\rm Good}_n\ge 1) \ge {\E[{\rm Good}_n]^2 \over \E[({\rm Good}_n)^2]}$ which is greater than  $c_{19}\ee^{-z}$ by (\ref{eq:killmoment1}) and (\ref{eq:killmoment2}). We conclude by observing that if ${\rm Good}_n \ge 1$ then $M_n^{\rm kill} < {3\over 2}\ln n -z$. \hfill $\Box$

\subsection{Proof of Proposition \ref{p:tailminkill}}

Corollary \ref{cor:tightnesskill-upper} and Lemma \ref{l:tightnesskill} already give the right rate of decay, but we want to strenghten it into an asymptotic as $z\to\infty$. We recall that $m^{{\rm kill},(n)}$ is chosen uniformly among the particles in $\T^{\rm kill}$ that achieve the minimum. We introduced the notation $\mathcal Z^{z,L}_n$ in Definition \ref{def:ZzL}.
By (\ref{cor:tightL}), we have that with high probability $m^{{\rm kill},(n)} \in\mathcal Z^{z,L}_n$ whenever $M_n^{\rm kill} \in I_n(z)$, where $L$ is a large constant. The first step of the proof is to give a representation of the probability $\P\left( M_n^{\rm kill} \in I_n(z),\, m^{{\rm kill},(n)} \in \mathcal Z^{z,L}_n \right)$ in terms of the spine decomposition presented in Section \ref{s:spine}. Recall that the notation $|u|^{\rm kill}=n$ is a short way to say that $u\in \T^{\rm kill}$ and $|u|=n$.

\begin{lemma}\label{l:step1}
For any $z\ge 0$, $L\ge 0$, and $n\ge 1$, we have
\begin{equation}\label{eq:step1}
\P\left( M_n^{\rm kill} \in I_n(z),\, m^{{\rm kill},(n)} \in \mathcal Z^{z,L}_n \right) =  \hat{\E}\left[{\ee^{V(w_n)}{\bf 1}_{\{ V(w_n)=M_n^{\rm kill}\}}\over \sum_{|u|^{\rm kill}=n} {\bf 1}_{\{ V(u)=M_n^{\rm kill}\}}}, \, w_n\in \mathcal Z^{z,L}_n \right].
\end{equation}
\end{lemma}

\noindent {\it Proof}. We observe that
\begin{eqnarray*}
\P\left( M_n^{\rm kill} \in I_n(z), \, m^{{\rm kill},(n)} \in \mathcal Z^{z,L}_n\right)
&=&
\E\left[\sum_{|u|=n} {\bf 1}_{\{ u=m^{{\rm kill},(n)},u\in \mathcal Z^{z,L}_n\}}\right]\\
&=&
\E\left[{\sum_{|u|=n} {\bf 1}_{\{V(u) = M_n^{\rm kill},\, u\in \mathcal Z^{z,L}_n\}} \over \sum_{|u|^{\rm kill}=n} {\bf 1}_{\{V(u)=M_n^{\rm kill} \}} }\right].
\end{eqnarray*}

\noindent  Using the measure $\Q$, it follows from  Proposition \ref{p:spine} (i) that
$$
\E\left[{\sum_{|u|=n} {\bf 1}_{\{V(u) = M_n^{\rm kill},\, u\in \mathcal Z^{z,L}_n\}} \over \sum_{|u|^{\rm kill}=n} {\bf 1}_{\{V(u)=M_n^{\rm kill} \}} }\right]
=
\hat {\E}\left[{\ee^{V(w_n)}\over \sum_{|u|^{\rm kill}=n} {\bf 1}_{\{ V(u)=M_n^{\rm kill}\}}}{\bf 1}_{\{V(w_n) = M_n^{\rm kill}, w_n\in \mathcal Z^{z,L}_n\}}\right],
$$

\noindent which completes the proof. \hfill $\Box$

\bigskip

We now study our branching random walk under $\Q$, which we identified with the branching random walk $\hat{\mathcal B}$ by the mean of Proposition \ref{p:spine1}. For $b\le n$ integers and $z\ge 0$, we define the event $\mathcal E_n(z,b)\in\hat{\mathscr F}_n$ by
\begin{equation}\label{def:En}
\mathcal E_n(z,b) :=\{ \forall\, k\le n-b,\,\forall\, v\in \Omega(w_k), \min_{u\ge v,|u|^{\rm kill}=n} V(u)\ge a_n(z)\}
\end{equation}

\noindent where, as before, $\Omega(w_k)$ denotes the set of siblings of $w_k$. On the event $\mathcal E_n(z,b)\cap\{M_n^{\rm kill} \in I_n(z)\}$, we are sure that any particle located at the minimum separated from the spine after the time $n-b$. The following lemma will be proved in Section \ref{s:proofEn}.

\begin{lemma}\label{l:En}
Let $\eta>0$ and $L\ge 0$. There exist $A>0$ and $B\ge 1$ such that for any integers $n\ge b\ge B$ and any real $z\ge A$,
\begin{equation}
\label{eq:En} \Q((\mathcal E_n(z,b))^c,\, w_n \in \mathcal Z^{z,L} _n) \le \eta n^{-3/2}.
\end{equation}
\end{lemma}

 Let, for $x\ge 0$, $L\ge 0$, and any integer $b\ge 1$
\begin{equation}\label{def:FLb}
F_{L,b}(x)
:=
\hat{\E}_x\left[{\ee^{V(w_b)-L}{\bf 1}_{\{V(w_b)=M_b\}} \over \sum_{|u|=b} {\bf 1}_{\{ V(u)=M_b \}}},\,\min_{k\in [0,b]} V(w_k)\ge -1,\,V(w_b)\in [L-1,L)\right].
\end{equation}

\noindent We stress that $M_b$ which appears in the definition of $F_{L,b}(x)$ is the minimum at time $b$ of the {\bf non-killed} branching random walk. Then, define
\begin{equation}\label{def:CLb}
C_{L,b} := {C_-C_+\sqrt{\pi} \over \sigma \sqrt{2}} \int_{x\ge 0} F_{L,b}(x)R_-(x) dx,
\end{equation}

\noindent where $C_-$, $C_+$ and $R_-(x)$ were defined in Section \ref{s:convlaw}. We recall that, by Proposition \ref{p:spine} (ii), the spine has the law of $(S_n)_{n\ge 0}$. In (\ref{def:FLb}), we see that ${{\bf 1}_{\{V(w_b)=M_b\}} \over \sum_{|u|=b} {\bf 1}_{\{ V(u)=M_b \}}}$ is smaller than $1$, and $\ee^{V(w_b)-L}\le 1$. Hence, $|F_{L,b}(x)|\le \P(S_b \le L-x)=:{\overline F}(x)$ which is non-increasing in $x$, and $\int_{x\ge 0} {\overline F}(x)xdx = {1\over 2}\E[ (L-S_b)^2{\bf 1}_{\{S_b \le L\}}]<\infty$. Moreover, changing the starting point from $x$ to $0$, we observe that
$$
F_{L,b}(x) = \ee^x \hat{\E}\left[{\ee^{V(w_b)-L}{\bf 1}_{\{V(w_b)=M_b\}} \over \sum_{|u|=b} {\bf 1}_{\{ V(u)=M_b \}}}{\bf 1}_{\{\min_{k\in [0,b]} V(w_k)\ge -x-1,\,V(w_b)\in [-x+L-1,-x+L)\}}\right].
$$

\noindent The fraction in the expectation is smaller than $1$. Using the identity $|{\bf 1}_E -a{\bf 1}_F|\le 1-a + |{\bf 1}_E - {\bf 1}_F|$ for $a\in(0,1)$, this yields that for $x\ge 0$, $\varepsilon>0$ and any $y\in [x,x+\varepsilon]$,
\begin{eqnarray*}
&& |F_{L,b}(y)-F_{L,b}(x)|\\
&\le& \ee^y\E\left[\Big|{\bf 1}_{\{\min_{k\in [0,b]} S_k\ge -y-1,\,S_b+y-L\in [-1,0)\}} - \ee^{x-y} {\bf 1}_{\{\min_{k\in [0,b]} S_k\ge -x-1,\,S_b+x-L\in [-1,0)\}}\Big|\right]\\
&\le&
\ee^y(1-\ee^{-\varepsilon}) +\ee^y \E\left[{\bf 1}_{\{\min_{k\in [0,b]} S_k+x+1 \in [-\varepsilon,0)\}} + {\bf 1}_{\{S_b+x-L \in [-1-\varepsilon,-1)\cup(-\varepsilon,0] \}}\right] 
\end{eqnarray*}

\noindent from which we  deduce that $x\to F_{L,b}(x)$ is Riemann integrable. Therefore, $F_{L,b}$ satisfies the conditions of Lemma \ref{l:RWconvlaw} for any $L\ge 0$ and integer $b\ge 1$.\\

We want to prove that the expectation in (\ref{eq:step1}) behaves like $\ee^{-z}$ with some constant factor, as $z\to\infty$. By Lemma \ref{l:En}, we can restrict to the event $\mathcal E_n(z,b)$. The next lemma shows that the expectation on this event is then equivalent to $C_{L,b}\ee^{-z}$.

\begin{lemma}\label{l:step2}
 Let $L\ge 0$ and $\eta>0$. Let $A$ and $B$ be as in Lemma \ref{l:En}. For any integer $b\ge B$, we can find a constant $H>0$ such that for $n$ large enough, and $z\in[A, (3/2)\ln (n) -L-H]$,
\begin{equation}\label{eq:Clbh}
\Bigg| \ee^z \hat{\E}\left[{\ee^{V(w_n)}{\bf 1}_{\{ V(w_n)=M_n^{\rm kill}\}}\over \sum_{|u|^{\rm kill}=n} {\bf 1}_{\{ V(u)=M_n^{\rm kill}\}}},\, w_n \in \mathcal Z^{z,L}_n,\,\mathcal E_n(z,b)\right]
 - C_{L,b}\Bigg| \le 3 \eta.
\end{equation}
\end{lemma}

\noindent {\it Proof}. Let $L$, $\eta$, $A$, $B$ be as in the lemma. Throughout the proof, $b$ is a fixed integer which is greater than $B$. We denote by $\hat {\E}_{\hbox{\scriptsize (\ref{eq:Clbh})}}$ the expectation in (\ref{eq:Clbh}). Recall that under $\Q$, our process is identified with $\hat{\mathcal B}$.  Applying the branching property at the vertex $w_{n-b}$ to $\hat{\mathcal B}$, we have for any $n\ge b$ and $z\ge 0$,
\begin{eqnarray*}
\hat{\E}_{\hbox{\scriptsize (\ref{eq:Clbh})}}
=
\hat{\E}\Big[F^{\rm kill}(V(w_{n-b})),\, V(w_\ell) \ge d_{\ell},\, \forall\, \ell \le n-b,\, \mathcal E_n(z,b)\Big]
\end{eqnarray*}

\noindent where $d_\ell:=d_\ell(n,z+L,1/2)$ (see (\ref{def:dk})) and  $F^{\rm kill}$ is defined for $x\ge 0$ by
\begin{equation}\label{def:tildeFkill}
 F^{\rm kill}(x)
:=
\hat{\E}_x\left[{\ee^{V(w_b)}{\bf 1}_{\{V(w_b)=M_b^{\rm kill}\}} \over \sum_{|u|^{\rm kill}=b} {\bf 1}_{\{ V(u)=M_b^{\rm kill} \}}}, \, \min_{k\in [0,b]} V(w_k)\ge a_n(z+L+1),\,V(w_b)\in I_n(z)\right].
\end{equation}

\noindent Notice that $F^{\rm kill}(x)\le n^{3/2}\ee^{-z}\Q_x(\min_{k\in [0,b]} V(w_k)\ge a_n(z+L+1),\,V(w_b)\in I_n(z))$. Hence
\begin{eqnarray*}\nonumber
&&\Big|\hat{\E}_{\hbox{\scriptsize (\ref{eq:Clbh})}} - \hat{ \E}\Big[ F^{\rm kill}(V(w_{n-b})),\, V(w_\ell) \ge d_{\ell},\, \forall \,\ell \le n-b\Big]\Big|\\
\nonumber &=&
\hat{\E}\Big[F^{\rm kill}(V(w_{n-b})),\, V(w_\ell) \ge d_{\ell},\, \forall\, \ell \le n-b,\,(\mathcal E_n)^c\Big] \\
&\le&
n^{3/2}\ee^{-z} \hat{\E}\Big[\Q_{V(w_{n-b})}\Big(\min_{k\in [0,b]} V(w_k)\ge a_n(z+L+1),\,V(w_b)\in I_n(z)\Big){\bf 1}_{\{ V(w_\ell) \ge d_{\ell},\, \forall\, \ell \le n-b\}\cap(\mathcal E_n)^c}\Big]
\end{eqnarray*}

\noindent where we wrote $\mathcal E_n$ for $\mathcal E_n(z,b)$. By the Markov property, the term $$\hat{\E}\Big[\Q_{V(w_{n-b})}\Big(\min_{k\in [0,b]} V(w_k)\ge a_n(z+L+1),\,V(w_b)\in I_n(z)\Big){\bf 1}_{\{ V(w_\ell) \ge d_{\ell},\, \forall\, \ell \le n-b\},\,(\mathcal E_n)^c}\Big]$$
is equal to $\Q\Big(w_n\in \mathcal Z^{z,L}_n,\,(\mathcal E_n(z,b))^c\Big)$ which is at most $\eta n^{-3/2}$ when $z\ge A$ and $n\ge b$ by Lemma \ref{l:En} and our choice of $A$ and $B$. Therefore, for any $n\ge b$ and $z\ge A$,
\begin{equation}
\Big|\hat{\E}_{\hbox{\scriptsize (\ref{eq:Clbh})}} - \hat{\E}\Big[F^{\rm kill}(V(w_{n-b})),\, V(w_\ell) \ge d_{\ell},\, \forall \,\ell \le n-b\Big]\Big|
\le
\eta\ee^{-z}. \label{eq:step2a}
\end{equation}

\noindent Recall the definition of $F_{L,b}$ in (\ref{def:FLb}). We would like to replace $F^{\rm kill}(x)$ by $n^{3/2}\ee^{-z} F_{L,b}(x -a_n(z+L))$. We notice that
\begin{eqnarray*}
&&n^{3/2}\ee^{-z} F_{L,b}(x - a_n(z+L))\\
&=&\hat {\E}_x\left[{\ee^{V(w_b)}{\bf 1}_{\{V(w_b)=M_b\}} \over \sum_{|u|=b} {\bf 1}_{\{ V(u)=M_b \}}}, \, \min_{k\in [0,b]} V(w_k)\ge a_n(z+L+1),\,V(w_b)\in I_n(z)\right].
\end{eqnarray*}

\noindent We observe that the only difference with (\ref{def:tildeFkill}) is that the branching random walk is not killed anymore. Since $\Big|{{\bf 1}_{\{V(w_b)=M_b\}} \over \sum_{|u|=b} {\bf 1}_{\{ V(u)=M_b \}}}-{{\bf 1}_{\{V(w_b)=M_b^{\rm kill}\}} \over \sum_{|u|^{\rm kill}=b} {\bf 1}_{\{ V(u)=M_b^{\rm kill} \}}}\Big|$ is at most $1$ and is equal to zero if no particle touched the barrier $0$, we have that, for any $H\ge 0$ such that $H\le a_n(z+L)$,
$$
\Big|{{\bf 1}_{\{V(w_b)=M_b\}} \over \sum_{|u|=b} {\bf 1}_{\{ V(u)=M_b \}}}-{{\bf 1}_{\{V(w_b)=M_b^{\rm kill}\}} \over \sum_{|u|^{\rm kill}=b} {\bf 1}_{\{ V(u)=M_b^{\rm kill} \}}}\Big| \le {\bf 1}_{\{\exists |u|\le b\,:\, V(u)\le a_n(z+L+H)\}}.
$$

\noindent Consequently,
\begin{eqnarray*}
&&\Big|F^{\rm kill}(x)- n^{3/2}\ee^{-z}F_{L,b}(x-a_n(z+L))\Big|\\
&\le&
\hat{\E}_x\left[\ee^{V(w_b)}{\bf 1}_{\{\exists |u|\le b\,:\, V(u)\le a_n(z+L+H)\}},\, \min_{k\in [0,b]} V(w_k)\ge a_n(z+L+1),\,V(w_b)\in I_n(z) \right]\\
&\le&
n^{3/2}\ee^{-z} \hat{\E}_x\left[{\bf 1}_{\{\exists |u|\le b\,:\, V(u)\le a_n(z+L+H)\}},\, \min_{k\in [0,b]} V(w_k)\ge a_n(z+L+1),\,V(w_b)\in I_n(z) \right]\\
&=&  n^{3/2}\ee^{-z}G_{H}(x-a_n(z+L))
\end{eqnarray*}

\noindent with for any $y\ge 0$,
$$
G_H(y):=\Q_y\left( \{\exists |u|\le b\,:\, V(u)\le -H\} \cap\{ \min_{k\in [0,b]} V(w_k)\ge -1,\,V(w_b)\in [L-1,L)\}\right).
$$

\noindent We do not write the dependency on $L$ and $b\ge B$ because they are fixed in this proof and so are considered as constants. This shows that, for any $z\ge 0$, $n\ge 1$ and $H\in[0,a_n(z+L)]$,
\begin{eqnarray*}
&& \hat{\E}\Big[\Big| F^{\rm kill}(V(w_{n-b}))- n^{3/2}\ee^{-z}F_{L,b}(V(w_{n-b}) - a_n(z+L))\Big|{\bf 1}_{\{ V(w_\ell) \ge d_\ell,\, \forall \, \ell\le n-b\}}\Big]\\
 &\le& n^{3/2}\ee^{-z} \hat{\E}\Big[ G_H(V(w_{n-b}) - a_n(z+L)){\bf 1}_{\{ V(w_\ell) \ge d_\ell,\, \forall \, \ell\le n-b\}}\Big].
\end{eqnarray*}

\noindent We choose $H$ such that ${C_-C_+ \sqrt{\pi} \over \sigma \sqrt{2}}\int_{y\ge 0}G_H(y)R_-(y)dy \le \eta/2$. We can check that the function $G_H$ satisfies the conditions of Lemma \ref{l:RWconvlaw} as we did for $F_{L,b}$. By Lemma \ref{l:RWconvlaw}, this yields that
$$
\hat{\E}\Big[\Big| F^{\rm kill}(V(w_{n-b}))- n^{3/2}\ee^{-z}F_{L,b}(V(w_{n-b}) - a_n(z+L))\Big|{\bf 1}_{\{ V(w_\ell) \ge d_\ell,\, \forall \, \ell\le n-b\}}\Big] \le \eta \ee^{-z}
$$

\noindent for $n$ large enough and $z\in [0,(3/2)\ln n-L-H]$. The cut-off at $(3/2)\ln(n)-L-H$ is here only to ensure that $H\le a_n(z+L)$. Combined with (\ref{eq:step2a}), we get that for $n$ large enough, and $z\in [A,(3/2)\ln(n)-L-H]$,
\begin{equation}\label{eq:step2b}
\Big|\hat{\E}_{\hbox{\scriptsize (\ref{eq:Clbh})}} -  n^{3/2}\ee^{-z} \hat{\E}\Big[F_{L,b}(V(w_{n-b}) - a_n(z+L)),\, V(w_\ell) \ge d_\ell,\, \forall \,0\le \ell \le n-b\Big]\Big|
\le 2\eta \ee^{-z}.
\end{equation}

\noindent  Recall the definition of $C_{L,b}$ in (\ref{def:CLb}). We apply again Lemma \ref{l:RWconvlaw} to see that
$$
\hat{\E}\Big[F_{L,b}(V(w_{n-b}) - a_n(z+L)),\, V(w_\ell) \ge d_\ell,\, \forall \,0\le \ell \le n-b\Big] \sim { C_{L,b} \over n^{3/2}}
$$

\noindent as $n\to\infty$ uniformly in $z\in[0, (3/2)\ln(n)-L]$. Consequently, we have for $n$ large enough and $z\in[0, (3/2)\ln(n)-L]$,
$$
\Big|n^{3/2}\ee^{-z}\hat{\E}\Big[F_{L,b}(V(w_{n-b}) - a_n(z+L)),\, V(w_\ell) \ge d_\ell,\, \forall \,0\le \ell \le n-b\Big]  - \ee^{-z}C_{L,b}\Big| \le \eta\ee^{-z}.
$$

\noindent The lemma follows from  (\ref{eq:step2b}). \hfill $\Box$

\bigskip

We now have the tools to prove Proposition \ref{p:tailminkill}.\\

\noindent {\it Proof of Proposition \ref{p:tailminkill}.} Let $\hat{\E}_{\hbox{\scriptsize (\ref{eq:Clbh})}}$ be the expectation in the left-hand side of (\ref{eq:Clbh}). We introduce for any $L\ge 0$ and any integer $b\ge 1$,
\begin{eqnarray*}
C_{L,b}^- &:=& \liminf_{z\to\infty} \liminf_{n\to\infty} \ee^z \hat{\E}_{\hbox{\scriptsize (\ref{eq:Clbh})}},\\
C_{L,b}^+ &:=&  \limsup_{z\to\infty} \limsup_{n\to\infty} \ee^z\hat{\E}_{\hbox{\scriptsize (\ref{eq:Clbh})}}.
\end{eqnarray*}

\noindent In particular, taking the limits in $n\to\infty$ then $z\to\infty$ in (\ref{eq:Clbh}), we have, for any $L\ge 0$, $\eta>0$ and $b\ge B(L,\eta)$ (with $B(L,\eta)$ as in Lemma \ref{l:En}),
\begin{equation}\label{eq:Clbeta}
C_{L,b} - 3\eta \le C_{L,b}^- \le C_{L,b}^+ \le C_{L,b} +3\eta.
\end{equation}

\noindent Notice that $\mathcal E_n(z,b)$ (hence $\hat{\E}_{\hbox{\scriptsize (\ref{eq:Clbh})}}$) is increasing in $b$. This implies that $C_{L,b}^-$ and $C_{L,b}^+$ are both increasing in $b$. For any $L\ge 0$, let $C_{L}^-$ and $C_{L}^+$ be respectively the (possibly zero or infinite) limits of $C_{L,b}^{-}$ and $C_{L,b}^+$ when $b\to\infty$.  By (\ref{eq:Clbeta}), we have for any $L\ge 0$ and $\eta>0$,
$$
\limsup_{b\to\infty} C_{L,b} - 3\eta \le C_{L}^- \le C_{L}^+ \le \liminf_{b\to\infty}C_{L,b} +3\eta.
$$

\noindent Letting $\eta$ go to $0$, this yields that $C_{L,b}$ has a limit as $b\to\infty$, that we denote by $C(L)=C_L^-=C_L^+$, this for any $L\ge 0$. Similarly, we see that  $\hat{\E}_{\hbox{\scriptsize (\ref{eq:Clbh})}}$ is increasing in $L$. This gives that $C(L)$ admits a limit as $L\to \infty$, that we denote by $C_2$. Beware that at this stage, we do not know whether $C_2\in(0,\infty)$. Let $\varepsilon>0$. By (\ref{cor:tightL}), there exists $L_0\ge 0$ such that for any $L\ge L_0$, $z\ge 0$ and $n\ge 1$,
$$
\P(m^{{\rm kill},(n)} \notin \mathcal Z^{z,L}_n,\, M_n^{\rm kill} \in I_n(z) )\le \varepsilon \ee^{-z}.
$$

\noindent By Lemma \ref{l:step1}, this yields that for $L\ge L_0$, $z\ge 0$ and any $n\ge 1$,
$$
\Big|\P(M_n^{\rm kill} \in I_n(z) )- \hat{\E}\left[{\ee^{V(w_n)}{\bf 1}_{\{ V(w_n)=M_n^{\rm kill}\}}\over \sum_{|u|^{\rm kill}=n} {\bf 1}_{\{ V(u)=M_n^{\rm kill}\}}},\, w_n \in \mathcal Z^{z,L}_n\right] \Big| \le \varepsilon \ee^{-z}.
$$

\noindent Take again $\eta>0$ and $L\ge L_0$,  and let $B=B(L,\eta)\ge 1$  and $A=A(L,\eta)>0$ as in Lemma \ref{l:En}. We have 
$$
\hat{\E}\left[{\ee^{V(w_n)}{\bf 1}_{\{ V(w_n)=M_n^{\rm kill}\}}\over \sum_{|u|^{\rm kill}=n} {\bf 1}_{\{ V(u)=M_n^{\rm kill}\}}},\, w_n \in \mathcal Z^{z,L}_n,\,\mathcal E_n(z,b)^c\right] \le n^{3/2}\ee^{-z} \Q(w_n \in \mathcal Z^{z,L}_n,\,\mathcal E_n(z,b)^c)\le \eta\ee^{-z}
$$

\noindent for any $n\ge b\ge B$ and $z\ge A$. Consequently, for any $L\ge L_0$, $n\ge b\ge B$ and $z\ge A$,
$$
\Big|\P(M_n^{\rm kill} \in I_n(z) )- \hat{\E}\left[{\ee^{V(w_n)}{\bf 1}_{\{ V(w_n)=M_n^{\rm kill}\}}\over \sum_{|u|^{\rm kill}=n} {\bf 1}_{\{ V(u)=M_n^{\rm kill}\}}},\, w_n \in \mathcal Z^{z,L}_n,\,\mathcal E_n(z,b)\right] \Big| \le (\varepsilon+\eta)\ee^{-z}.
$$

\noindent By Lemma \ref{l:step2}, we get that for $L\ge L_0$, $b\ge B(L,\eta)$,  $n$ large enough and $z\in[A(L,\eta),(3/2)\ln(n)-L-H(L,\eta,b)]$,
\begin{equation}\label{eq:boundCLb}
\Big|\ee^{z}\P(M_n^{\rm kill} \in I_n(z) )- C_{L,b} \Big| \le (\varepsilon+4\eta).
\end{equation}

\noindent We stress that $C_{L,b}$ depends actually on $\eta$ and $\varepsilon$ through the choice of $L_0$ and $B(L,\eta)$. By (\ref{eq:boundCLb}) and Corollary \ref{cor:tightnesskill-upper}, we know that for $L\ge L_0$ and $b\ge B(L,\eta)$, we have  $C_{L,b} \le c_{10} + \varepsilon+4\eta$. Taking the limit $b\to\infty$, this implies that for any $L\ge L_0$, we have $C(L)\le c_{10}+\varepsilon+4\eta$. Taking the limit $L\to\infty$, we deduce that $C_2\le c_{10}+\varepsilon+4\eta$ hence $C_2$ is finite. Let $L>L_0$ such that $|C_2-C(L)|\le \eta$ and $b\ge B(L)$ such that $|C_{L,b}-C(L)|\le \eta$. Then, by (\ref{eq:boundCLb}), we have for $n$ large enough and $z\in[A(L,\eta),(3/2)\ln(n)-L-H(L,\eta,b)]$,
$$
\Big|\ee^z \P(M_n^{\rm kill} \in I_n(z)) - C_2  \Big| \le \varepsilon +6\eta\le 2 \varepsilon 
$$

\noindent if we take $\eta:=\varepsilon/6$. It remains to show that $C_2>0$. We see that, necessarily,
$$
\limsup_{z\to\infty}\limsup_{n\to\infty} \Big|\ee^z \P(M_n^{\rm kill} < \frac32 \log n -z) - {C_2\over 1-\ee^{-1}}  \Big|=0.
$$

\noindent We know then that $C_2> 0$ by the lower bound obtained in Lemma \ref{l:tightnesskill}. \hfill $\Box$

\subsection{Proof of Lemma \ref{l:En}}
\label{s:proofEn}

We present here the postponed proof of Lemma \ref{l:En}.

\bigskip

\noindent {\it Proof of Lemma \ref{l:En}}.  We follow the same strategy as for Lemma \ref{l:tightnesskill}. Let $\eta>0$. To avoid superfluous notation, we prove the lemma for $L=0$ (the general case works similarly). Recall the definition of $\mathcal E_n(z,b)$ in (\ref{def:En}). We want to show that $\Q(\mathcal E_n(z,b)^c, \,w_n\in\mathcal Z^{z,0}_n)\le \eta n^{-3/2}$ when $b$ and $z$ are large enough. Let $d_k=d_k(n,z,1/2)$ as defined in (\ref{def:dk}) and
\begin{equation}\label{def:ek}
e_k=e_k^{(n)}:=
\begin{cases}
k^{1/12}, &\hbox{if $0\le k\le {n\over 2}$,} \cr
(n-k)^{1/12}, &\hbox{if ${n\over 2}<k\le n$.}\cr
\end{cases}
\end{equation}

\noindent We recall that $|u|=n$ is a $z$-good vertex if $u\in \mathcal Z^{z,0}_n$ and
$$
\sum_{v\in \Omega(u_{k})}\ee^{-(V(v)-d_k)}\Big\{1+(V(v)-d_k)_+\Big\} \le B\ee^{-e_k} \qquad  \forall\, 1\le k\le n
$$

\noindent with $B$ such that, for $n\ge 1$ and $z\ge 0$,
\begin{equation}\label{eq:goodvertexl}
\Q(w_n \in \mathcal Z^{z,0}_n,\, w_n \mbox{ is not a $z$-good vertex}) \le {\eta \over n^{3/2}}
\end{equation}

\noindent (see Lemma \ref{l:goodvertex}). Recall that $\Omega(w_k)$ is the set of siblings of $w_k$ and $\hat{\mathcal G}_\infty$ is defined in (\ref{def:Ginfty}). Recall the law of the branching random walk under $\Q$ which we identified with $\hat{\mathcal B}$ by the mean of Proposition \ref{p:spine1}. For $\mathcal E_n(z,b)$ to happen, every sibling of the spine at generation less than $n-b$ must have all its descendants at time $n$ at position greater than $a_n(z)$. In other words,
\begin{eqnarray}\label{eq:puz}
&& \Q((\mathcal E_n(z,b))^c,\, w_n \mbox{ is a $z$-good vertex})\\
\nonumber &=&
\hat{\E}\left[1-\prod_{k=1}^{n-b} \prod_{u\in \Omega(w_k)} (1-\Phi^{\rm kill}_{k,n}(V(u),z)),\, w_n \mbox{ is a $z$-good vertex}\right]
\end{eqnarray}

\noindent where $\Phi^{\rm kill}_{k,n}(V(u),z):= \P_{V(u)}(M_{n-k}^{\rm kill} < a_n(z))$ is the probability that the killed branching random walk rooted at $u$ has its minimum greater than $a_n(z)$ at time $n-k$. By Corollary \ref{cor:tightnesskill-upper}, we see that, if $|u|\le n/2$ (hence $a_{n}(z)=a_{n-|u|}(z)+O(1))$, then
$$
\Phi^{\rm kill}_{k,n}(V(u),z)
\le c_{20} (1+V(u)_+)\ee^{-z-V(u)}.
$$

\noindent On the event that $w_n$ is a $z$-good vertex, we have for $k\le n/2$ (hence $d_k =0$), $\sum_{u\in \Omega(w_k)}(1+V(u)_+)\ee^{-V(u)} \le B\ee^{-e_k}=B\ee^{-k^{1/12}}$. Using the inequality $x\ge \ee^{(x-1)/2}$ for $x$ close enough to $1$, we deduce that there exists $A_0\ge 0$ such that for $z\ge A_0$, $n\ge 1$, and $1\le k\le n/2$, on the event that $w_n$ is a $z$-good vertex, we have
$$
\prod_{u\in \Omega(w_k)} (1-\Phi^{\rm kill}_{k,n}(V(u),z)) \ge \exp\left(-c_{21} \ee^{-z}\ee^{-k^{1/12}}\right)
$$

\noindent with $c_{21}:=c_{20}B/2$. This yields that
$$
\prod_{k=1}^{\lfloor n/2 \rfloor} \prod_{u\in \Omega(w_k)} (1-\Phi^{\rm kill}_{k,n}(V(u),z)) \ge \exp\left(-c_{21} \ee^{-z}\sum_{k=1}^{\lfloor n/2 \rfloor}\ee^{-k^{1/12}}\right)\ge \exp(-c_{22}\ee^{-z}).
$$

\noindent Therefore, there exists $A_1>A_0$ such that for any $z\ge A_1$ and $n\ge 1$,
\begin{equation}\label{eq:En1}
\prod_{k=1}^{\lfloor n/2\rfloor} \prod_{u\in \Omega(w_k)} (1-\Phi^{\rm kill}_{k,n}(V(u),z)) \ge (1-\eta)^{1/2}.
\end{equation}

\noindent If $k> n/2$, we simply observe that if $M_{\ell}^{\rm kill}\le x$, a fortiori $M_\ell\le x$.  Since $W_n$ (defined in (\ref{def:Wn})) is a martingale, we have $1=\E[W_\ell]\ge \E[\ee^{-M_{\ell}}] \ge \ee^{-x}\P(M_\ell\le x)$ for any $\ell\ge 1$ and $x\in \r$. We get that
$$
\Phi^{\rm kill}_{k,n}(V(u),z)\le \P(M_{n-|u|} < a_n(z)-V(u))\le \ee^{ a_n(z)}\ee^{-V(u)} .
$$

\noindent We rewrite it  $\Phi^{\rm kill}_{k,n}(V(u),z) \le    \ee^{-(V(u) - d_k)}$ for $n/2<k\le n$. On the event that $w_n$ is a $z$-good vertex, we get that $\prod_{u\in \Omega(w_k)} (1-\Phi^{\rm kill}_{k,n}(V(u),z)) \ge \ee^{- c_{23} \ee^{-e_k}}= \ee^{- c_{23} (n-k)^{1/12}}$ for $k$ greater than some constant $b_1$. Consequently, for any $b\ge b_1$,
$$
\prod_{k=\lfloor n/2\rfloor +1}^{n-b} \prod_{u\in \Omega(w_k)} (1-\Phi^{\rm kill}_{k,n}(V(u),z)) \ge \ee^{-c_{23} \sum_{k=\lfloor n/2\rfloor +1}^{n-b} \ee^{-(n-k)^{1/12}}}.
$$

\noindent This yields that there exists  $B\ge 1$ such that for any $b\ge B$ and any $n\ge 1$, we have,
\begin{equation}\label{eq:En2}
\prod_{k=\lfloor n/2\rfloor +1}^{n-b} \prod_{u\in \Omega(w_k)} (1-\Phi^{\rm kill}_{k,n}(V(u),z)) \ge (1-\eta)^{1/2}.
\end{equation}

\noindent In view of (\ref{eq:En1}) and (\ref{eq:En2}), we have for $b\ge B$, $z\ge A_1$ and $n\ge 1$, $\prod_{k=1}^{n-b} \prod_{u\in \Omega(w_k)}  (1-\Phi^{\rm kill}_{k,n}(V(u),z)) \ge (1-\eta)$. Plugging it into (\ref{eq:puz}) yields that
$$
\Q((\mathcal E_n(z,b))^c,\, w_n \mbox{ is a $z$-good vertex})
\le
\eta \Q\left(w_n \mbox{ is a $z$-good vertex	}\right)\le \eta\Q\left(w_n \in \mathcal Z^{z,0}_n\right) .
$$

\noindent It follows from (\ref{eq:goodvertexl}) that
$$
\Q((\mathcal E_n(z,b))^c,\, w_n \in \mathcal Z^{z,0}_n)
\le
 \eta (\Q\left(w_n \in \mathcal Z^{z,0}_n\right) + n^{-3/2}).
$$

\noindent Recall that the spine behaves as a centered random walk. Then apply (\ref{lemmaA3}) to see that $\Q\left(w_n \in \mathcal Z^{z,0}_n\right)\le c_{24} n^{-3/2}$, which completes the proof of the lemma.  \hfill $\Box$

\section{Tail distribution of the minimum of the BRW}

\label{s:tailmin}
We prove a slightly stronger version of Proposition \ref{p:tailminmain}.
\begin{proposition}\label{p:tailmin}
Let $C_1$ be as in Proposition \ref{p:tailminkillmain} and $c_0$ as in (\ref{c0}). For any $\varepsilon>0$, there exist $N\ge 1$ and $A>0$ such that for any $n\ge N$ and $z\in[A,(3/2)\ln n -A]$,
$$
\Big| {\ee^z \over z}\P(M_n < \frac32 \ln n -z) - C_1c_0\Big| \le \varepsilon.
$$
\end{proposition}

\bigskip

We introduce some notation. To go from the tail distribution of  $M_n^{\rm kill}$ to the one of $M_n$, we have to control excursions inside the negative axis that can appear at the beginning of the branching random walk. For any real $r$, we define the set
\begin{equation}\label{def:SA}
\mathcal S^{\,r} :=\{ u \in \T \,:\,   \min_{k \le |u|-1} V(u_k)> V(u) \ge -r\}.
\end{equation}

\noindent Notice that $\mathcal S^{\,r}=\emptyset$ when $r< 0$. Let for $|v|\ge 1$,
 \begin{equation}\label{def:xi}
 \xi(v) := \sum_{w\in \Omega(v)} (1+(V(w)-V({\buildrel \leftarrow \over v}))_+)\ee^{-(V(w)-V({\buildrel \leftarrow \over v}))}
 \end{equation}

\noindent where ${\buildrel \leftarrow \over v}$ denotes the parent of $v$ (and $y_+:=\max(y,0)$). Notice that $\xi(v)$ is stochastically smaller than  $X+\tilde X$ as defined in (\ref{def:X}). To avoid some extra integrability conditions, we are led to consider vertices $u\in \mathcal S^{\,r}$ which behave 'nicely', meaning that $\xi(u_k)$ is not too big along the path $\{u_1,\ldots,u_{|u|}=u\}$.  Hence, for any real $r \ge 0$, we introduce 
\begin{equation}\label{def:T} 
\mathcal T^{\,r} := \{u\in \T\,:\, \forall 1\le k\le |u|: \xi(u_k) < \ee^{(V(u_{k-1})+r )/2}\}.
\end{equation}

\noindent For any integer $k\ge 0$, we denote by $\mathcal S^{\,r}_k$, resp. $\mathcal T^{\,r}_k$, the set $\mathcal S^{r} \cap\{|u|=k\}$, resp.  $\mathcal T^{\,r}\cap\{|u|=k\}$. Finally, for any integer $n\ge 1$, any $z\ge 0$ and any $u\in\T$, define
\begin{equation}\label{def:Bn}
B_n^z(u) :=
\begin{cases}
1 & \; \mbox{if }\exists\, v\ge u\,:\, |v|=n,\,\min_{ \ell \in[|u|,n]}V(v_\ell)\ge V(u), \mbox{ and } V(v)< a_n(z),\\ 
0 & \mbox{ otherwise}.
\end{cases}
\end{equation}

\noindent Notice that $B_n^z(u)=0$ if $|u|>n$. In words, $B_n^z(u)=1$ if there exists a descendant of $u$ which stays above $V(u)$ and is below level $a_n(z)$ at time $n$. Observe that if $M_n< a_n(z)$, then necessarily we can find such vertices $u$ and $v$.  The first subsection controls the set $\mathcal S^{r}$. Proposition \ref{p:tailminmain} is then proved in Section \ref{s:tailminmain}.

\begin{figure}
\begin{center}
\resizebox{8 cm}{4 cm}{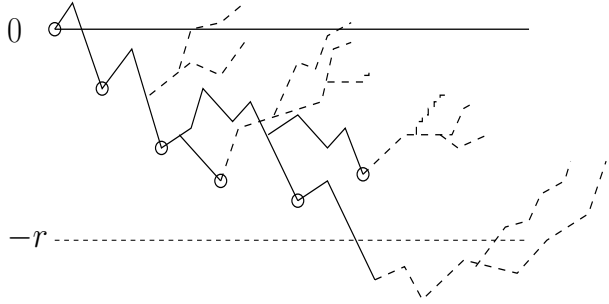}
\caption{ The set $\mathcal S^r$}
\label{f:Sr}
\end{center}
\end{figure}

\subsection{The branching random walk at the beginning}

 We will see that $\P(M_n < \frac32\ln n -z)$ is comparable to the probability that there exists $u \in \mathcal S^{z}$ such that $B_n^z(u)=1$. The lemmas in this section are used to give an asymptotic of this probability. As usual, we will use a second moment argument. Lemmas \ref{l:moment1} and \ref{l:moment2} give bounds  respectively on the first moment and second moment of the number of such vertices $u$. We recall that $M_{n}^{\rm kill}$ is the minimum at time $n$ of the branching random walk killed below zero. For any  integers $n\ge 1$, $k\in[0,n]$, and any reals $x,r$, we recall that
  \begin{equation}\label{def:Phikn}
 \Phi^{\rm kill}_{k,n}(x,r):= \P_x\left( M_{n-k}^{\rm kill} < a_n(r)\right).
\end{equation}

\noindent By Corollary \ref{c:tailminkill}, there exists  $N_0\ge 1$ and $A_0\ge 0$ such that for any $n\ge N_0$, $k\le n^{1/2}$ and $r\in[A_0,(3/2)\ln(n)-A_0]$, 
\begin{equation}\label{eq:phiC1}
\Big|\ee^{r}\Phi^{\rm kill}_{k,n}(0,r) -C_1 \Big|
\le
\varepsilon
\end{equation}

\noindent where we used the fact that $k=o(n)$ thus $\ln(n-k)=\ln(n)+o(1)$ (the same statement holds when replacing $n^{1/2}$ by any sequence $o(n)$).  Moreover, we know by Corollary \ref{cor:tightnesskill-upper} that for any integers $n\ge 1$, $k\in[0,n]$, and any reals $x,r\ge 0$,
\begin{equation}\label{eq:boundphikn}
 \Phi^{\rm kill}_{k,n}(x,r)\le c_{25} (1+x)\ee^{-x-r} \left( {n\over n-k+1} \right)^{3/2}.
\end{equation}

\begin{lemma}\label{l:moment1}
(i) Fix $\varepsilon>0$ and let $C_1$ be the constant in Proposition \ref{p:tailminkillmain}. There exists $A\ge 0$ such that for all $n$ sufficiently large, and all $z\in [A,(3/2)\ln(n)-A]$,
\begin{equation} \label{moment1}
\Bigg|{\ee^z \over R(z-A)} \E\left[  \sum_{u\in \mathcal S^{z-A}} B_n^z(u){\bf 1}_{\{ |u|\le  n^{1/2}\}} \right]
-
C_1 \Bigg| \le \varepsilon.
\end{equation}
(ii) There exists a constant $c$ such that for any $n\ge 1$ and any $z\in [0,(3/2)\ln(n)]$,
$$
 \E\left[  \sum_{u\in \mathcal S^{z}} B_n^z(u){\bf 1}_{\{ |u|>  n^{1/2}\}} \right]
\le c \ee^{-z}.
$$
(iii)  Uniformly in $A\ge 0$ and $n\ge 1$, we have
$$
\E\left[ \sum_{u\in \mathcal S^{z-A}\cap(\mathcal T^{z-A})^c} B_n^z(u){\bf 1}_{\{|u|\le n/2\}} \right] =  o(z)\ee^{-z}
$$

\noindent as $z\to\infty$, where the set $(\mathcal T^{z-A})^c$ denotes the complement of the set $\mathcal T^{z-A}$ in the set of vertices of $\T$.
\end{lemma}

\noindent {\bf Remark}. In (i) and (ii), we could replace $n^{1/2}$ by $n^a$ with $a\in (0,1)$. \\

\begin{figure}
\begin{center}
\resizebox{10 cm}{6 cm}{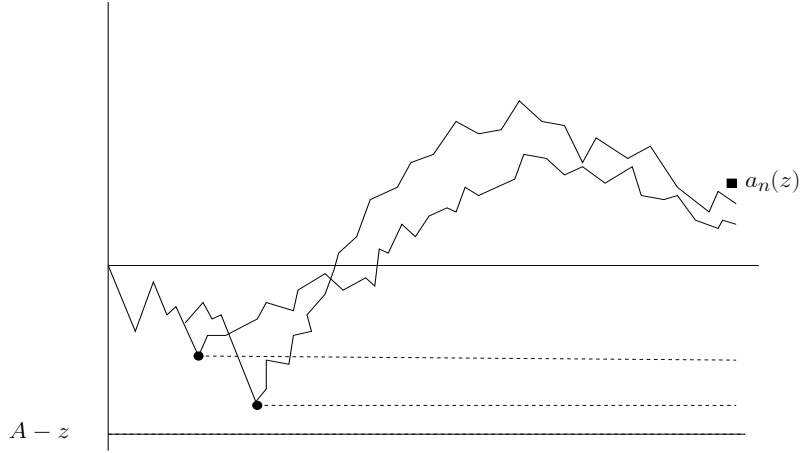}
\caption{ Particles in $\mathcal S^{z-A}$ such that $B_n^z(u)=1$.}
\label{f:Bnz}
\end{center}
\end{figure}

\noindent{\it Proof}. Let $k\le n$. By the Markov property at time $k$, we have
\begin{equation}\label{eq:momentsnl}
\E\left[ \sum_{u\in \mathcal S^{z-A}_k} B_n^z(u) \right]
=
 \E\left[ \sum_{u\in \mathcal S_k^{z-A}} \Phi^{\rm kill}_{k,n}(0,z+V(u))\right]
 \end{equation}
 
\noindent  with $\Phi^{\rm kill}_{k,n}$ as defined in (\ref{def:Phikn}). We want to apply equation (\ref{eq:phiC1}) to $r=z+V(u)$. We observe that $z+V(u) \in [A,z]$ when $u\in \mathcal S^{z-A}$. Hence, equation (\ref{eq:phiC1}) holds for $n\ge N_0$, $k\le n^{1/2}$, and $r=z+V(u)$, with $u\in \mathcal S^{z-A_0}$ and $z\in [A_0,(3/2)\ln(n)-A_0]$. It follows from (\ref{eq:momentsnl}) that for $n\ge N_0$, $k\le n^{1/2}$  and $z\in [A_0,(3/2)\ln(n)-A_0]$,
\begin{equation}\label{eq:moment1S(z-A)a}
\Bigg|\ee^z\E\left[ \sum_{u\in \mathcal S^{z-A_0}_k} B_n^z(u) \right] - C_1\E\left[ \sum_{u\in \mathcal S^{z-A_0}_k}\ee^{-V(u)}\right]\Bigg|
\le
\varepsilon \E\left[ \sum_{u\in \mathcal S_k^{z-A_0}}\ee^{-V(u)}\right].
\end{equation}

\noindent From the definition of $\mathcal S_k^{z-A}$ and (\ref{many-to-one}), we observe that, for any integer $k$, and any $z\ge A\ge 0$,
\begin{equation}\label{eq:manyto1Sk}
\E\left[ \sum_{u\in \mathcal S_k^{z-A}}\ee^{-V(u)}\right]= \P(S_k \ge A-z,\,S_k< S_\ell,\,\forall \,0\le \ell <k-1). 
\end{equation}

\noindent Summing over $k\ge 0$ yields  that
\begin{equation}\label{eq:manyto1S(z-A)}
\E\left[ \sum_{u\in \mathcal S^{z-A}}\ee^{-V(u)}\right]= R(z-A).
\end{equation}

\noindent In particular, summing equation (\ref{eq:moment1S(z-A)a}) over $k\le n^{1/2}$ gives that  for $n\ge N_0$  and $z\in [A_0,(3/2)\ln(n)-A_0]$,
\begin{equation}\label{eq:moment1S(z-A)b}
\Bigg|\ee^z\E\left[ \sum_{u\in \mathcal S^{z-A_0}} B_n^z(u){\bf 1}_{\{|u|\le n^{1/2}\}} \right] - C_1\E\left[ \sum_{u\in \mathcal S^{z-A_0}_k} \ee^{-V(u)}{\bf 1}_{\{|u|\le n^{1/2}\}}\right]\Bigg|
\le
\varepsilon R(z-A_0).
\end{equation}

\noindent Using the fact that $\P(S_k\ge -x,\, S_k < \min_{0\le j \le k-1} S_j)= \P((-S_k)\le x,\,  \min_{0\le j \le k-1} (-S_j)\ge 0)$, we have   by (\ref{lemmaA1}), for any integer $k\ge 0$ and any real $x\ge 0$, 
 \begin{equation}\label{eq:R(x)tail}
 \P(S_k\ge -x,\, S_k < \min_{0\le j \le k-1} S_j)  \le \alpha_2' (1+x)^2(1+k)^{-3/2}.
  \end{equation}

\noindent Therefore, we have for $n$ greater than some $N_1$ and $x \in[0,(3/2)\ln(n)]$,
$$
\sum_{k>n^{1/2}} \P(S_k \ge -x,\,S_k< S_\ell,\,\forall \,0\le \ell <k-1) \le \varepsilon.
$$

\noindent Going back to (\ref{eq:manyto1Sk}) with $A=A_0$, and summing over $k>n^{1/2}$, we obtain that for $n\ge N_1$ and $z\in[A_0,(3/2)\ln(n)-A_0]$,
$$
\E\left[ \sum_{u\in \mathcal S^{z-A_0}}\ee^{-V(u)} {\bf 1}_{\{|u| > n^{1/2} \}} \right]\le \varepsilon.
$$

\noindent In view of (\ref{eq:moment1S(z-A)b}) and (\ref{eq:manyto1S(z-A)}), this yields that for any $n\ge \max(N_0,N_1)$ and any $z\in[A_0,(3/2)\ln(n)-A_0]$,
$$
\Bigg|\ee^z\E\left[ \sum_{u\in \mathcal S^{z-A_0}} B_n^z(u){\bf 1}_{\{|u|\le n^{1/2}\}} \right] - C_1R(z-A_0)\Bigg|
\le
\varepsilon (R(z-A_0)+C_1).
$$

\noindent Since $R(x)\ge 1$ for any $x\ge 0$, this completes the proof of (i).  Let us prove (ii). The notation $c$ denotes a constant whose value can change from line to line. Using (\ref{eq:momentsnl}) with $A=0$, we find that
\begin{equation}
\E\left[ \sum_{u\in \mathcal S^{z}} B_n^z(u){\bf 1}_{\{ n/2\ge |u|> n^{1/2} \}} \right]
=
\sum_{k=\lfloor n^{1/2}\rfloor+1}^{\lfloor n/2\rfloor} \E\left[ \sum_{u\in \mathcal S_k^{z}} \Phi^{\rm kill}_{k,n}(0,z+V(u))\right].
\end{equation}

\noindent Equation (\ref{eq:boundphikn}) yields that
\begin{equation}
\E\left[ \sum_{u\in \mathcal S^{z}} B_n^z(u){\bf 1}_{\{ n/2\ge |u|> n^{1/2} \}} \right]
\le 
c_{25} \ee^{-z}  \sum_{k=\lfloor n^{1/2}\rfloor+1}^{\lfloor n/2\rfloor} \left( {n\over n-k+1} \right)^{3/2}  \E\left[ \sum_{u\in \mathcal S_k^{z}}  \ee^{-V(u)}   \right].
\end{equation}

\noindent Equations (\ref{eq:manyto1Sk}) and (\ref{eq:R(x)tail}) imply that
\begin{eqnarray} \nonumber 
\E\left[ \sum_{u\in \mathcal S^{z}} B_n^z(u){\bf 1}_{\{ n/2\ge |u|> n^{1/2} \}} \right]
&\le& 
c_{25} \alpha_2' \ee^{-z}(1+z)^2  \sum_{k=\lfloor n^{1/2}\rfloor+1}^{\lfloor n/2\rfloor} \left( {n\over n-k+1} \right)^{3/2} (1+k)^{-3/2} \\
&\le& c \ee^{-z} \label{eq:zhan0}
\end{eqnarray}

\noindent for any $n\ge 1$ and any $z\in [0,(3/2)\ln(n)]$. We deal now with vertices $u\in \mathcal S^z$ such that $|u|>n/2$, and split the case depending on whether $V(u)$ is greater or smaller than $-z +{3\over 2} \ln\left( {n\over n-|u|+1}\right)$.  Using the fact that $B_n^z(u)\le 1$, we get
\begin{eqnarray}\label{eq:zhan}
&& \E\left[ \sum_{u\in \mathcal S^{z}} B_n^z(u){\bf 1}_{\{  |u|> n/2 \}} \right] \\
\nonumber  &\le&
\E\left[ \sum_{u\in \mathcal S^{z}} B_n^z(u){\bf 1}_{\{  |u|> n/2,\, V(u)\ge -z +{3\over 2} \ln\left( {n\over n-|u|+1}\right) \}} \right]
+
\E\left[ \sum_{u\in \mathcal S^{z}} {\bf 1}_{\{  |u|> n/2,\, V(u)< -z +{3\over 2} \ln\left( {n\over n-|u|+1}\right) \}} \right] .
\end{eqnarray}
We bound the first term of the right-hand side. Equation (\ref{eq:boundphikn}) shows that
\begin{eqnarray*}
&&\E\left[ \sum_{u\in \mathcal S^{z}} B_n^z(u){\bf 1}_{\{  |u|> n/2,\, V(u)\ge -z +{3\over 2} \ln\left( {n\over n-|u|+1}\right) \}} \right] \\
&\le& 
c_{25} \ee^{-z}  \sum_{k=\lfloor n/2\rfloor+1}^{n} \left( {n\over n-k+1} \right)^{3/2}  \E\left[ \sum_{u\in \mathcal S_k^{z}}  \ee^{-V(u)} {\bf 1}_{\{   V(u)\ge -z +{3\over 2} \ln\left( {n\over n-k+1}\right) \}}   \right].
\end{eqnarray*}

\noindent From  (\ref{many-to-one}), we observe that,
\begin{eqnarray*}
&& \E\left[ \sum_{u\in \mathcal S_k^{z}}\ee^{-V(u)} {\bf 1}_{\{   V(u)\ge -z +{3\over 2} \ln\left( {n\over n-k+1}\right) \}}   \right]\\
&=&  \P\left(S_k \ge -z + {3\over 2} \ln\left( {n\over n-k+1}\right),\,S_k< S_\ell,\,\forall \,0\le \ell <k-1\right)
\end{eqnarray*}
which is $0$ if $-z + {3\over 2} \ln\left( {n\over n-k+1}\right)\ge 0$. Using (\ref{eq:R(x)tail}), we see that
\begin{eqnarray*}
&& \E\left[ \sum_{u\in \mathcal S^{z}} B_n^z(u){\bf 1}_{\{  |u|> n/2,\, V(u)\ge -z + {3\over 2} \ln\left( {n\over n-|u|+1}\right) \}} \right]  \\
&\le&
c \ee^{-z}  \sum_{k=\lfloor n/2\rfloor+1}^{n} \left( {n\over n-k+1} \right)^{3/2}  (1+k)^{-3/2} \max\left(z - {3\over 2} \ln\left( {n\over n-k+1}\right),1\right)^2.
\end{eqnarray*}
Since $z\le {3\over 2}\ln(n)$, we get that $\max\left(z - {3\over 2} \ln\left( {n\over n-k+1}\right),1\right) \le 1+{3\over 2}\ln\left(n-k+1\right)$. Consequently,
\begin{eqnarray*}
&& \E\left[ \sum_{u\in \mathcal S^{z}} B_n^z(u){\bf 1}_{\{  |u|> n/2,\, V(u)\ge -z + {3\over 2} \ln\left( {n\over n-|u|+1}\right) \}} \right]  \\
&\le&
c \ee^{-z}  \sum_{k=\lfloor n/2\rfloor+1}^{n} \left( {1\over n-k+1} \right)^{3/2}  (1+\ln(n-k+1))^2 \le c \ee^{-z}.
\end{eqnarray*}
Finally, let us consider the last term of (\ref{eq:zhan}). Equation (\ref{many-to-one}) implies that, for any $k$,
\begin{eqnarray}\nonumber
 \E\left[ \sum_{u\in \mathcal S^{z}_k} {\bf 1}_{\{  |u|> n/2,\, V(u)< -z +{3\over 2} \ln\left( {n\over n-|u|+1}\right) \}} \right] 
&=&\nonumber
 \E\left[\ee^{S_k} {\bf 1}_{\{S_k \in[-z, -z + {3\over 2} \ln\left( {n\over n-k+1}\right)),\,S_k< S_\ell,\,\forall \,0\le \ell <k-1\}}\right]\\
 &\le&\label{eq:zhan2}
  \E\left[\ee^{S_k} {\bf 1}_{\{S_k <-z + {3\over 2} \ln\left( {n\over n-k+1}\right),\,S_k< S_\ell,\,\forall \,0\le \ell <k-1\}}\right].
  \end{eqnarray}
 
 \noindent We notice that
 \begin{eqnarray*}
 && \E\left[\ee^{S_k} {\bf 1}_{\{S_k <-z + {3\over 2} \ln\left( {n\over n-k+1}\right),\,S_k< S_\ell,\,\forall \,0\le \ell <k-1\}}\right]
 \\
  &\le&
  \sum_{y\ge 0} \ee^{-y+1} \P(S_k \ge -y,\,S_k< S_\ell,\,\forall \,0\le \ell <k-1  ){\bf 1}_{\{ y> z - {3\over 2} \ln\left( {n\over n-k+1}\right)\}},
\end{eqnarray*}
which, in view of (\ref{eq:R(x)tail}) leads to
 \begin{eqnarray*}
 && \E\left[\ee^{S_k} {\bf 1}_{\{S_k <-z + {3\over 2} \ln\left( {n\over n-k+1}\right),\,S_k< S_\ell,\,\forall \,0\le \ell <k-1\}}\right]
 \\
  &\le&
 \alpha'_2 \sum_{y\ge 0} \ee^{-y+1} (1+y)^2(1+k )^{-3/2}{\bf 1}_{\{ y> z - {3\over 2} \ln\left( {n\over n-k+1}\right)\}}\\
 &\le&
  c \ee^{-z} \left({n\over n-k+1}\right)^{3/2} \max\left(z - {3\over 2} \ln\left( {n\over n-k+1}\right),1\right)^2 (1+k)^{-3/2}\\
  &\le&
   c \ee^{-z} \left({1\over n-k+1}\right)^{3/2} \left(1+ {3\over 2} \ln\left(n-k+1\right)\right)^2 
\end{eqnarray*}
for $k\in[n/2,n]$ and $z\le {3\over 2} \ln(n)$. Going back to (\ref{eq:zhan2}), it yields that
\begin{eqnarray*}
 &&\E\left[ \sum_{u\in \mathcal S^{z}} {\bf 1}_{\{  |u|> n/2,\, V(u)< -z +{3\over 2} \ln\left( {n\over n-|u|+1}\right) \}} \right] \\
 &\le&
  c \ee^{-z} \sum_{k=\lfloor n/2 \rfloor +1}^n \left({1\over n-k+1}\right)^{3/2} \left(1+ {3\over 2} \ln\left(n-k+1\right)\right)^2 \\
  &\le& c \ee^{-z}.
\end{eqnarray*}
Finally, by (\ref{eq:zhan}), $ \E\left[ \sum_{u\in \mathcal S^{z}} B_n^z(u){\bf 1}_{\{  |u|> n/2 \}} \right] \le c \ee^{-z}$
 which, combined with (\ref{eq:zhan0}), proves (ii). We prove now (iii). We have by the Markov property at time $k$,
$$
\E\left[ \sum_{u\in \mathcal S_k^{z-A}\cap {(\mathcal T^{z-A})^c} }  B_n^z(u) \right]
=
\E\left[ \sum_{u\in \mathcal S_k^{z-A}\cap {(\mathcal T^{z-A})^c}} \Phi^{\rm kill}_{k,n}(0,z+V(u))\right]
$$

\noindent where $\Phi^{\rm kill}_{k,n}$ is defined in (\ref{def:Phikn}). By (\ref{eq:boundphikn}), this implies that
\begin{equation}\label{eq:moment1a}
\E\left[ \sum_{u\in \mathcal S^{z-A}\cap {(\mathcal T^{z-A})^c}} B_n^z(u){\bf 1}_{\{|u|\le n/2\}} \right]
\le
 c_{26} \ee^{-z}\E\left[ \sum_{u\in \mathcal S^{z-A}\cap {(\mathcal T^{z-A})^c}} \ee^{-V(u)}{\bf 1}_{\{|u|\le n/2\}} \right].
\end{equation}

\noindent At this stage, we make use of the measure $\Q$, introduced in Section \ref{s:spine}. We recall that under $\Q$, we identified our branching random with $\hat {\mathcal B}$. By definition of $\Q$ then Proposition \ref{p:spine} (i), we have for any $k\le n/2$,
\begin{eqnarray}
\E\left[ \sum_{u\in \mathcal S^{z-A}_k\cap {(\mathcal T^{z-A})^c}}  \ee^{-V(u)}   \right]
&=&
\hat{\E}\left[{1\over W_k}  \sum_{u\in \mathcal S^{z-A}_k\cap {(\mathcal T^{z-A})^c}}  \ee^{-V(u)}   \right]  \nonumber \\
&=&
\Q(w_k\in \mathcal S_k^{z-A}\cap({\mathcal T}^{z-A})^c). \label{eq:moment1b}
\end{eqnarray}

\noindent The right-hand side is equal to $0$ when $k=0$ since $w_0 \in \mathcal T^{z-A}$ by definition. For $k\ge 1$, we observe that ${\bf 1}_{\{w_k\in (\mathcal T^{z-A})^c\}} \le \sum_{\ell=1}^{k}{\bf 1}_{\{\xi(w_\ell)\ge \ee^{(V(w_{\ell-1})+z-A)/2}\}}$. It follows that
$$
 \Q(w_k \in \mathcal S^{z-A}_k\cap (\mathcal T^{z-A})^c)
 \le
 \sum_{\ell=1}^k {\Q}\left(w_k\in \mathcal S^{z-A},\, \xi(w_\ell) \ge \ee^{(V(w_{\ell-1})+z-A)/2} \right).
$$

\noindent Together with equations (\ref{eq:moment1a}) and (\ref{eq:moment1b}), this gives that
$$
\E\left[ \sum_{u\in \mathcal S^{z-A}\cap {(\mathcal T^{z-A})^c}} B_n^z(u){\bf 1}_{\{|u|\le n/2\}} \right]
\le
c_{26}\ee^{-z} \sum_{\ell=1}^{\lfloor n/2 \rfloor}\sum_{k=\ell}^{\lfloor n/2\rfloor} \Q\left(w_k\in \mathcal S^{z-A},\, \xi(w_\ell) \ge \ee^{(V(w_{\ell-1})+z-A)/2} \right).
$$

\noindent In order to prove (iii), it is enough to show that
\begin{equation}\label{eq:proof(ii)}
\sum_{\ell\ge 1}\sum_{k\ge \ell} \Q\left(w_k\in \mathcal S^{z-A},\, \xi(w_\ell) \ge \ee^{(V(w_{\ell-1})+z-A)/2} \right) = o(z)
\end{equation}

\noindent uniformly in $A\ge 0$ as $z\to\infty$. The left-hand side of (\ref{eq:proof(ii)}) is $0$ if $z<A$. Therefore, we will assume that $z\ge A$. For $k\ge \ell$, notice that if  $w_k \in \mathcal S^{z-A}$, then necessarily $\min_{j\le \ell} V(w_j)\ge A-z$, $V(w_k)\ge A-z$ and $V(w_k)< \min_{\ell \le j \le k-1} V(w_j)$ (in particular, $k$ is a ladder epoch for the random walk started at $V(w_\ell)$). This implies that
\begin{eqnarray*}
&& \Q\left(w_k\in \mathcal S^{z-A},\, \xi(w_\ell) \ge \ee^{(V(w_{\ell-1})+z-A)/2} \right) \\
&\le&
\Q\left( \xi(w_\ell) \ge \ee^{(V(w_{\ell-1})+z-A)/2},\, \min_{j\le \ell} V(w_j)\ge A-z,\,A-z\le V(w_k)< \min_{\ell \le j \le k-1} V(w_j)\right).
\end{eqnarray*}

\noindent Summing over $k\ge \ell$, we get
\begin{eqnarray*}
&& \sum_{k\ge \ell} \Q\left(w_k\in \mathcal S^{z-A},\, \xi(w_\ell) \ge \ee^{(V(w_{\ell-1})+z-A)/2} \right)\\
&\le& \hat{\E}\left[ {\bf 1}_{\{ \xi(w_\ell) \ge \ee^{(V(w_{\ell-1})+z-A)/2} \}} {\bf 1}_{\{\min_{j\le \ell}V(w_j)\ge A-z \}}\sum_{k\ge \ell}{\bf 1}_{\{A-z\le V(w_k)<\min_{\ell \le j \le k-1}V(w_j)\}}  \right].
\end{eqnarray*}

\noindent By the Markov property at time $\ell$, we recognize in the term $\sum_{k\ge \ell}{\bf 1}_{\{A-z\le V(w_k)<\min_{\ell \le j \le k-1}V(w_j)\}}$ the number of strict descending ladder heights above level $A-z$ when starting from $V(w_\ell)$. Consequently,
\begin{eqnarray*}
&& \sum_{k\ge \ell} \Q\left(w_k\in \mathcal S^{z-A},\, \xi(w_\ell) \ge \ee^{(V(w_{\ell-1})+z-A)/2} \right)\\
 &\le& \hat{\E}\left[ {\bf 1}_{\{ \xi(w_\ell) \ge \ee^{(V(w_{\ell-1})+z-A)/2} \}}{\bf 1}_{\{\min_{j\le \ell} V(w_j) \ge A-z \}} R(z-A+V(w_\ell))\right].
\end{eqnarray*}

\noindent We know from (\ref{c0}) that there exists $c_{27}>0$ such that $R(x)\le c_{27}(1+ x)_+$ for any real $x$. Thus, $R(z-A+V(w_\ell)) \le c_{27} (1+z-A + V(w_{\ell -1}))_+ + c_{27}(V(w_{\ell}) - V(w_{\ell -1}))_+$. Also, we obviously have $\min_{j\le \ell} V(w_j) \le \min_{j\le \ell-1} V(w_j)$. This yields that
$$
 \sum_{k\ge \ell} \Q\left(w_k\in \mathcal S^{z-A},\, \xi(w_\ell) \ge \ee^{(V(w_{\ell-1})+z-A)/2} \right)
 \le c_{27}(f(\ell) + g(\ell))
$$

\noindent where
\begin{eqnarray*}
f(\ell) &:=& \hat{\E}\left[ {\bf 1}_{\{ \xi(w_\ell) \ge \ee^{(V(w_{\ell-1})+z-A)/2} \}}{\bf 1}_{\{\min_{j\le \ell-1} V(w_j) \ge A-z \}}(1+z-A+V(w_{\ell-1})) \right],\\
 g(\ell) &:=&  \hat{\E}\left[ {\bf 1}_{\{ \xi(w_\ell) \ge \ee^{(V(w_{\ell-1})+z-A)/2} \}}{\bf 1}_{\{\min_{j\le \ell-1} V(w_j) \ge A-z \}} (V(w_\ell) - V(w_{\ell}-1))_+\right].
\end{eqnarray*}

\noindent Equation (\ref{eq:proof(ii)}) boils down to
\begin{equation}\label{eq:proof(ii)bis}
\sum_{\ell \ge 1} (f(\ell) + g(\ell)) = o(z).
\end{equation}

\noindent Let $(\xi,\Delta)$ be a generic random variable independent of all the random variables used so far, and distributed as $(\xi(w_1),V(w_1))$ (under $\Q$). Using the Markov property at time $\ell-1$ in $f(\ell)$, we get
$$
f(\ell) = \hat{\E}\left[ {\bf 1}_{\{ \xi \ge \ee^{(V(w_{\ell-1})+z-A)/2} \}}{\bf 1}_{\{\min_{j\le \ell-1} V(w_j) \ge A-z \}}(1+z-A+V(w_{\ell-1})) \right].
$$

\noindent Summing over $\ell$ (and replacing $\ell-1$ by $\ell$) yields that
$$
\sum_{\ell\ge 1} f(\ell) =  \hat{\E}\left[ \sum_{\ell \ge 0}{\bf 1}_{\{ V(w_{\ell})+z-A \le 2\ln(\xi) \}}{\bf 1}_{\{\min_{j\le \ell} V(w_j) \ge A-z \}}(1+z-A+V(w_{\ell})) \right].
$$

\noindent By Lemma \ref{l:estRW} (i), there exists $c_{28}>0$ such that for any $x\ge 0$
\begin{eqnarray*}
&& \hat{\E}\left[ \sum_{\ell \ge 0}{\bf 1}_{\{ V(w_{\ell})+z-A \le x \}}{\bf 1}_{\{\min_{j\le \ell} V(w_j) \ge A-z \}}(1+z-A+V(w_{\ell})) \right] \\
&\le& (1+x)\sum_{\ell \ge 0} \Q\left( V(w_{\ell})+z-A \le x,\,\min_{j\le \ell} V(w_j) \ge A-z\right)\\
&\le& c_{28} (1+x)^2(1+\min(x,z-A))\\
&\le& c_{28} (1+x)^2(1+\min(x,z)).
\end{eqnarray*}

\noindent We deduce that, with the notation of (\ref{def:X}),
\begin{eqnarray}
\sum_{\ell\ge 1} f(\ell) &\le& c_{28} \hat{\E}[(1+2\ln_+ \xi)^2(1+\min(2 \ln_+ \xi,z))] \nonumber \\
 &\le& c_{28}\E[X(1+2\ln_+ (X+\tilde X))^2(1+\min(2 \ln_+ (X + \tilde X),z))] \nonumber \\
 &=& o(z)\label{eq:f(ell)}
\end{eqnarray}

\noindent under (\ref{cond-ln}) by Lemma \ref{l:estcond} (ii). We now consider $g(\ell)$. We have similarly
$$
\sum_{\ell\ge 1} g(\ell) =  \hat{\E}\left[ \Delta_+\sum_{\ell \ge 0}{\bf 1}_{\{ V(w_{\ell})+z-A \le 2\ln(\xi) \}}{\bf 1}_{\{\min_{j\le \ell} V(w_j) \ge A-z \}} \right].
$$

\noindent From Lemma \ref{l:estRW} (i), we get
\begin{eqnarray}
\sum_{\ell\ge 1} g(\ell) &\le& c_{28} \hat{\E}[\Delta_+(1+2\ln_+ \xi)(1+\min(2 \ln_+ \xi,z))] \nonumber \\
 &\le& c_{28}\E[\tilde X (1+ 2 \ln_+(X+\tilde X))(1+\min(2 \ln_+ (X+\tilde X),z))] \nonumber \\
 &=& o(z)\label{eq:g(ell)}
\end{eqnarray}

\noindent by Lemma \ref{l:estcond} (ii). Equations (\ref{eq:f(ell)}) and (\ref{eq:g(ell)}) imply (\ref{eq:proof(ii)bis}), and so complete the proof of (ii). \hfill $\Box$

\bigskip

\noindent {\bf Remark}. Equations (\ref{eq:momentsnl}), (\ref{eq:boundphikn}) and (\ref{eq:manyto1S(z-A)}) imply that for any $n\ge 1$ and $z\ge 0$,
\begin{equation}
\E\left[  \sum_{u\in \mathcal S^{z}} B_n^z(u){\bf 1}_{\{ |u|\le  n^{1/2}\}} \right]\le c_{25}2^{3/2} R(z)\ee^{-z}.
\end{equation}

\bigskip

We compute the second moment in the following lemma. 

\begin{lemma}\label{l:moment2}
 There exists a constant $c_{29}>0$ such that for any $z\ge A\ge 0$, and any integer $n\ge 1$,
\begin{equation}
\label{varbound2}
\E\left[U^2\right]-\E\left[ U \right]  \le c_{29} \ee^{-z}\ee^{-A}
\end{equation}

\noindent where $U:= \sum_{u\in \mathcal S^{z-A}\cap \mathcal{T}^{z-A}} B_n^z(u){\bf 1}_{\{|u|\le n/2\}}$.
\end{lemma}

\noindent {\it Proof}. Let $U$ be as in the lemma. We observe that
$$
U^2 - U = \sum_{u\neq v}B_n^z(u) B_n^z(v) {\bf 1}_{\{ u,v \in \mathcal S^{z-A}\}}{\bf 1}_{\{u,v \in \mathcal T^{z-A} \} }{\bf 1}_{\{|u|,|v| \le n/2\}}.
$$

\noindent  It follows that
\begin{eqnarray*}
\E[U^2-U]
&=&
\E\left[\sum_{u\neq v} B_n^z(u) B_n^z(v){\bf 1}_{\{ u,v \in \mathcal S^{z-A}\}}{\bf 1}_{\{u,v \in \mathcal T^{z-A} \} }{\bf 1}_{\{|u|,|v| \le n/2\}}\right]\\
&\le&
2 \E\left[\sum_{u\neq v,|u|\ge |v|} B_n^z(u) B_n^z(v){\bf 1}_{\{ u,v \in \mathcal S^{z-A} \}}{\bf 1}_{ \{ u\in \mathcal T^{z-A}\}}{\bf 1}_{\{|u|\le n/2\}}\right].
\end{eqnarray*}

\noindent For $|u|\ge |v|$, and $u\neq v$, notice that $B_n^z(u)$ depends on the branching random walk rooted at $u$, whereas $B_n^z(v){\bf 1}_{\{ u,v \in \mathcal S^{z-A} \}}{\bf 1}_{\{ u \in \mathcal T^{z-A}\}}$ is independent of it (even if $v$ is a (strict) ancestor of $u$). Therefore, by the branching property,
$$\E[U^2-U]
\le
2 \E\left[\sum_{u\neq v,|u|\ge |v|}\Phi^{\rm kill}_{|u|,n}(0,z+V(u)) B_n^z(v){\bf 1}_{\{ u,v \in \mathcal S^{z-A} \}}{\bf 1}_{\{ u \in \mathcal T^{z-A}\}}{\bf 1}_{\{|u|\le n/2\}}\right]
$$

\noindent where $\Phi^{\rm kill}_{k,n}$ is defined in (\ref{def:Phikn}). By (\ref{eq:boundphikn}), we have $\Phi^{\rm kill}_{|u|,n}(z+V(u)) \le c_{26} \ee^{-z-V(u)}$ for $|u|\le n/2$. This gives that
\begin{eqnarray}\nonumber
\E[U^2-U]
&\le&
c_{26}\ee^{-z} \E\left[\sum_{u\neq v,|u|\ge |v|}\ee^{-V(u)} B_n^z(v){\bf 1}_{\{ u,v \in \mathcal S^{z-A} \}}{\bf 1}_{\{u \in \mathcal T^{z-A} \}}{\bf 1}_{\{|u|\le n/2\}}\right]\\
&\le&
c_{26}\ee^{-z} \sum_{k=1}^{\lfloor n/2\rfloor}\E\left[\sum_{u \in \mathcal{S}^{z-A}_k\cap \mathcal T^{z-A}}\ee^{-V(u)} \sum_{v\neq u,|v|\le k} B_n^z(v){\bf 1}_{\{ v \in \mathcal S^{z-A} \}}\right].\label{eq:moment2}
\end{eqnarray}

\noindent The weight $\ee^{-V(u)}$ hints at a change of measure from $\P$ to $\Q$. For any $k\in [0,n/2]$, we have by Proposition \ref{p:spine} (i),
\begin{eqnarray}\nonumber
&& \E\left[\sum_{u \in \mathcal{S}^{z-A}_k\cap \mathcal T^{z-A}}\ee^{-V(u)}  \sum_{v\neq u,|v|\le k} B_n^z(v){\bf 1}_{\{ v \in \mathcal S^{z-A} \}}\right]\\
\label{eq:moment2'}&=&
\hat{\E}\left[{\bf 1}_{\{ w_k \in \mathcal S^{z-A}\cap \mathcal T^{z-A} \}}  \sum_{v\neq w_k,|v|\le k}  B_n^z(v){\bf 1}_{\{ v \in \mathcal S^{z-A} \}} \right].
\end{eqnarray}

\noindent We have to split the cases depending on the location of the vertex $v$ with respect to $w_k$. We say that $u\nsim v$ if neither $v$ nor $u$ is an ancestor of the other. If $v\neq w_k$ and $|v|\le k$, then either $v\nsim u$, or $v=w_\ell$ for some $\ell<k$. In view of (\ref{eq:moment2}) and (\ref{eq:moment2'}), the lemma will be proved once the following two estimates are shown:
\begin{eqnarray}\label{eq:moment2a}
\sum_{k= 1}^{ \lfloor n/2 \rfloor } \hat {\E} \left[ \sum_{v\nsim w_k} B_n^z(v){\bf 1}_{\{v \in \mathcal S^{z-A} \}},w_k\in \mathcal S^{z-A} \cap \mathcal T^{z-A} \right] &\le& c_{31}\ee^{-A},\\
\sum_{k=1}^{ \lfloor n/2 \rfloor }\sum_{\ell=0}^{k-1} \hat {\E} \left[  B_n^z(w_{\ell}), \,w_\ell \in \mathcal S^{z-A},\,w_{k} \in \mathcal S^{z-A}\cap \mathcal T^{z-A} \right] &\le& c_{32}\ee^{-A}.\label{eq:moment2b}
\end{eqnarray}

\noindent {\it Proof of equation (\ref{eq:moment2a})}.\\ 
Decomposing the sum $\sum_{v \nsim w_k}$ along the spine, we see that for any $k\in [ 1,n/2]$,
\begin{equation}\label{eq:decomp}
\sum_{v\nsim w_k}  B_n^z(v){\bf 1}_{\{ v \in \mathcal S^{z-A} \}}
=
\sum_{\ell =1}^{k}\sum_{x\in \Omega(w_\ell)} \sum_{v\ge x} B_n^z(v){\bf 1}_{\{ v \in \mathcal S^{z-A} \}}.
\end{equation}

\noindent  The branching random walk rooted at $x\in \Omega(w_\ell)$ has the same law under $\P$ and $\Q$. Recall the definition of $\hat{\mathcal G}_\infty$ in (\ref{def:Ginfty}). We have for $\ell\le n/2$ and $x\in \Omega(w_\ell)$,
\begin{equation*}
\hat{\E}\left[\sum_{v\ge x} B_n^z(v){\bf 1}_{\{ v \in \mathcal S^{z-A} \}} \,\Big|\, \hat{\mathcal G}_\infty \right]
=
\hat{\E}\left[\sum_{v\ge x} \Phi^{\rm kill}_{|v|,n}(0,z+V(v)){\bf 1}_{\{ v \in \mathcal S^{z-A} \}} \,\Big|\, \hat{\mathcal G}_\infty \right]
\end{equation*}

\noindent with the notation of (\ref{def:Phikn}), and (\ref{eq:boundphikn}) implies that
\begin{equation}\label{eq:moment2c}
\hat{\E}\left[\sum_{v\ge x} B_n^z(v){\bf 1}_{\{ v \in \mathcal S^{z-A} \}} \,\Big|\, \hat{\mathcal G}_\infty \right]\le c_{26} \ee^{-z}\hat{\E}\left[\sum_{v\ge x} \ee^{-V(v)}{\bf 1}_{\{ v \in \mathcal S^{z-A} \}}  \,\Big|\, \hat{\mathcal G}_\infty\right].
\end{equation}

\noindent We observe now that if $v\ge x$ and $v \in \mathcal S^{z-A}$, then  $\min_{|x| \le j\le |v|-1} V(v_j)>V(v) \ge A-z$. Therefore, by the Markov property,
\begin{eqnarray*}
\hat {\E} \left[\sum_{v\ge x} \ee^{-V(v)}{\bf 1}_{\{ v \in \mathcal S^{z-A} \}}  \,\Big|\, \hat{\mathcal G}_\infty\right]
\le
\E_{V(x)}\left[\sum_{v\in \T} \ee^{-V(v)}{\bf 1}_{\{ \min_{ j\le |v|-1} V(v_j)>V(v) \ge A-z \}} \right].
\end{eqnarray*}

\noindent By (\ref{many-to-one}) and the definition of the renewal function $R(x)$ in  (\ref{def:R(x)}), we observe that
$$
\E_{V(x)}\left[\sum_{v\in \T} \ee^{-V(v)}{\bf 1}_{\{ \min_{ j\le |v|-1} V(v_j) >V(v) \ge A-z \}} \right]
=
\ee^{-V(x)} R(z-A+V(x)).
$$

\noindent Going back to (\ref{eq:moment2c}), we get that for any $\ell\le n/2$ and $x\in \Omega(w_\ell)$,
$$
\hat{\E}\left[\sum_{v\ge x}  B_n^z(v){\bf 1}_{\{ v \in \mathcal S^{z-A} \}} \,\Big|\, \hat{\mathcal G}_\infty \right]
\le
c_{26}\ee^{-z}\ee^{-V(x)} R(z-A+V(x)).
$$

\noindent In view of (\ref{eq:decomp}), we obtain that
\begin{eqnarray}\label{eq:moment2a0}
&&\sum_{k= 1}^{\lfloor n/2\rfloor} \hat{\E}\left[ \sum_{v\nsim w_k} B_n^z(v) {\bf 1}_{\{v \in \mathcal S^{z-A} \}},w_k\in \mathcal S^{z-A}\cap\mathcal T^{z-A} \right]\\
&\le& c_{26}\ee^{-z} \sum_{k\ge 1} \sum_{\ell=1}^k \hat{\E} \left[\sum_{x\in \Omega(w_\ell)} \ee^{-V(x)} R(z-A+V(x)),\,w_k\in \mathcal S^{z-A}\cap \mathcal T^{z-A}\right].\nonumber
\end{eqnarray}

\noindent We look at $R(z-A+V(x))$ for $x\in \Omega(w_\ell)$. If $V(x)\le V(w_{\ell-1})$ and $z-A+V(w_{\ell-1})\ge 0$, we have 
$$
R(z-A+V(x))\le R(z-A+V(w_{\ell-1}))\le c_{27}(1+z-A+V(w_{\ell-1})).$$

\noindent  If $V(x)> V(w_{\ell-1})$ and $z-A+V(w_{\ell-1})\ge 0$, we write that
\begin{eqnarray*}
R(z-A+V(x)) &\le& c_{27}(1+z-A+V(x))\\
&\le&
c_{27}(1+z-A+V(w_{\ell-1}))(1+V(x)-V(w_\ell-1)).
\end{eqnarray*}

\noindent  Therefore, for any $\ell\le k$, we have on the event that $w_k\in \mathcal S^{z-A}$,
\begin{eqnarray*}
&&\sum_{x\in \Omega(w_\ell)} \ee^{-V(x)} R(z-A+V(x)) \\
&\le& 
c_{27} (1+z-A+V(w_{\ell-1}) )  \sum_{x\in \Omega(w_\ell)} (1+(V(x)-V(w_{\ell-1}))_+) \ee^{-V(x)} \\
&=&
c_{27} (1+z-A+V(w_{\ell-1}) )\ee^{-V(w_{\ell-1})} \xi(w_\ell)
\end{eqnarray*}

\noindent by definition (\ref{def:xi}). On the event that $w_k\in \mathcal S^{z-A}\cap \mathcal T^{z-A}$, we conclude that
$$
\sum_{x\in \Omega(w_\ell)} \ee^{-V(x)} R(z-A+V(x))  \le c_{27}  \ee^{(z-A)/2}(1+z-A+V(w_{\ell-1}) )\ee^{-V(w_{\ell-1})/2}.
$$

\noindent Therefore, we have by (\ref{eq:moment2a0}),
\begin{eqnarray}\label{eq:moment2a1}
&&\sum_{k= 1}^{\lfloor n/2\rfloor} \hat{\E}\left[ \sum_{v\nsim w_k} B_n^z(v) {\bf 1}_{\{v \in \mathcal S^{z-A} \}},w_k\in \mathcal S^{z-A}\cap\mathcal T^{z-A} \right]\\
&\le& c_{26}c_{27} \ee^{-(z+A)/2}\sum_{k\ge 1} \sum_{\ell=1}^k \hat{\E} \left[(1+z-A+V(w_{\ell-1}) )\ee^{-V(w_{\ell-1})/2},\,w_k\in \mathcal S^{z-A}\right].\nonumber
\end{eqnarray}

\noindent  Proposition \ref{p:spine} (ii) says that
\begin{eqnarray*}
&& \hat{\E} \left[(z-A+V(w_{\ell-1}) +1)\ee^{-V(w_{\ell-1})/2},\,w_k\in \mathcal S^{z-A}\right]\\
&=&
\E\left[\ee^{-S_{\ell-1}/2} (1+z-A+S_{\ell-1} ),\, \min_{j\le k-1} S_j>S_k \ge A-z\right].
\end{eqnarray*}

\noindent We observe that
\begin{eqnarray*}
&& \sum_{k\ge 1} \sum_{\ell=1}^k \E\left[\ee^{-S_{\ell-1}/2} (1+z-A+S_{\ell-1} ),\, \min_{j\le k-1} S_j> S_k \ge A-z\right]\\
&=& \sum_{\ell \ge 1} \E\left[\ee^{-S_{\ell-1}/2} (1+z-A+S_{\ell-1} ) \sum_{k\ge \ell }{\bf 1}_{\{ \min_{j\le k-1} S_j>S_k \ge A-z\}}\right].
\end{eqnarray*}

\noindent Since
$$
\sum_{k\ge \ell }{\bf 1}_{\{ \min_{j\le k-1} S_j>S_k \ge A-z\}}
\le 
{\bf 1}_{\{ \min_{j\le \ell-1} S_j \ge A-z\}}\sum_{k\ge \ell }{\bf 1}_{\{ \min_{j\in[\ell-1, k-1]} S_j>S_k \ge A-z\}},
$$

\noindent we deduce by the Markov property at time $\ell-1$ that
\begin{eqnarray*}
&& \sum_{k\ge 1} \sum_{\ell=1}^k \E\left[\ee^{-S_{\ell-1}/2} (1+z-A+S_{\ell-1} ),\,\min_{j\le k-1} S_j>S_k \ge A-z\right]\\
&\le& \sum_{\ell \ge 1} \E\left[\ee^{-S_{\ell-1}/2} (1+z-A+S_{\ell-1} ) R(S_{\ell-1} +z-A),\, \min_{j\le \ell-1}S_j\ge A-z\right]\\
&\le&
c_{27}\sum_{\ell \ge 1} \E\left[\ee^{-S_{\ell-1}/2} (1+z-A+S_{\ell-1} )^2,\, \min_{j\le \ell-1}S_j\ge A-z\right].
\end{eqnarray*}

\noindent Using this bound in (\ref{eq:moment2a1}) yields that
\begin{eqnarray}\label{eq:moment2a2}
&&
\sum_{k= 1}^{\lfloor n/2\rfloor} \hat{\E}\left[ \sum_{v\nsim w_k} B_n^z(v) {\bf 1}_{\{v \in \mathcal S^{z-A} \}},w_k\in \mathcal S^{z-A}\cap\mathcal T^{z-A} \right]\\
\nonumber &\le&
c_{26}c_{27}^2\ee^{-(z+A)/2}\sum_{\ell \ge 1}  \E\left[\ee^{-S_{\ell-1}/2} (1+z-A+S_{\ell-1} )^2,\, \min_{j\le \ell-1} S_j\ge A-z\right].
\end{eqnarray}

\noindent By Lemma \ref{l:estRW} (iii), we have
$$
\sum_{\ell \ge 1} \E\left[\ee^{-S_{\ell-1}/2} (1+z-A+S_{\ell-1} )^2,\, \min_{j\le \ell-1}S_j\ge A-z\right] \le c_{33}\ee^{(z-A)/2}.
$$

\noindent Consequently, by (\ref{eq:moment2a2})
$$
\sum_{k= 1}^{\lfloor n/2\rfloor} \hat{\E}\left[ \sum_{v\nsim w_k} B_n^z(v) {\bf 1}_{\{v \in \mathcal S^{z-A} \}},w_k\in \mathcal S^{z-A}\cap\mathcal T^{z-A} \right]\
\le
c_{34} \ee^{-(z+A)/2}\ee^{(z-A)/2} = c_{34}\ee^{-A}.
$$

\noindent Equation (\ref{eq:moment2a}) follows. \\

\noindent{\it Proof of equation (\ref{eq:moment2b})}\\  
We have
\begin{eqnarray}\nonumber
 &&\sum_{k= 1}^{\lfloor n/2 \rfloor }\sum_{\ell=0}^{k-1} \hat{\E}\left[  B_n^z(w_{\ell}),\, w_\ell \in \mathcal S^{z-A},\,w_{k} \in \mathcal S^{z-A}\cap \mathcal T^{z-A} \right]\\
\nonumber &=&
\sum_{\ell= 0}^{\lfloor n/2\rfloor-1} \sum_{k=\ell+1}^{\lfloor n/2 \rfloor}  \hat{\E}\left[  B_n^z(w_{\ell}), \,w_\ell\in \mathcal S^{z-A},\,w_{k} \in \mathcal S^{z-A} \cap \mathcal T^{z-A} \right]\\
\label{eq:moment2b0}&= &
\sum_{\ell= 0}^{\lfloor n/2 \rfloor -1}  \hat{\E}\left[  B_n^z(w_{\ell}){\bf 1}_{\{w_{\ell} \in \mathcal S^{z-A}\}}\sum_{k=\ell+1}^{\lfloor n/2\rfloor} {\bf 1}_{\{ w_k \in \mathcal S^{z-A} \cap \mathcal T^{z-A}\}} \right].
\end{eqnarray}

\noindent Let $t_\ell$ be the first time $t$ after $\ell$ such that $V(w_t)< V(w_\ell)$. If $k>\ell$ and $w_k \in \mathcal S^{z-A}$, then $V(w_k) < V(w_\ell)$, which means that necessarily $k\ge t_\ell$. Moreover, if $k\ge t_\ell$  and $w_k\in \mathcal T^{z-A}$, then $w_{t_\ell} \in \mathcal T^{z-A}$. Thus, for any $\ell\ge 0$,
\begin{eqnarray*}
\sum_{k=\ell +1}^{\lfloor n/2\rfloor} {\bf 1}_{\{ w_k \in \mathcal S^{z-A}\cap \mathcal T^{z-A}\}} &=&\sum_{k= t_\ell}^{\lfloor n/2\rfloor} {\bf 1}_{\{ w_k \in \mathcal S^{z-A}\cap \mathcal T^{z-A}\}}\\
&\le&
{\bf 1}_{\{ w_{t_\ell} \in \mathcal T^{z-A}, t_\ell \le n/2 \}}\sum_{k\ge t_\ell} {\bf 1}_{\{  \min_{t_\ell \le j<k} V(w_j) > V(w_k) \ge A-z\}}.
\end{eqnarray*}

\noindent We observe that $B_n^z(w_{\ell})$ is a function of the branching random walk killed below $V(w_\ell)$ and therefore is independent of the subtree rooted at $w_{t_\ell}$. As a result, applying the branching property, we get that for any $\ell\in [0,n/2]$,
\begin{eqnarray*}
&& \hat{\E}\left[ B_n^z(w_{\ell}){\bf 1}_{\{w_{\ell} \in \mathcal S^{z-A}\}}\sum_{k=\ell+1}^{\lfloor n/2 \rfloor} {\bf 1}_{\{ w_k \in \mathcal S^{z-A}\cap \mathcal T^{z-A} \}} \right]\\
&\le&
\hat{\E}\left[ B_n^z(w_{\ell}){\bf 1}_{\{w_{\ell} \in \mathcal S^{z-A}\}}{\bf 1}_{\{w_{t_\ell} \in \mathcal T^{z-A},t_\ell \le n/2 \}}\sum_{k\ge t_\ell} {\bf 1}_{\{  \min_{t_\ell \le j<k} V(w_j) > V(w_k) \ge A-z\}}\right]\\
&=&
\hat{\E}\left[ B_n^z(w_{\ell}){\bf 1}_{\{w_{\ell} \in \mathcal S^{z-A}\}}{\bf 1}_{\{ w_{t_\ell} \in \mathcal T^{z-A},t_\ell\le n/2\}}R(z-A+V(w_{t_\ell}))\right].
\end{eqnarray*}

\noindent By equation (\ref{eq:moment2b0}), we deduce that
\begin{eqnarray*}
&&\sum_{k= 1}^{\lfloor n/2 \rfloor }\sum_{\ell=0}^{k-1} \hat{\E}\left[  B_n^z(w_{\ell}),\, w_\ell \in \mathcal S^{z-A},\,w_{k} \in \mathcal T^{z-A} \right]\\
&\le& 
\sum_{\ell= 0}^{\lfloor n/2 \rfloor -1}\hat{\E}\left[ B_n^z(w_{\ell}){\bf 1}_{\{w_{\ell} \in \mathcal S^{z-A}\}}{\bf 1}_{\{ w_{t_\ell} \in \mathcal T^{z-A}, t_\ell\le n/2\}}R(z-A+V(w_{t_\ell}))\right].
\end{eqnarray*}

\noindent We have $V(w_{t_\ell})<V(w_\ell)$. Since $R$ is a non-decreasing function, we obtain that
\begin{eqnarray}
\nonumber && \sum_{k= 1}^{\lfloor n/2 \rfloor }\sum_{\ell=0}^{k-1} \hat{\E}\left[  B_n^z(w_{\ell}),\, w_\ell \in \mathcal S^{z-A},\,w_{k} \in \mathcal S^{z-A}\cap \mathcal T^{z-A} \right]\\
\label{eq:moment2b00}&\le&
\sum_{\ell= 0}^{\lfloor n/2 \rfloor -1}\hat{\E}\left[ B_n^z(w_{\ell}){\bf 1}_{\{w_{\ell} \in \mathcal S^{z-A}\}}{\bf 1}_{\{ w_{t_\ell} \in \mathcal T^{z-A},t_\ell\le n/2\}}R(z-A+V(w_{\ell}))\right].
\end{eqnarray}

\noindent Recall from (\ref{def:Gell})  that $\hat {\mathcal G}_\ell$ is the $\sigma$-algebra generated by the spine and its siblings up to time $\ell$. The Markov property at time $\ell$ shows that
\begin{eqnarray*}
\hat {\E}\left[B_n^z(w_\ell), w_{t_\ell} \in \mathcal{T}^{z-A},\,t_\ell\le n/2 \, | \, \hat{\mathcal G}_\ell\right]
&=&
{\bf 1}_{\{w_\ell \in \mathcal T^{z-A}\}} \tilde \Phi^{\rm kill}_{\ell,n,A}(z+V(w_\ell))\\
&\le& 
\tilde \Phi^{\rm kill}_{\ell,n,A}(z+V(w_\ell))
\end{eqnarray*}

\noindent where, if $\tau_0^-:=\min\{j\ge 0\,:\, V(w_j)<0\}$, then for any integers $n\ge 1$, $\ell\le n/2$, any $A,r\ge 0$, we defined
$$
\tilde \Phi^{\rm kill}_{\ell,n,A}(r):= \Q\left(\tau_0^- \le (n/2)-\ell,\,M_{n-\ell}^{\rm kill} < a_n(r), \,\xi(w_j) \le \ee^{(r+V(w_{j-1})-A)/2},\, \forall\, 1\le j\le \tau_0^-\right).
$$

\noindent We deduce that, for any $n\ge 1$, any $\ell< n/2$, any $z\ge A\ge 0$,
\begin{eqnarray*}
&& \hat{\E}\left[ B_n^z(w_{\ell}){\bf 1}_{\{w_{\ell} \in \mathcal S^{z-A}\}}{\bf 1}_{\{ w_{t_\ell} \in \mathcal T^{z-A},t_\ell\le n/2\}}R(z-A+V(w_{\ell}))\right] \\
&\le&
\hat{\E}\left[{\bf 1}_{\{w_{\ell} \in \mathcal S^{z-A}\}}R(z-A+V(w_\ell)) \tilde \Phi^{\rm kill}_{\ell,n,A}(z+V(w_\ell))\right]\\
&=&
 \hat{\E}\left[  {\bf 1}_{\{ \min_{j<\ell} V(w_j) > V(w_\ell)\ge A-z \}} R(z-A+V(w_\ell)) \tilde \Phi^{\rm kill}_{\ell,n,A}(z+V(w_\ell))\right].
\end{eqnarray*}

\noindent By Proposition \ref{p:spine} (ii), this implies that
\begin{eqnarray}
&& \hat{\E}\left[ B_n^z(w_{\ell}){\bf 1}_{\{w_{\ell} \in \mathcal S^{z-A}\}}{\bf 1}_{\{ w_{t_\ell} \in \mathcal T^{z-A},t_\ell \le n/2\}}R(z-A+V(w_{\ell}))\right]\label{eq:moment2b1}\\
&\le&
\E\left[  {\bf 1}_{\{ \min_{j<\ell} S_j>S_\ell\ge A-z \}} R(z-A+S_\ell) \tilde \Phi^{\rm kill}_{\ell,n,A}(z+S_\ell)\right].\nonumber
\end{eqnarray}

\noindent Let us estimate $\tilde \Phi^{\rm kill}_{\ell,n,A}(r)$ for $\ell<n/2$. We have to decompose along the spine. Notice that if $M_{n-\ell}^{\rm kill} < a_n(r)$, and $\tau_0^-\le \lfloor n/2\rfloor-\ell$ then there must be some $j<\tau_0^-\le \lfloor n/2\rfloor-\ell$ and $x\in \Omega(w_j)$ such that there exists a line of descent from $x$ which stays above $0$ and ends below $a_n(r)$ at time $n-\ell$. Therefore, for any $n\ge 1$, $\ell< n/2$, and $A,r\ge 0$,
$$
 \tilde \Phi^{\rm kill}_{\ell,n,A}(r)
\le
\sum_{j=1}^{\lfloor n/2 \rfloor -\ell} \hat{\E}\left[\sum_{x\in \Omega(w_j)} \Phi^{\rm kill}_{\ell+j,n}(V(x),r), \, \xi(w_j) \le \ee^{(r+V(w_{j-1})-A)/2},\,j \le \tau_0^- \right]
$$

\noindent with the notation of (\ref{def:Phikn}). By (\ref{eq:boundphikn}), we get that
\begin{eqnarray}
\label{eq:tildephia}&& \tilde \Phi^{\rm kill}_{\ell,n,A}(r)\\
&\le&
c_{35}\ee^{-r}\sum_{j=1}^{\lfloor n/2 \rfloor -\ell} \hat{\E}\left[\sum_{x\in \Omega(w_j)} (1+V(x)_+)\ee^{-V(x)}, \, \xi(w_j) \le \ee^{(r+V(w_{j-1})-A)/2},\,j\le \tau_0^- \right].
\nonumber
\end{eqnarray}

\noindent We observe that
\begin{eqnarray*}
&&\sum_{x\in \Omega(w_j)} (1+V(x)_+)\ee^{-V(x)} \\
&\le& 
(1+V(w_{j-1})_+)\ee^{-V(w_{j-1})} \sum_{x\in \Omega(w_j)} (1+(V(x)-V(w_{j-1}))_+)\ee^{-(V(x)-V(w_{j-1}))} \\
&=&
(1+V(w_{j-1})_+)\ee^{-V(w_{j-1})} \xi(w_{j})
\end{eqnarray*}

\noindent by definition (\ref{def:xi}). We deduce from (\ref{eq:tildephia}) that
\begin{eqnarray*}
&& \tilde \Phi^{\rm kill}_{\ell,n,A}(r)\\
&\le&
c_{35}\ee^{-r}\sum_{j=1}^{\lfloor n/2 \rfloor -\ell} \hat{\E}\left[\ee^{-V(w_{j-1})}(1+V(w_{j-1}))\ee^{(r+V(w_{j-1}) -A)/2},\,j\le \tau_0^- \right].
\end{eqnarray*}

\noindent It follows that, for any $n\ge 1$, $\ell< n/2$, and  $A,r\ge 0$,
\begin{eqnarray*}
\tilde \Phi^{\rm kill}_{\ell,n,A}(r)
&\le&
c_{35}\ee^{-A}\ee^{-(r-A)/2}\sum_{j\ge 1} \E\left[\ee^{-S_{j-1}/2}(1+S_{j-1}),\,j\le \tau_0^- \right]\\
&=&
c_{36}\ee^{-A}\ee^{-(r-A)/2},
\end{eqnarray*}

\noindent by Lemma \ref{l:estRW} (ii). Going back to (\ref{eq:moment2b1}), we obtain that for any $n\ge 1$, $\ell< n/2$, $z\ge A \ge 0$,
\begin{eqnarray*}
&& \hat{\E}\left[ B_n^z(w_{\ell}){\bf 1}_{\{w_{\ell} \in \mathcal S^{z-A}\}}{\bf 1}_{\{ w_{t_\ell} \in \mathcal T^{z-A},t_\ell \le n/2\}}R(z-A+V(w_{\ell}))\right]\\
&\le&
c_{36}\ee^{-A}\E\left[  {\bf 1}_{\{ \min_{j<\ell} S_j > S_\ell\ge A-z \}} R(z-A+S_\ell) \ee^{-(S_\ell+z-A)/2}\right].
\end{eqnarray*}

\noindent Equation (\ref{eq:moment2b00}) yields that for any $n\ge 1$, and $z\ge A\ge 0$,
\begin{eqnarray*}
 && \sum_{k= 1}^{\lfloor n/2 \rfloor }\sum_{\ell=0}^{k-1} \hat{\E}\left[  B_n^z(w_{\ell}),\, w_\ell \in \mathcal S^{z-A},\,w_{k} \in \mathcal S^{z-A}\cap \mathcal T^{z-A} \right]\\
&\le&
c_{36}\ee^{-A}\E\left[  \sum_{\ell\ge 0} {\bf 1}_{\{ \min_{j<\ell} S_j > S_\ell\ge A-z \}} R(z-A+S_\ell) \ee^{-(S_\ell+z-A)/2}\right].
\end{eqnarray*}

\noindent Applying Lemma \ref{l:estRW} (iii) implies  (\ref{eq:moment2b}) and thus completes the proof of the lemma. \hfill $\Box$

\subsection{Proof of Proposition \ref{p:tailmin}}

\label{s:tailminmain}

We can now prove Proposition \ref{p:tailmin}.

\bigskip

\noindent {\it Proof of Proposition \ref{p:tailmin}}. Let $\varepsilon>0$.  For any $r\ge 0$, we observe that by (\ref{many-to-one})
\begin{eqnarray*}
\P(\exists \, u \in \T\,:\, V(u)\le -r)
&\le& \sum_{n\ge 0} \E\left[ \sum_{|u|=n}{\bf 1}_{\{ V(u)\le -r,V(u_k)>-r,\forall\, k<n \}} \right]\\
\nonumber &=&
\sum_{n\ge 0} \E\left[\ee^{S_n},S_n\le -r,S_k >-r \, \forall k<n \right]\\
\nonumber &\le& \ee^{-r}.
\end{eqnarray*}

\noindent Therefore
\begin{equation*}\label{eq:Z[-r]}
\P(\exists \, u \in \T\,:\, V(u)\le A-z) \le \ee^{A-z} .
\end{equation*}

\noindent For any $z\ge A\ge 0$, we observe that on the event $ \{\forall u\in\T,\, V(u)\ge A-z\}$, we have $M_n < \frac32\ln n -z$ if and only if $\sum_{u\in\mathcal S^{z-A}} B_n^z(u) \ge 1$ (recall the definition of $\mathcal S^{r}$ and $B_n^z$  in   (\ref{def:SA}) and in (\ref{def:Bn})). Therefore, for $n\ge 1$ and $z\ge A$, we have 
$$
0\le \P\left(M_n \le {3\over 2} \ln n -z \right) - \P\left(\sum_{u\in\mathcal S^{z-A}} B_n^z(u) \ge 1\right)
\le
 \ee^{A-z}.
$$

\noindent  We notice that $\P\left(\sum_{u\in\mathcal S^{z-A}} B_n^z(u) \ge 1\right)\le \E\left[\sum_{u\in\mathcal S^{z-A}} B_n^z(u)\right]$. Hence,
$$
\P\left(M_n \le {3\over 2} \ln n -z \right)
\le
\ee^{A-z} + \E\left[\sum_{u\in\mathcal S^{z-A}} B_n^z(u)\right].
$$

\noindent Lemma \ref{l:moment1} (i) and (ii) implies that for $n\ge N_1$ and $z\in [A_1,(3/2)\ln(n)-A_1]$,
$$
{\ee^z \over R(z-A_1)} \P\left(M_n \le {3\over 2} \ln n -z \right)- C_1
\le
 {\ee^{A_1} +c \over R(z-A_1)}+\varepsilon.
$$

\noindent Since $R(x)\sim c_0 x$ at infinity by (\ref{c0}), we have for $n\ge N_1$ and $z\in [A_2,(3/2)\ln(n)-A_1]$,
$$
{\ee^z \over c_0 z} \P\left(M_n \le {3\over 2} \ln n -z \right)- C_1
\le
 {\ee^{A_1}+c \over c_0 z}+ 2\varepsilon.
$$

\noindent We deduce that for $n\ge N_1$ and $z\in[A_3,(3/2)\ln(n)-A_1]$,
$$
{\ee^z \over c_0 z} \P\left(M_n \le {3\over 2} \ln n -z \right)- C_1 \le 3\varepsilon.
$$

\noindent This proves the upper bound. Similarly, we have for the lower bound
\begin{eqnarray*}
\P\left(M_n \le {3\over 2} \ln n -z \right) &\ge&
\P\left(\sum_{u\in\mathcal S^{z-A}}  B_n^z(u) \ge 1\right)\\
&\ge&
\P\left(\sum_{u\in\mathcal S^{z-A}\cap\mathcal T^{z-A}} B_n^z(u){\bf 1}_{\{ |u|\le n/2 \}} \ge 1\right) .
\end{eqnarray*}

\noindent If we write $U(A):=\sum_{u\in\mathcal S^{z-A}\cap \mathcal T^{z-A}} B_n^z(u){\bf 1}_{\{ |u|\le n/2 \}}$, then by the Paley-Zygmund formula, we have $\P(U(A)\ge 1) \ge {\E[U(A)]^2 \over \E[U(A)^2]}$. By Lemma \ref{l:moment1}, we know that ${\ee^z\over R(z-A_4)}\E[U(A_4)] \ge C_1 - \varepsilon$ for $n\ge N_2$ and $z\in [A_5,(3/2)\ln(n)-A_4]$. By Lemma \ref{l:moment2}, we have that $\E[U(A_4)^2] \le (1+\varepsilon)\E[U(A_4)]$ if $A_5$ is taken large enough. Hence, ${\ee^z\over R(z-A_4)}\P(U(A_4)\ge 1) \ge {\ee^z\over R(z-A_4)}(1+\varepsilon)^{-1}\E[U(A_4)] \ge (1+\varepsilon)^{-1}( C_1-\varepsilon)$. This yields that
$$
{\ee^z\over R(z-A_4)}\P\left(M_n \le {3\over 2} \ln n -z \right) \ge  (1+\varepsilon)^{-1}( C_1-\varepsilon).
$$

\noindent From here, we proceed as before to see that for $n\ge N_2$ and $z\in[A_6,(3/2)\ln(n)-A_4]$,
$$
{\ee^z\over c_0z}\P\left(M_n \le {3\over 2} \ln n -z \right) \ge C_1 - c_{37} \varepsilon.
$$

\noindent The proposition follows. \hfill $\Box$

\section{Proof of Theorem \ref{main}}
\label{s:main}

For $\beta\ge 0$, we look at the branching random walk killed below $-\beta$. The population at  time $n$ of this process is $\{|u|=n\,:\, V(u_k) \ge -\beta,\,\forall \, k\le n\}$. We define the associated martingale (see Appendix \ref{s:derivative})
\begin{equation}\label{def:Dnbeta}
D_n^{(\beta)} := \sum_{|u|=n} R(\beta +V(u))\ee^{-V(u)}{\bf 1}_{\{V(u_k)\ge -\beta,\, k\le n\}}.
\end{equation}

\noindent Since $D_n^{(\beta)}$ is non-negative, it has a limit almost surely that we denote by $D_\infty^{(\beta)}$. Under (\ref{cond-2}) and (\ref{cond-ln}), we know by Proposition \ref{l:Dbeta} that $D_\infty^{(\beta)}>0$ almost surely on the event of non-extinction for the  branching random walk killed below $-\beta$. For $A\ge 0$, let $\mathcal Z[A]$ denote the set of particles absorbed at level $A$, i.e.
$$
\mathcal Z[A] := \{u\in\T\,:\, V(u)\ge A,\, V(u_k)<A\, \forall\, k<|u|\}.
$$

\noindent In the words of Section 6 in \cite{biggins-kyprianou04}, this set is a very simple optional line. By Theorem 6.1 (and Lemma 6.1) of \cite{biggins-kyprianou04}, we know that $\sum_{u\in\mathcal Z[A]} R(\beta +V(u))\ee^{-V(u)}{\bf 1}_{\{V(u_k)\ge -\beta\}}$ converges to $D_\infty^{(\beta)}$ almost surely as $A\to\infty$. Recall that $R(x)\sim c_0 x$ at infinity by (\ref{c0}). Recall from (\ref{def:Wn}) that the martingale $W_n$ is defined by
$$
W_n:=\sum_{|x|=n}\ee^{-V(x)}
$$

\noindent and we know from \cite{lyons} that $W_n$ converges to $0$ almost surely as $n\to\infty$ under (\ref{cond-bound}). On the event $\{\min_{u\in\T} V(u)\ge -\beta\}$, we see that necessarily $D_\infty^{(\beta)} = c_0 D_\infty$ almost surely, and $\sum_{u\in\mathcal Z[A]} R(\beta +V(u))\ee^{-V(u)}{\bf 1}_{\{V(u_k)\ge -\beta\}} \sim c_0 \sum_{u\in\mathcal Z[A]} (\beta+V(u))\ee^{-V(u)}$ as $A\to\infty$. Again by Theorem 6.1 (and Lemma 6.1) of \cite{biggins-kyprianou04}, we have $\lim_{A\to\infty} \sum_{u\in\mathcal Z[A]} \ee^{-V(u)}=W_\infty=0$ almost surely. We deduce that
\begin{equation}\label{eq:limitZA}
\lim_{A\to\infty} \sum_{u\in\mathcal Z[A]} V(u)\ee^{-V(u)} = D_\infty
\end{equation}

\noindent on the event $\{\min_{u\in\T} V(u)\ge -\beta\}$, and therefore almost surely by making $\beta\to\infty$. We can now prove the convergence in law.\\

{\it Proof of Theorem \ref{main}}.  Fix $x\in \r$ and  let $\varepsilon>0$.  For any $A>0$, we have for $n$ large enough,
\begin{eqnarray}
\P(\exists \, u\in \mathcal Z[A]\,:\, |u|\ge n^{1/2}) &\le& \varepsilon. \nonumber 
\end{eqnarray}

\noindent Take $A>0$. Let $\mathcal Y_A:=\{\max_{u\in\mathcal Z[A]}|u|\le n^{1/2}\}$. Let $\mathcal Z_n[A]:=\{u\in\mathcal Z[A]\,:\, V(u)\le \ln(n)\}$. We observe that, for $n$ large enough,
\begin{eqnarray*}
\P(M_n \ge (3/2) \ln n +x) &\ge & \P(M_n \ge (3/2) \ln n +x, \mathcal Y_A)\\
&=& \E\left[\prod_{u\in \mathcal Z[A]} (1-\Phi_{|u|,n}(V(u)-x)), \mathcal Y_A\right]
\end{eqnarray*}

\noindent where for any integers $n\ge 1$, $k\in[0,n]$ and any real $r\ge 0$,
$$
 \Phi_{k,n}(r) := \P(M_{n-k} < (3/2) \ln(n) - r).
$$

\noindent By Proposition \ref{p:tailmin}, there exists $A$ large enough and $N\ge 1$ such that for any $n\ge N$, $k \le n^{1/2}$ and $z\in [A-x,(3/2) \ln(n)-A-x]$,
\begin{equation}\label{eq:thmmain}
\Big|{\ee^z \over z}  \Phi_{k,n}(z) - C_1c_0\Big| \le \varepsilon.
\end{equation}

\noindent Therefore, for $n$ large enough, we have that if $u\in\mathcal Z_{n}[A]$, then on the event $\mathcal Y_A$,
$$
\Phi_{|u|,n}(V(u)-x) \le (C_1c_0+\varepsilon)(V(u)-x)\ee^{x-V(u)}.
$$

\noindent On the other hand, if $u\in \mathcal Z[A]\backslash \mathcal Z_{n}[A]$, then  Corollary \ref{c:infinite} implies that for $n$ large enough, on the event $\mathcal Y_A$, 
$$
\Phi_{|u|,n}(V(u)-x) \le c(1+V(u))\ee^{-V(u)}
$$
where the constant $c$ depends on $x$.  We get that, for $n$ large enough,
$$
\P(M_n \ge (3/2) \ln n +x)
\ge
\E\left[\prod_{u\in \mathcal Z[A]}(1- g_{u}(V(u))),\mathcal Y_A\right]
$$
where, if $u\in \mathcal Z_{n}[A]$, then $g_{u}(V(u)) =  (C_1c_0+\varepsilon)(V(u)-x)\ee^{x-V(u)})$, and if $u\in \mathcal Z[A]\backslash \mathcal Z_{n}[A]$, then $g_{u}(V(u)) = c (1+V(u))\ee^{-V(u)}$.

\noindent Since $\P(\mathcal Y_A^c)\le \varepsilon$ for $n$ large enough, we have for $n$ large enough
$$
\P(M_n \ge (3/2) \ln n +x)
\ge
\E\left[\prod_{u\in \mathcal Z[A]}(1- g_{u}(V(u)))\right] -\varepsilon.
$$
\noindent Notice that, almost surely,
$$
\lim_{n\to\infty} \prod_{u\in \mathcal Z[A]}(1- g_{u}(V(u))) = \prod_{u\in \mathcal Z[A]}(1- (C_1c_0+\varepsilon)(V(u)-x)\ee^{x-V(u)}).
$$

\noindent In particular, by dominated convergence,
$$
\liminf_{n\to \infty} \P(M_n \ge (3/2) \ln n +x) \ge \E\left[\prod_{u\in \mathcal Z[A]}(1- (C_1c_0+\varepsilon)(V(u)-x)\ee^{x-V(u)})\right]
-  \varepsilon.
$$

\noindent We let $A$ go to infinity. We have almost surely by (\ref{eq:limitZA}) and the fact that $\sum_{u\in\mathcal Z[A]} \ee^{-V(u)}$ vanishes
\begin{eqnarray}
\lim_{A\to\infty} \sum_{u\in \mathcal Z[A]} \ln (1- (C_1c_0+\varepsilon) (V(u)-x)\ee^{x-V(u)})
=
-(C_1c_0+\varepsilon) \ee^x D_\infty  .
\end{eqnarray}

\noindent By dominated convergence, we deduce that
$$
\liminf_{n\to \infty} \P(M_n \ge (3/2) \ln n +x) \ge \E\left[\exp(-(C_1c_0+\varepsilon) \ee^x D_\infty)\right]
-  \varepsilon,
$$

\noindent which gives the lower bound by letting $\varepsilon \to 0$.  The upper bounds works similarly. Let $A$ be such that (\ref{eq:thmmain}) is satisfied for $n\ge N$, $k \le n^{1/2}$ and $z\in [A-x,(3/2) \ln(n)-A-x]$. We observe that, for $n$ large enough,
\begin{eqnarray*}
\P(M_n \ge (3/2) \ln n +x)
&\le&
\P(M_n \ge (3/2) \ln n +x,\mathcal Y_A) + \varepsilon\\
&=&
\E\left[\prod_{u\in \mathcal Z_{n}[A]}(1-\Phi_{|u|,n}(V(u)-x)), \mathcal Y_A\right] +  \varepsilon.
\end{eqnarray*}

\noindent Using (\ref{eq:thmmain}), we end up with
$$
\limsup_{n\to \infty} \P(M_n \ge (3/2) \ln n +x)
\le
\E\left[\prod_{u\in \mathcal Z[A]}(1-(C_1c_0-\varepsilon)(V(u)-x)\ee^{x-V(u)})\right] +  \varepsilon.
$$

\noindent From here, we proceed as for the lower bound. \hfill $\Box$

\appendix

\section{The derivative martingale}
\label{s:derivative}

We work under (\ref{cond-bound}), (\ref{cond-2}) and (\ref{cond-ln}) but we drop the assumption that $\mathcal L$ is non-lattice. We recall from (\ref{def:R(x)}) that the renewal function $R(x)$ is defined by
$$
R(x)= \sum_{k\ge 0} \P(S_k \ge -x,\, S_k<\min_{0\le j\le k-1} S_j).
$$

\noindent The duality lemma says that $R(x)$ is also the expected  number of visits of the random walk $(S_n)_{n\ge 0}$ to the interval $(-x,0]$ before hitting $[0,\infty)$ (after time $1$). For any $\beta\ge 0$, we introduce for $n\ge 0$
$$
D_n^{(\beta)}:= \sum_{|u|=n} R(V(u)+\beta)\ee^{-V(u)} {\bf 1}_{\{V(u_k) \ge -\beta,\,\forall \,k\le n\}}.
$$

\noindent The following lemma is Lemma 10.2 in  \cite{biggins-kyprianou04}.  The analog in the case of the Brownian motion is Theorem 9 in \cite{kyp04}.
\begin{lemma}{\rm (\cite{biggins-kyprianou04})}
For any $\beta\ge 0$, the process $(D_n^{(\beta)},\, n\ge 0)$ is a non-negative martingale with respect to $(\mathscr F_n,n\ge 0)$. 
\end{lemma}
{\it Proof}. We recall that under $\P_a$, the branching random walk $(V(v),\,v\in\T)$ and the one-dimensional random walk $(S_k,\,k\ge 0)$ start at $a$. By the Markov property, we have
$$
\E\left[D_{n+1}^{(\beta)}\, |\, \mathscr F_n\right]= \sum_{|u|=n}   {\bf 1}_{\{V(u_k) \ge -\beta,\,\forall \,k\le n\}}  \E_{V(u)}\left[  \sum_{|v|=1} R(V(v)+\beta)\ee^{-V(v)} {\bf 1}_{\{ V(v)\ge -\beta \}}\right ].
$$

\noindent By (\ref{many-to-one}), we see that for any $u\in \T$ with $|u|=n$,
$$
\E_{V(u)}\left[  \sum_{|v|=1} R(V(v)+\beta)\ee^{-V(v)} {\bf 1}_{\{ V(v)\ge -\beta \}} \right]
=
\E_{V(u)}\left[  R(S_1+\beta)    {\bf 1}_{\{ S_1\ge -\beta \}} \right]\ee^{-V(u)}
$$

\noindent which is $R(V(u)+\beta)\ee^{-V(u)}$ by Lemma 1 of \cite{tanaka}. Therefore,
$$
\E\left[D_{n+1}^{(\beta)}\, |\, \mathscr F_n\right]= \sum_{|u|=n}   {\bf 1}_{\{V(u_k) \ge -\beta,\,\forall \,k\le n\}} R(V(u)+\beta)\ee^{-V(u)}
$$
which completes the proof.
 \hfill $\Box$

\bigskip

\noindent Since $(D_n^{(\beta)},\, n\ge 0)$ is a non-negative martingale, we can define for any $a\ge 0$ a probability measure $\Q_a^{(\beta)}$ on $\mathscr F_\infty$ such that for any $n\ge 1$,
\begin{equation}\label{def:Qbeta}
{d\Q^{(\beta)}_a \over d\P_a}\Big|_{\mathscr F_n}= {D_n^{(\beta)}\over R(a+\beta)\ee^{-a}}
\end{equation}

\noindent and we write as usual  $\Q^{(\beta)}$ for $\Q^{(\beta)}_0$, and $\hat{\E}^{(\beta)}_a$ (resp. $\hat{\E}^{(\beta)}$) for the expectation associated with $\Q_a^{(\beta)}$ (resp. $\Q^{(\beta)}$).  Let $\hat {\mathcal B}_a^{(\beta)}$ be the branching random walk with a spine defined as follows:
The spine $w_0^{(\beta)}$ starts at $V(w_0^{(\beta)})=a$. At time $1$ it gives birth to a point process distributed as $(V(x),|x|=1)$ under $\Q^{(\beta)}_a$. Then the spine element$w_1^{(\beta)}$ at time $1$ is chosen proportionally to $R(V(u)+\beta)\ee^{-V(u)}{\bf 1}_{\{V(u) \ge -\beta\}}$ among the children $u$ of $w_0^{(\beta)}$. At each time $n$, the spine element  $w_n^{(\beta)}$ produces an independent point process distributed as $(V(x),|x|=1)$ under $\Q^{(\beta)}_{V(w_n^{(\beta)})}$, while the other particles $|u|=n$ generate independent point processes distributed as $(V(x),|x|=1)$ under $\P_{V(u)}$. The spine $w_{n+1}^{(\beta)}$ at time $n+1$ is chosen proportionally to the weight $R(V(u)+\beta)\ee^{-V(u)}{\bf 1}_{\{V(u_k) \ge -\beta,\,\forall \,k\le n\}}$ among the children of $w_n^{(\beta)}$.  
We write $\hat{\mathscr F}_n^{(\beta)}$ for the $\sigma$-algebra obtained from $\mathscr F_n$ by including the information on the spine up to time $n$.  We write $\mathcal B_a^{(\beta)}$ for the (non-marked) branching random walk obtained from $\hat {\mathcal B}_a^{(\beta)}$ by ignoring the location of the spine, and note that  $\mathcal B_a^{(\beta)}$ is measurable with respect to $\mathscr F_\infty$. 
\begin{lemma}({\rm \cite{biggins-kyprianou04}})
The branching random walk under $\Q^{(\beta)}_a$ is distributed as $\mathcal B^{(\beta)}_a$.
\end{lemma}

\noindent {\it Proof}. We give a sketch of the proof. Let $n\ge 1$ and $T_n$ be a deterministic tree of height less than $n$.   We denote by $\T_{|n}$ the (random) tree $\T$ truncated at level $n$. Let $\P_{\hat{\mathcal B}^{(\beta)}_a}$ be a probability measure associated with ${\hat{\mathcal B}^{(\beta)}_a}$. We want to prove that the projection of  $\P_{\hat{\mathcal B}^{(\beta)}_a}$ on the space of non-marked branching random walks is $\Q_a^{(\beta)}$. Given deterministic infinitesimal intervals $(dz_u,u\in T_n)$, we compute that
\begin{eqnarray*}
&& \P_{\hat{\mathcal B}^{(\beta)}_a}(\T_{|n}=T_n, V(u) \in dz_u\, \forall \,u\in T_n)\\
&=&
\sum_{u\in T_n,|u|=n} \P_{\hat{\mathcal B}^{(\beta)}_a}(\T_{|n}=T_n, V(u) \in dz_u\, \forall \,u\in T_n ,\, w_n^{(\beta)}=u).
\end{eqnarray*}

\noindent For any $u\in T_n$ with $|u|=n$, we check that, by construction of our process $\hat{\mathcal B}_a^{(\beta)}$,
\begin{eqnarray*}
&& \P_{\hat{\mathcal B}^{(\beta)}_a}(  \T_{|n}=T_n, V(u) \in dz_u\, \forall \,u\in T_n ,\,  w_n^{(\beta)}=u   )\\
&=&
     \P_a (  \T_{|n}=T_n, V(u) \in dz_u\, \forall \,u\in T_n  )  {R(V(u_j)+\beta)\ee^{-V(u_j)}{\bf 1}_{\{\min_{j\le n}  V(u_j) \ge -\beta\}}  \over R(a+\beta)\ee^{-a}}
\end{eqnarray*}

\noindent where $u_j$ denotes the ancestor of $u$ in $T_n$ at generation $j$. Therefore,
\begin{eqnarray*}
 \P_{\hat{\mathcal B}^{(\beta)}_a}(\T_{|n}=T_n, V(u) \in dz_u\, \forall \,u\in T_n)
&=&
\E_a\left[{\bf 1}_{\{   \T_{|n}=T_n, V(u) \in dz_u\, \forall \,u\in T_n   \}} {D_n^{(\beta)} \over D_0^{(\beta)} }\right]
\end{eqnarray*}

\noindent which is $\Q_a^{(\beta)}(    \T_{|n}=T_n, V(u) \in dz_u\, \forall \,u\in T_n    )$ by definition. \hfill $\Box$

\bigskip

\noindent From now on, we will identify our branching random walk under $\Q^{(\beta)}_a$ with $\hat{\mathcal B}_a^{(\beta)}$. Notice that the proof shows that, for any vertex $u\in\T $ such that $|u|=n$,
$$
\Q_a^{(\beta)}(w_n^{(\beta)}=u\,|\, \mathscr F_n) = {  R(V(u)+\beta)\ee^{-V(u)} {\bf 1}_{ \{   \min_{j\le n} V(u_j) \ge -\beta\}} \over D_n^{(\beta)}  }.
$$

\noindent For $F$ a measurable function from $\mathbb{R}^{n+1}$ to $\mathbb{R}_+$, we notice that
\begin{eqnarray*}
&& \hat {\E}^{(\beta)}_a\left[F(V(w_0^{(\beta)}),\ldots,V(w_n^{(\beta)}))\right] \\
&=&
 \hat{ \E}^{(\beta)}_a\left[ \sum_{  |u|=n }  F(V(u_0),\ldots,V(u_n)) { R(V(u)+\beta)\ee^{-V(u)}{\bf 1}_{\{\min_{j\le n}V(u_j) \ge -\beta\} } \over D_n^{(\beta)}} \right] \\
 &=&
  {1\over D_0^{(\beta)}}\E_a\left[ \sum_{  |u|=n }  F(V(u_0),\ldots,V(u_n))  R(V(u)+\beta)\ee^{-V(u)}{\bf 1}_{\{\min_{j\le n}V(u_j) \ge -\beta\} }  \right] .
\end{eqnarray*}

\noindent Therefore, (\ref{many-to-one}) yields that
\begin{equation}\label{eq:Stanaka}
\hat{\E}^{(\beta)}_a\left[F(V(w_0^{(\beta)}),\ldots,V(w_n^{(\beta)}))\right] = {1\over R(a+\beta)}\E_a\left[F(S_0,\ldots,S_n)R(S_n+\beta),\, \min_{k\le n} S_k\ge -\beta\right].
\end{equation}

\noindent Under $\Q^{(\beta)}_y$, the spine process $(V(w_n),\,n\ge 0)$ is distributed as the random walk $(S_n)_{n\ge 0}$ conditioned to stay above $-\beta$, in the sense of \cite{tanaka} or \cite{bertoin-doney}. It is the Markov chain with transition probabilities, for any $x\ge -\beta$,
$$
\hat p^{(\beta)}(x,dy):= {R(y+\beta)\over R(x+\beta)} {\bf 1}_{\{ y\ge -\beta \}} p(x,dy)
$$

\noindent where $p(x,dy)=\P_x(S_1\in dy)$. The fact that this defines a transition probability comes from the equality $\E_x[R(S_1){\bf 1}_{\{ S_1\ge 0 \}}]=R(x)$ for any $x\ge 0$ by Lemma 1 in \cite{tanaka}. This Markov chain then never hits the region $(-\infty,-\beta)$, hence its name.

\bigskip

Since $(D_n^{(\beta)},\, n\ge 0)$ is a (non-negative) martingale, it has a limit that we denote by $D_\infty^{(\beta)}$. The question of the convergence in $L^1$ was adressed in \cite{biggins-kyprianou04}, where the authors give almost optimal conditions for the convergence to hold. However, we deal with slightly weaker conditions, so we have to prove the convergence in our case.

\begin{proposition}\label{l:Dbeta}
Assume (\ref{cond-bound}), (\ref{cond-2}) and (\ref{cond-ln}). Then, \\
(i) for any $\beta\ge 0$, $D_n^{(\beta)}$ converges in $L^1$ to $D_\infty^{(\beta)}$. \\
(ii) We have $D_\infty^{(\beta)}>0$ almost surely on the event of non-extinction of the branching random walk killed below $-\beta$.\\
(iii) We have $D_\infty>0$ almost surely on the event of non-extinction of $\T$.
\end{proposition}

\noindent {\it Proof}. We adapt the proof of \cite{biggins-kyprianou04} (see \cite{lyons} for the case of the additive martingale).  We observe that if $\sup_{n\to\infty} D_n^{(\beta)}<\infty$, $\Q^{(\beta)}$-a.s, then the family $(D_n^{(\beta)})_{n\ge 0}$ under $\P$ is uniformly integrable, hence converges in $L^1$.  Let 
$$
\hat {\mathcal G}_\infty^{(\beta)}:=\sigma\{ w_j^{(\beta)},V(w_j^{(\beta)}),\Omega(w_j^{(\beta)}),(V(u))_{u\in \Omega(w_j^{(\beta)})},\,j\ge 1\}
$$ 
be the $\sigma$-algebra of the spine and its brothers. Using the martingale property of $D_n^{(\beta)}$ for the subtrees rooted at brothers of the spine, we have
$$
\hat{\E}^{(\beta)}[D_n^{(\beta)}\,|\,\hat{\mathcal G}_\infty^{(\beta)}]
=
R(V(w_n^{(\beta)}) + \beta)\ee^{-V(w_n^{(\beta)})} + \sum_{k= 1}^n \sum_{x\in \Omega(w_k^{(\beta)})} R(V(x)+\beta)\ee^{-V(x)}{\bf 1}_{\{V(x_j)\ge -\beta,\,\forall\, j\le k\}}.
$$

\noindent It is well-known (see for example the construction available in \cite{tanaka} for the random walk conditioned to stay positive) that $V(w_n^{(\beta)})\to\infty$ $\Q^{(\beta)}$-almost surely, therefore $R(V(w_n^{(\beta)})+\beta)\ee^{-V(w_n^{(\beta)})}$ goes to zero as $n\to\infty$. Furthermore, we see that $1/D_n^{(\beta)}$ is under $\Q^{(\beta)}$ a positive supermartingale, and therefore converges as $n\to\infty$. We still denote by $D_\infty^{(\beta)}$ the (possibly infinite) limit of $D_n^{(\beta)}$ under $\Q^{(\beta)}$. We already know that there exists $c_{27}>0$ such that $R(x)\le c_{27}(1+x)_+\le c_{27}(1+x_+)$ for any $x\in\r$. Then, by Fatou's lemma
\begin{eqnarray}\nonumber
\hat {\E}^{(\beta)}[D_\infty^{(\beta)}\,|\, \hat {\mathcal G}_\infty^{(\beta)} ] &\le& \liminf_{n\to\infty} \hat {\E}^{(\beta)}[D_n^{(\beta)}\,|\,\hat {\mathcal G}_\infty^{(\beta)}]\\
&\le&
c_{27}\sum_{k\ge 1} \sum_{x\in \Omega(w_k^{(\beta)})} (1+(\beta + V(x))_+)\ee^{-V(x)}.\label{eq:A1A2a}
\end{eqnarray}

\noindent To prove (i), it remains to show that the right-hand side of the last inequality is finite $\Q^{(\beta)}$-almost surely (which implies that $D_\infty^{(\beta)}$ is finite $\Q^{(\beta)}$-a.s). We observe that
\begin{equation}\label{eq:A1A2b}
\sum_{k\ge 1} \sum_{x\in \Omega(w_k^{(\beta)})} (1+(\beta+V(x))_+)\ee^{-V(x)}
\le
A_1 + A_2
\end{equation}

\noindent with
\begin{eqnarray}
\label{def:A1} A_1 &:=& \sum_{k\ge 1}  (1+\beta+V(w_{k-1}^{(\beta)}))\ee^{-V(w_{k-1}^{(\beta)})}\sum_{x\in \Omega(w_k^{(\beta)})}\ee^{-(V(x)-V(w_{k-1}^{(\beta)}))},\\
\label{def:A2} A_2 &:=& \sum_{k\ge 1} \ee^{-V(w_{k-1}^{(\beta)})}\sum_{x\in \Omega(w_k^{(\beta)})} (V(x)-V(w_{k-1}^{(\beta)}))_+\ee^{-(V(x)-V(w_{k-1}^{(\beta)}))}.
\end{eqnarray}

\noindent Let us consider $A_1$. We recall that $X:=\sum_{|x|=1}\ee^{-V(x)}$, $\tilde X:=\sum_{|x|=1}V(x)_+ \ee^{-V(x)}$ and we introduce $X':=\sum_{|x|=1} R(\beta+V(x))\ee^{-V(x)}{\bf 1}_{\{V(x)\ge -\beta\}}$. We observe that, for any $a\ge -\beta$,
$$
X' \le c_{27} \sum_{|x|=1} \ee^{-V(x)} \left((1+ a +\beta) + (V(x)-a)_+\right).
$$

\noindent Therefore, we have for any $z\in\r$ and $a\ge -\beta$,
\begin{eqnarray}\nonumber
\Q_a^{(\beta)}\left(\sum_{|x|=1} \ee^{-(V(x)-a)}>z\right)
&=&
{1\over R(a+\beta)\ee^{-a}}\E_a\left[X'{\bf 1}_{\{ \sum_{|x|=1} \ee^{-(V(x)-a)} >z\}}\right]\\
\nonumber  &\le&
c_{40}\ee^a  \E_a\left[ \sum_{|x|=1} \ee^{-V(x)}   \left(1+ {(V(x)-a)_+\over 1+a +\beta} \right)   {\bf 1}_{\{\sum_{|x|=1} \ee^{-(V(x)-a)}>z\}}\right]\\
\label{eq:Qabeta} &=&
c_{40}\E\left[X{\bf 1}_{\{X>z\}}\right]
+
c_{40}{1\over 1+a+\beta}\E\left[\tilde X{\bf 1}_{\{X>z\}}\right]\\
\nonumber &=:&
c_{40}h_1(z) + c_{40}{1\over 1+a+\beta}h_2(z)
\end{eqnarray}

\noindent where $h_1$ and $h_2$ are defined by the last equation. We deduce by the Markov property at time $k-1$ that
\begin{eqnarray*}
&&\Q^{(\beta)}\left(\sum_{x\in \Omega(w_k^{(\beta)})}\ee^{-(V(x)-V(w_{k-1}^{(\beta)}))} \ge \ee^{V(w_{k-1}^{(\beta)})/2}\right)\\
&\le&
c_{40}\hat{\E}^{(\beta)}\left[h_1(\ee^{V(w_{k-1}^{(\beta)})/2}) + {1\over 1+V(w_{k-1}^{(\beta)})+\beta}h_2(\ee^{V(w_{k-1}^{(\beta)})/2})\right].
\end{eqnarray*}

\noindent Hence,
\begin{eqnarray}\label{eq:dnbeta1}
&&\sum_{k\ge 1} \Q^{(\beta)}\left(\sum_{x\in \Omega(w_k^{(\beta)})}\ee^{-(V(x)-V(w_{k-1}^{(\beta)}))} \ge \ee^{V(w_{k-1}^{(\beta)})/2}\right)\\
&\le&
c_{40}\sum_{\ell\ge 0} \hat{\E}^{(\beta)}\left[h_1(\ee^{V(w_{\ell}^{(\beta)})/2})\right]+ c_{40}\sum_{\ell\ge 0} \hat{\E}^{(\beta)}\left[{1\over 1+V(w_{\ell}^{(\beta)})+\beta}h_2(\ee^{V(w_{\ell}^{(\beta)})/2})\right]. \nonumber
\end{eqnarray}

\noindent We next estimate $\sum_{\ell\ge 0} \hat{\E}^{(\beta)}\left[h_1(\ee^{V(w_{\ell}^{(\beta)})/2})\right]$. By (\ref{eq:Stanaka}), we have
\begin{eqnarray*}
\hat{\E}^{(\beta)}\left[h_1(\ee^{V(w_{\ell}^{(\beta)})/2})\right]
&=&
{1\over R(\beta)}\E\left[ R(\beta+S_\ell)h_1(\ee^{S_{\ell}/2}),\, \min_{j\le \ell} S_j \ge -\beta\right]\\
&=&
{1\over R(\beta)}\E\left[ R(\beta+S_\ell)X {\bf 1}_{\{ S_{\ell}\le 2\ln X\}},\, \min_{j\le \ell} S_j \ge -\beta\right],\nonumber
\end{eqnarray*}

\noindent where $X$ and the random walk $(S_n,\,n\ge 0)$ are taken independent. Conditioning on $X$, then using Lemma \ref{l:estRW} (i), we get that
\begin{eqnarray}\nonumber
\sum_{\ell \ge 0} \hat{\E}^{(\beta)}\left[h_1(\ee^{V(w_{\ell}^{(\beta)})/2})\right] 
&\le&
{1\over R(\beta)}\E\left[ X R(\beta+ 2\ln(X)) \sum_{\ell\ge 0} {\bf 1}_{\{ S_{\ell}\le 2\ln X, \min_{j\le \ell} S_j \ge -\beta\}}\right]
\\
&\le& { c_{41} \over R(\beta)} \E[X(1+ \ln_+ X)^2)] \label{eq:h1}
\end{eqnarray}

\noindent which is finite by (\ref{cond-ln}). Similarly,
$$
\hat{\E}^{(\beta)}\left[ {1\over 1+V(w_\ell^{(\beta)})+\beta} h_2(\ee^{V(w_{\ell}^{(\beta)})/2})\right]
\le
c_{42}\E\left[\tilde X {\bf 1}_{\{S_\ell \le 2\ln X\}},\,\min_{j\le \ell} S_j\ge -\beta\right].
$$

\noindent Lemma \ref{l:estRW} (i) implies that
\begin{equation}\label{eq:h2}
\sum_{\ell \ge 0} \hat{\E}^{(\beta)}\left[ {1\over 1+V(w_\ell^{(\beta)})+\beta} h_2(\ee^{V(w_{\ell}^{(\beta)})/2})\right]
\le  c_{43}\E\left[ \tilde X (1+\ln_+ X)\right]<\infty
\end{equation}

\noindent under (\ref{cond-ln}) by Lemma \ref{l:estcond} (i). Equations (\ref{eq:dnbeta1}) , (\ref{eq:h1}) and (\ref{eq:h2}) give that
\begin{equation}\label{eq:borel}
\sum_{k\ge 1} \Q^{(\beta)}\left(\sum_{x\in \Omega(w_k^{(\beta)})}\ee^{-(V(x)-V(w_{k-1}^{(\beta)}))} \ge \ee^{V(w_{k-1}^{(\beta)})/2}\right)<\infty.
\end{equation}

\noindent By the Borel-Cantelli lemma, we obtain that
$$
(1+\beta+V(w_{k-1}^{(\beta)}))\ee^{-V(w_{k-1}^{(\beta)})}\sum_{x\in \Omega(w_k^{(\beta)})}\ee^{-(V(x)-V(w_{k-1}^{(\beta)})} \le (1+\beta+V(w_{k-1}^{(\beta)}))\ee^{-V(w_{k-1}^{(\beta)})/2}
$$

\noindent for $k$ large enough almost surely. It is known that, for any $a\in(0,1/2)$, we have $V(w_k^{(\beta)})\ge k^a$ for $k$ large enough. From (\ref{def:A1}), we deduce that $A_1<\infty$. We proceed similarly for $A_2$, replacing in (\ref{eq:Qabeta}) ${\bf 1}_{\{X>z\}}$ by ${\bf 1}_{\{\tilde X>z\}}$. By analogy, we find that $A_2<\infty$ if $\E[X(1+\ln_+ \tilde X)^2]$ and $\E[\tilde X(1+\ln_+ \tilde X)]$ are finite. This is the case by (\ref{cond-ln}) and Lemma \ref{l:estcond} (i). Equations (\ref{eq:A1A2a}) and (\ref{eq:A1A2b}) yield that $D_{\infty}^{(\beta)}<\infty$ $\Q^{(\beta)}$-a.s, which ends the proof of (i). We prove now (iii).  We see that, for any $x\in\T$ with $|x|=1$,
$$
D_\infty \ge \ee^{-V(x)}D_{\infty,x}\ge 0
$$
where for any $x\in\T$, $ D_{n,x}:=\sum_{|u|=n,u\ge x} (V(u)-V(x))\ee^{-(V(u)-V(x))} $ and $D_{\infty,x}:=\lim_{n\to\infty} D_{n,x}$. We used the fact that the martingale $\sum_{|u|=n,u\ge x}\ee^{-V(x)}$ converges to $0$ as $n\to\infty$. This implies that if $D_\infty=0$, then $D_{\infty,x}=0$. Notice that $D_{\infty,x}$ is distributed as $D_\infty$. Therefore, writing $p:=\P(D_\infty=0)$, we have that $p>0$  implies that $p\le \E[p^{\sum_{|x|=1} 1}]$. Consequently, $p=1$ or $p\le \P(\mbox{extinction of } \T)$. 
On the other hand, observe that $p \ge \P(\mbox{extinction of } \T)$, since the sum in (\ref{def:dWn}) is empty for large $n$ when the tree $\T$ is finite. Finally, we get that $\P(D_\infty=0)$ is $\P(\mbox{extinction of } \T)$ or $1$.  Now, notice that $\P(D_\infty^{(0)}>0)>0$ by (i). Since $R(x)\le c_{27}(1+x_+)$, we see that $D^{(0)}_\infty \le c_{27}D_\infty$, and therefore $\P(D_\infty>0)>0$. Hence, we have $D_\infty>0$ $\P$-a.s. on the event of non-extinction. We can now prove (ii). Let $\beta \ge 0$. On the event of non-extinction of the branching random walk killed below $\beta$, we can find a vertex $u$ (in the killed branching random walk) such that there is an infinite line of descent from $u$ which stays above $V(u)$. For such a vertex $u$, we have
$$
\sum_{v\ge u,|v|=n} R(V(v)+\beta)\ee^{-V(v)}{\bf 1}_{\{V(v_k)\ge -\beta,\, \forall k\le n\}}=\sum_{v\ge u,|v|=n} R(V(v)+\beta)\ee^{-V(v)}.
$$

\noindent The sum $\sum_{v\ge u,|v|=n} R(V(v)+\beta)\ee^{-V(v)}$ converges  to $c_0\ee^{-V(u)}D_{\infty,u}$ as $n\to\infty$ . We know from (iii) that $D_{\infty,u}>0$, hence $\sum_{v\ge u,|v|=n} R(\beta+V(v))\ee^{-V(v)}{\bf 1}_{\{V(v_k)\ge -\beta,\, \forall k\le n\}}$ has a positive limit as $n\to\infty$. Since $D_n^{(\beta)}\ge \sum_{v\ge u,|v|=n} R(\beta+V(v))\ee^{-V(v)}{\bf 1}_{\{V(v_k)\ge -\beta,\, \forall k\le n\}}$, we have that $D_\infty^{(\beta)}>0$. \hfill $\Box$

\section{Auxiliary estimates}

\begin{lemma}\label{l:estcond}
Let $X$ and $\tilde X$ be non-negative random variables such that (\ref{cond-ln}) holds. \\
(i) We have
$$
\E\left[X (\ln_+\tilde X)^2\right] < \infty, \qquad
\E\left[ \tilde X \ln_+ X\right] < \infty.
$$

\noindent (ii) As $z\to\infty$,
\begin{eqnarray*}
\E\left[ X (\ln_+(X+\tilde X))^2\min(\ln_+(X+\tilde X),z)\right]&=&o(z),\\
\E\left[ \tilde X \ln_+(X+\tilde X)\min(\ln_+(X+\tilde X),z)\right]&=& o(z).
\end{eqnarray*}
\end{lemma}

\noindent {\it Proof}. We first prove (i). We claim that for any $x,\tilde x\ge 0$
\begin{equation}\label{eq:estcond}
x(\ln_+ \tilde x)^2 \le  4 x(\ln_+ x)^2 + 2 \tilde x \ln_+\tilde x.
\end{equation}

\noindent We can assume that $\tilde x \ge 1$. If $\tilde x<x^2$, then $x(\ln_+ \tilde x)^2 \le  4 x(\ln_+ x)^2$. If $\tilde x\ge x^2$, we check that
$x (\ln_+ \tilde x)^2 \le 2 \tilde x \ln_+\tilde x$ since $\ln(y) \le 2 \sqrt{y}$ for any $y\ge 1$. This gives (\ref{eq:estcond}). It follows that
$$
\E\left[X(\ln_+ \tilde X)^2\right] \le  4\E\left[X(\ln_+ X)^2\right] + 2 \E\left[\tilde X \ln_+\tilde X\right]
$$

\noindent  which is finite under (\ref{cond-ln}).  Also, $\tilde X \ln_+ X\le \max(\tilde X\ln_+ \tilde X, X\ln_+ X)$, hence $\E[\tilde X\ln_+(X)]<\infty$. We turn to the proof of (ii). Let $\varepsilon>0$. We observe that
\begin{eqnarray*}
&&\E\left[ X (\ln_+(X+\tilde X))^2\min(\ln_+(X+\tilde X),z)\right]\\
&=&
\E\left[ X (\ln_+(X+\tilde X))^2\min(\ln_+(X+\tilde X),z),\, \ln_+ (X+\tilde X)\ge \varepsilon z\right] \\
&&+\, \E\left[ X (\ln_+(X+\tilde X))^2\min(\ln_+(X+\tilde X),z),\, \ln_+ (X+\tilde X)<\varepsilon z\right].
\end{eqnarray*}

\noindent On one hand,
\begin{eqnarray*}
&&\E\left[ X (\ln_+(X+\tilde X))^2\min(\ln_+(X+\tilde X),z),\, \ln_+ (X+\tilde X)\ge \varepsilon z\right]\\
&\le&
z\E\left[ X (\ln_+(X+\tilde X))^2,\, \ln_+ (X+\tilde X)\ge \varepsilon z\right]\\
&=&
z o_z(1)
\end{eqnarray*}

\noindent since $\E\left[ X (\ln_+(X+\tilde X))^2\right]<\infty$. On the other hand,
$$
\E\left[ X (\ln_+(X+\tilde X))^2\min(\ln_+(X+\tilde X),z),\, \ln_+ (X+\tilde X)<\varepsilon z\right]
\le
\varepsilon z \E\left[X (\ln_+ X+\tilde X)^2\right].
$$

\noindent Thus, $\E\left[ X (\ln_+(X+\tilde X))^2\min(\ln_+(X+\tilde X),z)\right]\le (1+\E[X (\ln_+ X+\tilde X)^2] )\varepsilon z$ for $z$ large enough, and is therefore $o(z)$. We show similarly that $\E\left[ \tilde X \ln_+(X+\tilde X)\min(\ln_+(X+\tilde X),z)\right]=o(z)$. \hfill $\Box$ \\

Let $(S_n)_{n\ge 0}$ be a one-dimensional random walk, with $\E[S_1]=0$ and $\E[(S_1)^2]<\infty$.
\begin{lemma}\label{l:estRW}
(i) There exists a constant $c_{45}>0$ such that for any  $z\ge 0$ and $x\ge 0$
$$
\sum_{\ell \ge 0} \P_z\left( S_{\ell} \le x,\, \min_{j\le \ell} S_j \ge 0 \right)
\le c_{45} (1+ x)(1+\min(x,z)).
$$
(ii) Let $a>0$. We have
$$
\E\left[\sum_{\ell\ge 0} \ee^{-a S_\ell}{\bf 1}_{\{\min_{j\le \ell} S_j \ge 0\}} \right] =c_{46}(a)<\infty .
$$
(iii) Let $a>0$. There exists a constant $c_{47}(a)>0$ such that for any $z\ge 0$,
$$
\E_z\left[\sum_{\ell\ge 0} \ee^{-a S_{\ell}}{\bf 1}_{\{\min_{j\le \ell}S_j\ge 0\}}\right] \le c_{47}(a).
$$
\end{lemma}

\noindent {\it Proof}. Suppose that $x< z$. If $\tau_{x}^-$ denotes the first passage time below level $x$ of $(S_n)_{n\ge 0}$, we have
\begin{eqnarray*}
\sum_{\ell\ge 0}\P_z\left( S_{\ell} \le x,\, \min_{j\le \ell} S_j \ge 0 \right)
&=&
\E_z\left[\sum_{\ell\ge \tau_x^-}{\bf 1}_{\{ S_{\ell} \le x,\, \min_{j\le \ell} S_j \ge 0}\} \right]\\
&\le&
\E\left[\sum_{\ell\ge 0}{\bf 1}_{\{ S_{\ell} \le x,\, \min_{j\le \ell} S_j \ge -x\}} \right]
\end{eqnarray*}

\noindent where we used the Markov property at time $\tau_{x}^-$. We have
\begin{eqnarray}
\nonumber \sum_{\ell \ge 0} \P\left( S_{\ell} \le x,\, \min_{j\le \ell} S_j \ge -x \right)
&\le&
1+x^2 +\sum_{\ell >x^2} \P\left( S_{\ell} \le x,\, \min_{j\le \ell} S_j \ge -x \right)\\
\nonumber &\le& 1+x^2+c_{48}\sum_{\ell >x^2} (1+x)^3 \ell^{-3/2}\\
\label{eq:estRW1}&\le& c_{49}(1+x)^2
\end{eqnarray}

\noindent by (\ref{lemmaA1}). Suppose now that $x\ge z$. Then,
\begin{eqnarray*}
&&\sum_{\ell \ge 0} \P_z\left( S_{\ell} \le x,\, \min_{j\le \ell} S_j \ge 0 \right)\\
&\le&
\sum_{\ell \le x^2} \P_z\left(\min_{j\le \ell} S_j \ge 0 \right) +\sum_{\ell >x^2} \P_z\left( S_{\ell} \le x,\, \min_{j\le \ell} S_j \ge 0 \right).
\end{eqnarray*}

\noindent From (\ref{eq:kozlov}), we know that $\P_z\left(\min_{j\le \ell} S_j \ge 0 \right) \le c_{50} (1+z)(1+\ell)^{-1/2}$, whereas, by (\ref{lemmaA1}), 
$$\P_z\left( S_{\ell} \le x,\, \min_{j\le \ell} S_j \ge 0 \right)\le c_{51}(1+z)(1+x)^2 (1+\ell)^{-3/2}.$$

\noindent We get
\begin{eqnarray}
\nonumber \sum_{\ell \ge 0} \P_z\left( S_{\ell} \le x,\, \min_{j\le \ell} S_j \ge 0 \right)
&\le& c_{50} \sum_{\ell \le x^2}{1+z \over \sqrt{1+\ell}} + c_{51} \sum_{\ell >x^2} (1+z)(1+x)^2 (1+\ell)^{-3/2}\\
\label{eq:estRW2} &\le& c_{52}(1+z)(1+x).
\end{eqnarray}

\noindent From (\ref{eq:estRW1}) when $x<z$ and (\ref{eq:estRW2}) when $x\ge z$, we have for $x,z\ge 0$,
$$
\sum_{\ell \ge 0} \P_z\left( S_{\ell} \le x,\, \min_{j\le \ell} S_j \ge 0 \right)
\le (c_{49} + c_{52})(1+x)(1+\min(x,z)).
$$

\noindent  This ends the proof of (i). We turn to the statement (ii). Without loss of generality, we assume that $a=1$ (in (ii) and in (iii)). We have 
\begin{eqnarray*}
\sum_{\ell\ge 0}\E\left[ \ee^{-S_\ell}{\bf 1}_{\{\min_{j\le \ell} S_j \ge 0\}} \right]
=
\sum_{\ell\ge 0} \sum_{i\ge 0} \ee^{-i} \P(S_{\ell} \in[i,i+1),\, \min_{j\le \ell} S_j \ge 0).
\end{eqnarray*}

\noindent By (\ref{lemmaA1}), $\P(S_{\ell} \in[i,i+1),\, \min_{j\le \ell} S_j \ge 0) \le c_{53}(1+i)(1+\ell)^{-3/2}$, which completes the proof of (ii). Finally, we prove (iii). Let $( T_k,  H_k,k\ge 0)$ be the strict descending ladder epochs and heights of $(S_n)_{n\ge 0}$, i.e. $ T_0:=0$, $ H_0:=S_0$ and for any $k\ge 1$, $ T_k:=\min\{j>  T_{k-1}\,:\, S_j< H_{k-1}\}$, $H_k:=S_{ T_k}$.  By applying the Markov property at the times $( T_k,\,k\ge 0)$, we observe that
$$
\E_z\left[\sum_{\ell\ge 0} \ee^{-S_{\ell}} {\bf 1}_{\{\min_{j\le \ell} S_j\ge 0\}}\right]
=
c_{46} \E_z\left[\sum_{k\ge 0} \ee^{-H_k}{\bf 1}_{\{H_k\ge 0\}} \right]
$$

\noindent where $c_{46}$ is the constant of (ii). The fact that $Z(z):=\E_z\left[\sum_{k\ge 0} \ee^{- H_k}{\bf 1}_{\{ H_k\ge 0\}} \right]$ is bounded in $z\ge 0$ then comes from the renewal theorem: let $U(dy)$ denote the renewal measure of $( H_k,\,k\ge 0)$, i.e. $U(dy):= \sum_{k\ge 0} \P( H_k \in dy)$. Then $Z(z) = \int_{-z}^0 \ee^{-(z+y)}U(dy)$. In  Section XI.1 of \cite{feller}, combine Lemma p.359 with the renewal Theorem p.363 to conclude that $Z$ is bounded. This completes the proof of (iii). \hfill $\Box$

\bigskip

For $\alpha>0$, $a\ge 0$, $n\ge 1$ and $0\le i\le n$, we define
\begin{equation}\label{defk}
k_i :=
\begin{cases}
i^{\alpha} , &\hbox{if $0\le i\le \lfloor n/2\rfloor$,} \cr
a+(n-i)^{\alpha} , &\hbox{if $\lfloor n/2 \rfloor <i\le n$.}\cr
\end{cases}	
\end{equation}

\begin{lemma}\label{l:a6}
Let $\alpha \in (0,1/6)$ and $\varepsilon>0$. \\
(i) There exist $d>0$ and  $c_{53}>0$ such that for any $u\ge 0$, $a\ge 0$ and any integer $n\ge 1$,
\begin{eqnarray}\label{eq:a6}
&&\P\Big\{\exists\, 0\le i\le   n :\, S_i\le  k_i-d,\, \min_{j\le n} S_j \ge 0, \;
 \min_{\lfloor  n/2\rfloor < j\le n} S_j \ge a, \;S_n \le a +u\Big\}\\
\nonumber &\le& (1+u)^2\Big\{{\varepsilon\over n^{3/2}} + c_{53}{(n^{\alpha} + a)^2\over n^{2-\alpha}}\Big\},
\end{eqnarray}

\noindent where $k_i$ is given by (\ref{defk}). 
\end{lemma}

\noindent {\it Proof}. We treat $n/2$ as an integer. Let $E$ be the event in (\ref{eq:a6}). We have
$\P(E)\le \sum_{i=1}^n \P(E_i)$ where
$$
E_i :=\{  S_i\le  k_i-d,\, \min_{j\le n} S_{j} \ge 0, \;
 \min_{ n/2 <j\le n} S_j \ge a, \;S_n \le a +u\}.
$$

\noindent We first treat the case $i\le n/2$, so that $k_i=i^{\alpha}$. By the Markov property at time $i\ge 1$ and  (\ref{lemmaA3}), we have
$$
\P(E_i) \le {c_{54}(1+u)^2\over n^{3/2}}\E\left[ (1+S_i){\bf 1}_{\{  S_i\le  i^{\alpha},\, \min_{j\le i} S_j \ge 0\}}\right]
$$

\noindent which is smaller than ${c_{55}(1+u)^2\over n^{3/2}} \, \,{(1+i^{\alpha})^3\over i^{3/2}}$ by (\ref{lemmaA1}). It yields that, if $K$ is greater than some constant $K_0$ (which does not depend on $d$), we have
\begin{equation}\label{eq:a6sum1}
\sum_{i=K}^{n/2} \P(E_i) \le (1+u)^2{\varepsilon \over n^{3/2}},
\end{equation}

\noindent [$\sum_{i=x}^y :=0$ if $x>y$.] We treat the case $n/2 < i\le n$. We have by the Markov property at time $i$ and (\ref{lemmaA1}),
$$
\P(E_i) \le {c_{56}(1+u)^2\over (n-i+1)^{3/2}}\E\left[ (1+S_i-a){\bf 1}_{\{  S_i\le  a+(n-i)^{\alpha},\, \min_{j\le i} S_j \ge 0,\,\min_{  n/2 <j\le i} S_j \ge a \}}\right].
$$

\noindent If $i\ge  2n/3$, we use (\ref{lemmaA3}) to see that $ \P(E_i) \le c_{57}(1+u)^2{(1+n-i)^{3\alpha-{3\over 2}}\over n^{3/2}}$. Therefore, if $K\ge K_1$, ($K_1$ does not depend on $d$),
\begin{equation}\label{eq:a6sum2}
\sum_{i=\lfloor 2n/3 \rfloor}^{n-K} \P(E_i) \le (1+u)^2{\varepsilon \over n^{3/2}}.
\end{equation}

\noindent If $n/2< i< 2 n/3$, we simply write
\begin{eqnarray*}
\P(E_i)
&\le&
{c_{56}(1+u)^2\over (n-i+1)^{3/2}}\E\left[ (1+S_i-a){\bf 1}_{\{  a\le S_i\le  a+(n-i)^{\alpha},\, \min_{j\le i} S_j \ge 0\}}\right]\\
&\le& c_{59} (1+u)^2 {(n-i)^{\alpha} \over (n-i+1)^{3/2}} \P(a\le S_i\le  a+(n-i)^{\alpha},\, \min_{j\le i} S_j \ge 0)\\
&\le&
c_{60} (1+u)^2 {n^{\alpha}(a+n^{\alpha})^2\over n^3}
\end{eqnarray*}

\noindent by (\ref{lemmaA1}). We deduce that
\begin{equation}\label{eq:a6sum3}
\sum_{i=n/2}^{\lfloor 2n/3\rfloor} \P(E_i) \le c_{61}(1+u)^2{ (n^{\alpha} +a)^2 \over n^{2-\alpha}}.
\end{equation}

\noindent Notice that our choice of $K$ does not depend on the constant $d$. Thus, we are allowed to choose  $d\ge K^{\alpha}$, for which  $\P(E_i)=0$ if $i\in [1,K]\cup [n-K,n]$. We obtain by (\ref{eq:a6sum1}),(\ref{eq:a6sum2}) and (\ref{eq:a6sum3})
$$
\sum_{i=1}^n \P(E_i) \le (1+u)^2\Big\{2{\varepsilon \over n^{3/2}}+c_{61}{ (n^{\alpha} +a)^2 \over n^{2-\alpha}}\Big\},
$$

\noindent hence $\P(E)\le (1+u)^2\Big\{2{\varepsilon \over n^{3/2}}+c_{61}{ (n^{\alpha} +a)^2 \over n^{2-\alpha}}\Big\}$ indeed. \hfill $\Box$

\section{The good vertex}
\label{s:goodvertex}

  Let $z\ge 0$ and $L\ge 0$. 
Let $d_k=d_k(n,z+L,1/2)$ as defined in (\ref{def:dk}). Let also
$$
e_k=e_k^{(n)}:=
\begin{cases}
k^{1/12}, &\hbox{if $0\le k\le {n\over 2}$,} \cr
(n-k)^{1/12}, &\hbox{if ${n\over 2}<k\le n$.}\cr
\end{cases}
$$

\noindent We recall from definition \ref{def:ZzL} that $u \in \mathcal Z^{z,L}_n$ if $|u|=n$, $V(u_k)\ge d_k$ for $k\le n$ and $V(u)\in I_n(z)$. We say that $u$ such that $|u|=n$ is a $(z,L)$-good vertex if $u\in\mathcal Z^{z,L}_n$ and for any $1\le k\le n$,
\begin{equation}\label{eq:goodappendix}
\sum_{v\in \Omega(u_k)} \ee^{-(V(v)-d_k)}\Big\{ 1+ (V(v)-d_k)_+ \Big\}\le B\ee^{-e_k}.
\end{equation}

\noindent Note that a $(z,0)-$good vertex is a $z$-good vertex as introduced in Section \ref{ss:tightness}. We defined the probability $\Q$ in (\ref{defQ}) and the spine $(w_n,n\ge 0)$ in Section \ref{s:spine}.

\begin{lemma}\label{l:goodvertex}
Fix $L\ge 0$.  For any $\varepsilon>0$, we can find $B$ large enough in (\ref{eq:goodappendix}) such that $\Q(w_n \mbox{ is not a $(z,L)-$good vertex},\,w_n \in \mathcal Z^{z,L}_n)\le \varepsilon n^{-3/2}$ for any $n\ge 1$ and $z\ge 0$.
\end{lemma}

\noindent {\it Proof}. Fix $L\ge 0$ and let $\varepsilon>0$.  We  have
  \begin{eqnarray}\label{eq:goodvertex}
&&\Q(w_n \mbox{ is not a $(z,L)-$good vertex},\,w_n \in \mathcal Z^{z,L}_n) \\
&\le&
 \Q\left(\exists \, k\in[1,n]\,:\,  \sum_{v\in \Omega(w_k)} \ee^{-(V(v)-d_k)}\Big\{ 1+ (V(v)-d_k)_+ \Big\}> B \ee^{-e_k} , \, w_n\in\mathcal Z^{z,L}_n \right)  \nonumber.
\end{eqnarray}

\noindent We want to show that we can find $B$ large enough such that
\begin{equation}\label{eq:3n4a}
 \Q\left(\exists \, k\in[1,n]\,:\,  \sum_{v\in \Omega(w_k)} \ee^{-(V(v)-d_k)}\Big\{ 1+ (V(v)-d_k)_+ \Big\}> B \ee^{-e_k} , \, w_n\in\mathcal Z^{z,L}_n \right) \le {\varepsilon \over n^{3/2}}.
\end{equation}

 \noindent We see that, for any $1\le k\le n$,
\begin{eqnarray*}
&&\left\{ \sum_{v\in \Omega(w_k)} \ee^{-(V(v)-d_k)}\Big\{ 1+ (V(v)-d_k)_+ \Big\}> B\ee^{-e_k},\, V(w_{k-1}) \ge d_k+2 e_k-c_{62}\right\}\\
&\subset&
\left\{ \sum_{v\in \Omega(w_k)} \ee^{-(V(v)-d_k)}\Big\{ 1+ (V(v)-d_k)_+ \Big\}> B\ee^{- {V(w_{k-1})-d_{k}+c_{62}\over 2}}\right\}.
\end{eqnarray*}

 \noindent By Lemma \ref{l:a6}, there exists $c_{62}=c_{62}(L)>0$ and $N=N(L)$ such that for $n\ge N$ and $z\ge 0$
$$
\Q\left( w_n \in \mathcal Z^{z,L}_n,\, \exists\, 0\le j\le n-1\,:\, V(w_j) \le d_{j+1} + 2 e_{j+1}-c_{62} \right) \le { \varepsilon \over n^{3/2}}.
$$

\noindent Consequently, it is enough to show that for $B$ large enough,
\begin{equation}\label{eq:gooda-2}
\sum_{k=1}^{n }  \Q\left(  \sum_{v\in \Omega(w_k)} \ee^{-(V(v)-d_k)}\Big\{ 1+ (V(v)-d_k)_+ \Big\}> B \ee^{-{V(w_{k-1})-d_k\over 2}} , \, w_n\in\mathcal Z^{z,L}_n \right) \le \varepsilon n^{-3/2}.
\end{equation}

 \noindent We see that
\begin{eqnarray*}
&&\sum_{v\in \Omega(w_k)} \ee^{-(V(v)-d_k)}\Big( 1+ (V(v)-d_k)_+ \Big)  \\
&\le&
\ee^{-(V(w_{k-1})-d_k)}\sum_{v\in \Omega(w_k)} \ee^{-(V(v)-V(w_{k-1}))}\Big\{ 1+ (V(w_{k-1})-d_k)_+ + (V(v)-V(w_{k-1}))_+ \Big\}\\
&\le&
\ee^{-(V(w_{k-1}) -d_k)}( 1+ (V(w_{k-1})-d_k)_+)\sum_{v\in \Omega(w_k)} \ee^{-(V(v)-V(w_{k-1}))}\Big\{1 + (V(v)-V(w_{k-1}))_+ \Big\}.
\end{eqnarray*}

\noindent With the notation of (\ref{def:xi}), we have then
$$
\sum_{v\in \Omega(w_k)} \ee^{-(V(v)-d_k)}\Big( 1+ (V(v)-d_k)_+ \Big)
\le 
\ee^{-(V(w_{k-1})-d_k)}( 1+  (V(w_{k-1})-d_k)_+) \xi(w_k).
$$

\noindent Equation (\ref{eq:gooda-2}) boils down to showing that, for $B$ large enough,
$$
\sum_{k=1}^{n }  \Q\left(  \xi(w_k) > B {\ee^{{V(w_{k-1})-d_k\over 2}} \over  1+  (V(w_{k-1})-d_k)_+}  , \, w_n\in\mathcal Z^{z,L}_n \right) \le \varepsilon n^{-3/2}.
$$

 \noindent Actually, we are going to  show that, for $B$ large enough,
\begin{equation}\label{eq:n}
\sum_{k=1}^{n }  \Q\left(  \xi(w_k) > B \ee^{(V(w_{k-1})-d_k)/3}   , \, w_n\in\mathcal Z^{z,L}_n \right) \le \varepsilon n^{-3/2}.
\end{equation}

\noindent First, we deal with the case $k\in[1,3n/4]$. We notice that
$$
 \Q\left(  \xi(w_k) > B \ee^{(V(w_{k-1})-d_k)/3}   , \, w_n\in\mathcal Z^{z,L}_n \right) \le \Q\left(  \xi(w_k) > B \ee^{V(w_{k-1})/3}   , \, w_n\in\mathcal Z^{z,L}_n \right) .
$$

\noindent By the Markov property at time $k$, we get
$$
   \Q\left(  \xi(w_k) > B \ee^{V(w_{k-1})/3}   , \, w_n\in\mathcal Z^{z,L}_n \right)     
=
\hat \E\left[\lambda(V(w_k),k,n){\bf 1}_{\{\xi(w_k) > B\ee^{V(w_{k-1})/3},\, V(w_j)\ge 0,\, \forall \, j\le k\}}  \right]
$$

\noindent where $\lambda(r,k,n):= \Q_r( V(w_j)\ge d_{j+k},\, \forall\, j\le n-k,\, V(w_{n-k}) \in I_n(z))$. We get by (\ref{lemmaA3}), $\lambda(r,k,n)\le c_{63} n^{-3/2}(1+r_+)$ (since $k\le 3n/4)$. This yields that
\begin{eqnarray}\label{eq:goodappendix1}
&& \Q\left(  \xi(w_k) > B \ee^{V(w_{k-1})/3}   , \, w_n\in\mathcal Z^{z,L}_n \right)\\
&\le&
c_{63} n^{-3/2}  \hat{\E}\left[(1+V(w_k)_+) {\bf 1}_{\{ \xi(w_k) > B\ee^{V(w_{k-1})/3},\, V(w_j)\ge 0,\,\forall\, j\le k\}}\right].\nonumber
\end{eqnarray}

\noindent On the other hand, we have
$$
1+V(w_k)_+ \le 1+ V(w_{k-1})_+ + (V(w_k) - V(w_{k-1}))_+ .
$$

\noindent Let  $(\xi,\Delta)$ be generic random variables distributed as $\Big(\xi(w_1),V(w_1)_+\Big)$ under $\Q$, and independent of the other random variables. 
By the Markov property at time $k-1$, we obtain that
$$
 \hat{\E}\left[  (1+V(w_k)_+) {\bf 1}_{\{ \xi(w_k) > B\ee^{V(w_{k-1})/3},\, V(w_j)\ge 0,\,\forall\, j\le k\}}  \right]
\le
\hat{\E}\left[ \kappa (V(w_{k-1})) {\bf 1}_{\{ V(w_j) \ge 0,\, \forall \, j\le k-1 \}} \right]
$$

\noindent with, for $x\ge 0$, $\kappa(x):=(1+x){\bf 1}_{\{\xi > B\ee^{x/3}\}} + \Delta_+ {\bf 1}_{\{ \xi> B\ee^{x/3}\}}$. In view of (\ref{eq:goodappendix1}), it follows that
\begin{eqnarray*}
\sum_{k=1}^{3n/4} \Q\left( \xi(w_k) > B \ee^{V(w_{k-1})/3}  ,\,w_n\in\mathcal Z^{z,L}_n \right)
\le
c_{63}n^{-3/2}(D_1 + D_2)
\end{eqnarray*}

\noindent where
\begin{eqnarray*}
D_1 &:=& \sum_{k\ge 0} \hat{\E}\left[ (1+V(w_{k})){\bf 1}_{\{V(w_{k}) \le 3 (\ln \xi -\ln  B)\}},\, \min_{j\le k} V(w_j)\ge 0\right]\\
D_2 &:=& \sum_{k\ge 0} \hat{\E}\left[ \Delta_+ {\bf 1}_{\{V(w_{k}) \le 3(\ln \xi  -\ln  B)\}},\, \min_{j\le k} V(w_j) \ge 0\right].
\end{eqnarray*}

\noindent We recall that by Proposition \ref{p:spine} $(V(w_n),\,n\ge 0)$ is distributed as $(S_n,\, n\ge 0)$ (under $\P$). Notice that in the definition of $D_1$, the term inside the expectation is $0$ if $ B> \xi$. Therefore, we can add the indicator that $ B\le \xi$. By Lemma \ref{l:estRW} (i), we get that
$$
D_1 \le c_{65}\hat{\E}\left[{\bf 1}_{\{ B\le \xi \}}(1+ (\ln\xi -\ln  B)_+)^2\right]\le c_{65} \hat{\E}\left[{\bf 1}_{\{ B\le \xi \}}(1+ \ln_+\xi)^2\right].
$$

\noindent Observe that $\xi \le X + \tilde X$ with the notation of (\ref{def:X}). Going back to the measure $\P$, we get
$$
D_1 \le  c_{65}\E\left[ X {\bf 1}_{\{ B\le X+\tilde X \}}(1+ \ln_+(X+\tilde X))^2\right]\le \varepsilon
$$

\noindent for $B$  large enough since $\E\left[ X (1+ \ln_+(X+\tilde X))^2\right]<\infty$ by (\ref{cond-ln}) and Lemma \ref{l:estcond} (i). Similarly,
$$
D_2 \le c_{66}\E\left[ \tilde X {\bf 1}_{\{ B\le X+\tilde X \}}(1+ \ln_+(X+\tilde X))\right]\le \varepsilon
$$

\noindent for $B$ large enough. Therefore, for $B$ large enough
\begin{equation}\label{eq:3n4}
\sum_{k=1}^{3n/4} \Q\left( \xi(w_k)> B\ee^{-(V(w_{k-1})-d_k)/3},\,w_n\in\mathcal Z^{z,L}_n \right)
\le
2{ \varepsilon \over n^{3/2}}.
\end{equation}

\noindent In order to prove (\ref{eq:n}), it remains to treat the case $3n/4\le k\le n$. We want to show that for $B$ large enough,
\begin{equation}\label{eq:nb}
 \sum_{k=3n/4}^n \Q\left( \xi(w_k) >B \ee^{(V(w_{k})-d_k)/3}, w_n\in\mathcal Z^{z,L}_n \right) \le \varepsilon.
\end{equation}

\noindent We want to condition the point process $\mu(w_1):=\sum_{u\in \I(w_1)}\delta_{V(u)}$ on the value of $V(w_1)$. To do this, we make a disintegration (see for example 15.3.3 pp. 164 of \cite{kallenberg}). This gives the existence of probabilities ${\bf Q}_r$ on the space of locally finite measures $\mathcal{M}$ on $\mathbb{R}$, such that:
\begin{itemize}
\item for any set $A$ in the canonical $\sigma$-algebra of $\mathcal{M}$, the map $r\in\mathbb{R} \to  {\bf Q}_r(A)$ is measurable with respect to the Borelian $\sigma$-algebra of $\mathbb{R}$. Here, the canonical $\sigma$-algebra of $\mathcal{M}$ refers to the one generated by the mappings $\mu \in\mathcal{M}\to \mu(I)$ for $I$ intervals of $\mathbb{R}$ (see chapter 1 of \cite{kallenberg}). 
\item For any bounded measurable function $F$, we have 
$$\hat{\E}\left[F\left(\mu(w_1),V(w_1)\right)\right]= \int_{\mathbb{R}} \Q(V(w_1)\in dr)\int_{{\mathcal M}} F\left(\mu,r\right) {\bf Q}_r(d\mu)  .$$
\end{itemize}

\noindent We deduce that
\begin{eqnarray*}
&& \Q\left(\xi(w_k) > B \ee^{(V(w_{k})-d_k)/3} ,\,w_n\in\mathcal Z^{z,L}_n \right)\\
&=&
  \Q\left({\overline \xi}\left(V(w_k)-V(w_{k-1})\right) > B \ee^{(V(w_{k})-d_k)/3} ,\,w_n\in\mathcal Z^{z,L}_n \right)
\end{eqnarray*}

\noindent where, given $(V(w_k),\,k\le n)$, the random variable ${\overline \xi}(V(w_k)-V(w_{k-1}))\in\r$ \ has the distribution of $\int_{x\in\r} (1+x_+)\ee^{-x}\mu(dx)$ under ${\bf Q}_{V(w_k)-V(w_{k-1})}(d\mu)$. The last line is equal to
\begin{eqnarray*}
\P\left({\overline \xi}\left(S_k-S_{k-1}\right) > B \ee^{(S_k-d_k)/3},\,S_n\in I_n(z),\,{\underline S}_n \ge 0,\, {\underline S}_{(n/2,n]}\ge a_n(z+L+1) \right)
\end{eqnarray*}

\noindent where ${\underline S}_n:= \min\{ S_k,\,k\le n \}$, ${\underline S}_{(\ell_1,\ell_2]}:= \min\{ S_k,\,\ell_1<k\le \ell_2 \}$ and, under $\P$, and conditionally on  $(S_k,\,k\le n)$, the random variable ${\overline \xi}(S_k-S_{k-1})$ has the distribution of $\int_{x\in\r} (1+x_+)\ee^{-x}\mu(dx)$ under ${\bf Q}_{S_k-S_{k-1}}(d\mu)$.  We return time, that is we replace $S_k$ by $S_n-S_{n-k}$. We check that
\begin{eqnarray*}
&& \P\left({\overline \xi}\left(S_k-S_{k-1}\right) > B\ee^{(S_k-d_k)/3},\,S_n\in I_n(z),\,{\underline S}_n \ge 0,\, {\underline S}_{(n/2,n]}\ge a_n(z+L+1) \right)\\
&\le&
 \P\Big({\overline \xi}\left(S_{n-k+1}-S_{n-k}\right) > B\ee^{L+1}\ee^{-S_{n-k}/3} ,\,S_n\in I_n(z), {\underline {-S}_n} \ge -a_n(z),\, {\underline {-S}}_{[0,n/2)}\ge -L-1 \Big)
\end{eqnarray*}

\noindent where ${\underline {-S}}_n:= \min\{ -S_k,\,k\le n \}$ and ${\underline {-S}}_{[\ell_1,\ell_2)}:= \min\{ -S_k,\,\ell_1 \le k< \ell_2 \}$. We use the Markov property at time $n-k+1$. There exists a constant $c_{67}>0$ such that, for any $r\le L+1$, any $n\ge 1$, and any $k\in [3n/4,n]$,
$$
\P_r\left( S_{k-1}\in I_n(z),\,  {\underline {-S}_{k-1}} \ge -a_n(z),\, {\underline {-S}}_{[0,k-1-{n\over 2})}\ge -L-1  \right) \le c_{67}(2+L-r)n^{-3/2}.
$$

\noindent The last inequality comes from (\ref{lemmaA3}), after a time reversal. This yields that, for any $n\ge 1$ and $k\in [3n/4,n]$,
\begin{eqnarray*}
&& 
\P\left({\overline \xi}(S_k-S_{k-1}) > B\ee^{(S_k-d_k)/3},\,S_n\in I_n(z),\,{\underline S}_n \ge 0,\, {\underline S}_{(n/2,n]}\ge a_n(z+L+1) \right)
\\
&\le&
c_{67}n^{-3/2} \E\left[  (2+L-S_{n-k+1}){\bf 1}_{\{  {\overline \xi}(S_{n-k+1}-S_{n-k}) > \tilde B \ee^{-S_{n-k}/3} ,  {\underline {-S}}_{n-k+1}\ge -L-1\}} \right]\\
&=&
c_{67}n^{-3/2} {\hat \E}\left[  (2+(L-V(w_{n-k+1}))_+){\bf 1}_{\{  \xi(w_{n-k}) > \tilde B \ee^{-V(w_{n-k})/3} ,   -V(w_j)\ge -L-1,\forall j\le n-k+1\}} \right]
\end{eqnarray*}

\noindent where $\tilde B:= B\ee^{L+1}$. Beware that we reintegrated the measures $({\bf Q_r},r\in\mathbb{R})$ in the last line. We find that, for any $k\in[3n/4,n]$,
\begin{eqnarray*}
&& \Q\left(\xi(w_k) > B \ee^{(V(w_{k})-d_k)/3} ,\,w_n\in\mathcal Z^{z,L}_n \right)\\
&\le&
c_{67}n^{-3/2} {\hat \E}\left[  (2+(L-V(w_{n-k+1}))_+){\bf 1}_{\{  \xi(w_{n-k}) > \tilde B \ee^{-V(w_{n-k})/3} ,   -V(w_j)\ge -L-1,\forall j\le n-k+1\}} \right].
\end{eqnarray*}

\noindent This is the analog of (\ref{eq:goodappendix1}), replacing there $V(w_j)$ by $-V(w_j)$, $k$ by $n-k+1$ and $\{V(w_j)\ge 0,\, \forall\,j\le k\}$ by $\{-V(w_j)\ge -L-1,\, \forall\,j\le n-k+1\}$. Then (\ref{eq:nb}) follows as in the case $k\in[1,3n/4]$. This with (\ref{eq:3n4}) prove (\ref{eq:n}) hence the lemma. \hfill $\Box$

\section{Notation}
\label{s:notation}

\noindent {\it Branching random walk: } 

$\mathcal L$ : the point process 

$X,\tilde X$ : defined in (\ref{def:X})

$\T$ : the genealogical tree 

$V(x)$ : position of particle $x$

$|x|$ : generation of vertex $x$

$\I(x)$ : siblings of vertex $x$

$x_k$ : ancestor at generation $k$ of vertex $x$

$M_n$ : minimum at generation $n$ of the {\it non-killed} branching random walk 

$\mathscr F_n$ : $\sigma-$algebra of the branching random walk up to time $n$

$\xi(x)$ : defined in (\ref{def:xi})

$ \Phi_{k,n}(r) := \P(M_{n-k}\ge (3/2) \ln(n) - r).$

$\mathcal Z[A]$: set of particles freezed when going above level $A$.

\bigskip

\noindent {\it Killed branching random walk}

$\T^{\rm kill}$ : the genealogical tree of the ${\it killed}$ branching random walk 

$|u|^{\rm kill}$ : generation of a vertex $u$ when $u\in\T^{\rm kill}$

$M_n^{\rm kill }$ : minimum at  generation $n$ of the {\it killed} branching random walk

$m^{{\rm kill},n}$ : uniform particle in $\T^{\rm kill}$ among those achieving $M_n^{\rm kill}$

 $\Phi^{\rm kill}_{k,n}(x,r):= \P_x\left( M_{n-k}^{\rm kill} \le (3/2)\ln(n) -r \right)$
\bigskip

\noindent {\it Random walk}

$(S_n)_{n\ge 0}$ : non-lattice centered random walk with finite variance, defined by (\ref{many-to-one}). The non-lattice assumption is dropped in the Appendix

$\sigma^2$ : variance of $S_1$

$R(x)$: renewal function of $S$

$R_-(x)$ : renewal function of $-S$

Many-to-one lemma : equation (\ref{many-to-one})

$H_k,T_k$ : strict descending ladder heights and epochs of $S$

$H_k^-,T_k^-$ : strict descending ladder heights and epochs of $-S$

\bigskip

\noindent {\it Martingales}

$W_n$ : additive martingale at time $n$

$D_n$ : derivative martingale at time $n$

$D_n^{(\beta)}$: martingale of the branching random walk killed below $-\beta$

\bigskip 

\noindent {\it Probability measures}

$\P_a$ : probability under which the branching random walk $(V(x))_{x\in \T}$ and the random walk $(S_n)_n$ starts at $a$ ($\P_0=\P$). Expectation $\E_a$

$\Q_a$ : tilted probability I defined by (\ref{defQ}). Expectation $\hat{\E}_a$

$\Q_a^{(\beta)}$ : tilted probability II defined by (\ref{def:Qbeta}). Expectation $\hat{\E}_a^{(\beta)}$

\bigskip

\noindent {\it  Spine decomposition I}

$w_n$ : spine at generation $n$

$(V(w_n))_n$ : centered random walk distributed as $(S_n)_n$ 

$\hat{\mathcal L}$ : Radon-Nykodim derivative $\sum_{i\in \mathcal L} \ee^{-V(i)}$ with respect to $\mathcal L$

$\hat{\mathcal B}$ : branching random walk  with a spine. Under $\P$, we identify $(V(x))_{x\in \T}$ with $\hat{\mathcal B}$ 

$\hat{\mathscr F}_n$ : $\sigma-$algebra of $\hat{\mathcal B}$ up to time $n$

$\hat{\mathcal G}_n$ : $\sigma-$algebra of the spine and its siblings up to time $n$ ($\hat{\mathcal G}_n\subset  \hat{\mathscr F}_n$)

\bigskip

\noindent {\it  Spine decomposition II }

$w_n^{(\beta)}$ : spine at generation $n$

$(V(w_n^{(\beta)}))_n$ : random walk conditioned to stay above $-\beta$ 

$\hat{\mathcal B}^{(\beta)}$ : branching random walk  with a spine

$\hat{\mathscr F}_n^{(\beta)}$ : $\sigma-$algebra of $\hat{\mathcal B}^{(\beta)}$ up to time $n$

$\hat{\mathcal G}_n^{(\beta)}$ : $\sigma-$algebra of the spine and its siblings up to time $n$ ($\hat{\mathcal G}_n^{(\beta)}\subset  \hat{\mathscr F}_n^{(\beta)}$)

\bigskip 

\noindent {\it Paths of particles}

$a_n(z)={3\over 2} \ln(n)-z$

$d_k(n,z,\lambda)$ : defined in (\ref{def:dk})

$e_k$ : defined in (\ref{def:ek})

$I_n(z)=[a_n(z)-1,a_n(z))$

$\mathcal Z^{z,L}_n$ : in Definition \ref{def:ZzL}, see Figure \ref{f:ZzL}. Particles of generation $n$ that stayed above $d_k(n,z+L,1/2)$ and end in $I_n(z)$

$\mathcal S^r$ : defined in (\ref{def:SA}), see Figure \ref{f:Sr}. Set of particles that achieve a new minimum (on their ancestral line) 

$B_n^z(u)$ : defined in (\ref{def:Bn}), see Figure \ref{f:Bnz}. Equal to $1$ if there is a line of descent from $u$ to a vertex at generation $n$ which stays above $V(u)$ and ends below $a_n(z)$

$\mathcal T^{\,r}$ : defined in (\ref{def:T}) 

$z$-good vertex : defined in (\ref{def:good})

$\mathcal E_n(z,b)$ : defined in (\ref{def:En}). Good event on which the particles at generation $n$ which are located below $a_n(z)$ have a common ancestor with the spine at generation greater than $n-b$

$F_{L,b}$ : defined in (\ref{def:FLb})

$C_{L,b}$ : defined in (\ref{def:CLb})

\bigskip

{\it Acknowledgements.} I would like to thank the referee for his helpful comments. I thank also Loic de Raphelis and Zhan Shi for spotting mistakes in an earlier version.

\bigskip


{\footnotesize

\baselineskip=12pt

\hskip110pt Department of Mathematics and Computer science

\hskip110pt Eindhoven University of Technology

\hskip110pt P.O. Box 513

\hskip110pt 5600 MB Eindhoven

\hskip110pt The Netherlands

\hskip110pt {\tt elie.aidekon@gmail.com}

}

\end{document}